\def\version{5.2.2023}
\def\users{us}  %
\def\users{final-layout}   
\documentclass[12pt]{article}
\usepackage{upgreek}

\usepackage{xcolor}
\usepackage{bm,amsmath,amsthm,hyperref,amsfonts,amssymb,color}
\usepackage{mathrsfs} 
\usepackage{cite}

\usepackage{ifthen}
\ifthenelse{\equal{\users}{final-layout}}{}{
\usepackage{fancyhdr}
\pagestyle{fancy}
\headheight=28pt\headwidth=17cm
\definecolor{gray}{gray}{0.5}
\rhead{\color{gray}Landau PT viscoelastic magnets\\
T.Roub\'\i\v cek}
\chead{}
\lhead{Version\,\version, file:\,\jobname.tex
\\
compiled:
\number\day.\number\month.\number\year\ at
\the\hour:\ifnum\minute<10 0\fi\the\minute\ h\ \ \ \ \ 
}
}

\definecolor{labelkey}{rgb}{1.,.2,0.}

\newcount\hour \newcount\minute
\hour=\time
\divide \hour by 60
\minute=\time
\loop \ifnum \minute > 59 \advance \minute by -60 \repeat

\usepackage[normalem]{ulem}

\usepackage[normalem]{ulem}
\usepackage{ifthen}
\usepackage{color}

\ifthenelse{\equal{\users}{final-layout}}{

	\newcommand{\COMMENT}[1]{}
	\newcommand{\COMMENTGT}[1]{}
	\newcommand{\TODO}[1]{}
	\newcommand{\INTERNAL}[1]{}
	\newcommand{\QUESTION}[1]{}
	\newcommand{\DELETE}[1]{}

	\newcommand{\REM}[1]{\marginpar{\bfseries\tiny{\color{blue}}}}
    \newcommand{\MARGINOTE}[1]{}
}
{
	
	\newcommand{\COMMENT}[1]{{\color{red}\uuline{#1}\color{black}}}
	\newcommand{\COMMENTGT}[1]{{\hfill\large\color{red}***{#1}***\color{black}\hfill}\\}
	\newcommand{\TODO}[1]{{\color{red}\uuline{#1}\color{black}}}
	\newcommand{\INTERNAL}[1]{\footnote{#1}}
	\newcommand{\QUESTION}[1]{{\color{brown}\uuline{#1}\color{black}}}
	\newcommand{\DELETE}[1]{{\color{red}\sout{#1}\color{black}}}

	\newcommand{\REM}[1]{\marginpar{\bfseries\tiny{\color{blue}#1}}}
\newcommand{\MARGINOTE}[1]{\marginpar{\color{red}\tiny\texttt{#1}}}
}
\newcommand\ZJ[1]{\mathchoice
                 {{\buildrel{\hspace*{.1em}{_{\,\boldsymbol\circ}}}\over{#1}}}
                 {{\buildrel{\hspace*{.1em}{_{\,\boldsymbol\circ}}}\over{#1}}}
                 {{\buildrel{\hspace*{.1em}{\boldsymbol\circ}}\over{#1}}}
                 {{\buildrel{\hspace*{.1em}{\boldsymbol\circ}}\over{#1}}}}
\newcommand\DT[1]{\mathchoice
                 {{\buildrel{\hspace*{.1em}\text{\LARGE.}}\over{#1}}}
                 {{\buildrel{\hspace*{.1em}\text{\LARGE.}}\over{#1}}}
                 {{\buildrel{\hspace*{.1em}\text{\Large.}}\over{#1}}}
                 {{\buildrel{\hspace*{.1em}\text{\large.}}\over{#1}}}}
\newcommand\pdt[1]{\frac{\partial{#1}}{\partial t}} 
\newcommand{\lineunder}[2]{\LU{\begin{array}[t]{c}\underbrace{#1}\vspace*{.5em}\end{array}}{\mbox{\footnotesize\rm #2}}}
\newcommand{\LU}[2]{\begin{array}[t]{c}#1\vspace*{-1em}\\_{#2}\end{array}}
\newcommand{\linesunder}[3]{\LSU{\begin{array}[t]{c}\underbrace{#1}\vspace*{.5em}\end{array}}{\mbox{\footnotesize\rm #2}}{\mbox{\footnotesize\rm #3}}}
\newcommand{\LSU}[3]{\begin{array}[t]{c}#1\vspace*{-1em}\\_{#2}\vspace*{-.5em}\\_{#3}\end{array}}
\newcommand{\morelinesunder}[4]{\LSUU{\begin{array}[t]{c}\underbrace{#1}\vspace*{.5em}\end{array}}{\mbox{\footnotesize\rm #2}}{\mbox{\footnotesize\rm #3}}{\mbox{\footnotesize\rm #4}}}
\newcommand{\LSUU}[4]{\begin{array}[t]{c}#1\vspace*{-1em}\\_{#2}\vspace*{-.5em}\\_{#3}\vspace*{-.5em}\\_{#4}\end{array}}
\newcommand{\Item}[2]{\parbox[t]{.055\textwidth}{#1}\hfill%
      \parbox[t]{.945\textwidth}{#2}\vspace*{.8mm}} 
\newcommand{\divS}{\mathrm{div}_{\scriptscriptstyle\textrm{\hspace*{-.1em}S}}^{}}
\newcommand{\nablaS}{\nabla_{\scriptscriptstyle\textrm{\hspace*{-.3em}S}}^{}}
\newcommand{\NablaS}{\Nabla_{\scriptscriptstyle\textrm{\hspace*{-.3em}S}}^{}}
\def\Vdots{\!\mbox{\setlength{\unitlength}{1em}
\begin{picture}(0,0)
\put(-.07,0){.}
\put(-.07,.3){.}
\put(-.07,.6){.}
\end{picture}\hspace*{.2em}}}
\usepackage{mathrsfs}   
\usepackage{eucal}      

  \def\bbI{{\mathbb I}}

\def\FG{\boldsymbol}
 \def\bb{{\FG b}}  
\def\dd{{\FG d}} \def\ee{{\FG e}} \def\ff{{\FG f}} 
 \def\hh{{\FG h}} 
\def\jj{{\FG j}}   
\def\mm{{\FG m}} \def\nn{{\FG n}}  
  \def\rr{{\FG r}} 
\def\tt{{\FG t}}  
\def\vv{{\FG v}}  \def\xx{{\FG x}}

\def\DD{{\FG D}} 
\def\FF{{\FG F}} 

 \def\KK{{\FG K}}  
\def\MM{{\FG M}}   
   
\def\SS{{\FG S}} \def\TT{{\FG T}} 

\newcommand{\R}{\mathbb R}
\newcommand{\N}{\mathbb N}
\newcommand{\Nabla}{{\nabla}}
\newcommand{\Fe}{\FF}
\newcommand{\FFeps}{\FF_{\!\EPS}}
\newcommand{\FFepsk}{\FF_{\!\EPS k}}
\newcommand{\EE}{{\bm e}}
\newcommand{\pl}{\partial}
\newcommand{\eq}[1]{(\ref{#1})}
\renewcommand{\d}{\mathrm d}  
\newcommand{\barOmega}{\hspace*{.2em}{\overline{\hspace*{-.2em}\varOmega}}}
\newcommand{\bulet}{\text{\footnotesize$\,\bullet\,$}}

\newtheorem{theorem}{Theorem}[section]
\newtheorem{lemma}[theorem]{Lemma}
\newtheorem{definition}[theorem]{Definition}
\newtheorem{example}[theorem]{Example}
\newtheorem{proposition}[theorem]{Proposition}

\newtheorem{remark}[theorem]{Remark}

\numberwithin{equation}{section}

\usepackage{graphicx}
\usepackage{psfrag} 
\newcounter{myfigure}
\newenvironment{my-picture}[3]{\refstepcounter{myfigure}\label{#3}\setlength{\unitlength}{1cm}\begin{picture}(#1,#2)}{\end{picture}}

\renewcommand{\varsigma}{c_1}

\newcommand{\UUU}{\color{blue}}
\newcommand{\EEE}{\color{black}}

\begin{document}


\allowdisplaybreaks

\def\upyy{\text{\bf{y}}}
\def\upxx{\text{\bf{X}}}
\def\upvv{\text{\bf{v}}}
\def\upFF{\text{\bf{F}}}
\def\upTT{\text{\bf{T}}}
\def\upDD{\text{\bf{D}}}
\def\uprhoxi{\uprho^{\bm\xi}}
\def\uprhoxik{\uprho^{\bm\xi_k}}

\vspace*{2em}

\begin{center}

{\Large\bf Landau theory for ferro-paramagnetic phase transition\\[.2em]
in finitely-strained viscoelastic magnets}

\bigskip\bigskip\bigskip

{\large Tom\'a\v s Roub\'\i\v cek\,}

\bigskip

{\small\it Mathematical Institute, Charles University,\\
Sokolovsk\'a 83, CZ-186~75~Praha~8,  
\\email: ${\small\texttt{tomas.roubicek@mff.cuni.cz}}$}

\medskip

{\small\it 
Institute of Thermomechanics, Czech Academy of Sciences,\\
Dolej\v skova 5, CZ-18200~Praha~8, Czech Republic}

\bigskip\bigskip\bigskip

{\small\bf Abstract}.

\vspace*{.3em}

\begin{minipage}[t]{37em}{\small
\baselineskip=13pt
The thermodynamic model of visco-elastic deformable magnetic materials
at finite strains is formulated in a fully Eulerian way in rates.
The Landau theory applies for ferro-to-para-magnetic phase transition,
the gradient theory (leading exchange energy) for magnetization with
general mechanically dependent coefficient, hysteresis in magnetization
evolution by Landau-Lifshitz-Gilbert equation involving objective corotational
time derivative of magnetization, and demagnetizing field are considered in
the model. The Kelvin-Voigt viscoelastic rheology with a higher-order viscosity
(exploiting the concept of multipolar materials) is used,
allowing for physically relevant frame-indifferent stored energies and
for local invertibility of deformation. The model complies with
energy conservation and Clausius-Duhem entropy inequality. 
Existence and a certain regularity of weak solutions is proved
by a Faedo-Galerkin semi-discretization and a suitable regularization.
\medskip

\medskip


}
\end{minipage}
\end{center}

\bigskip

\baselineskip=16pt

\def\TRACTION{\bm{k}}
\def\GRAVITY{\bm{g}}
\def\rhoR{\uprho}
\def\M{m}
\def\MM{M}
\def\COUPLING{\upzeta}
\def\LAM{\lambda}
\def\W{w}
\def\OMEGA{\omega}
\def\EPS{\varepsilon}
\def\COND{\mbox{\footnotesize$\mathcal{K}$}}
\def\wt{\widetilde}
\def\wh{\widehat}

\def\ALPH{\alpha}
\def\ONEALPH{2}
\def\DIS{\nu_1}
\def\TWO{2}
\def\TWOprime{2}
\def\EXP{\mu}
\def\NU{\nu}
\def\GM{G_\text{\sc m}}

\def\Ee{{\bm E}_{\rm e}}
\def\Ep{{\bm E}_{\rm p}}
\def\HC{h_\text{\sc c}^{}}
\def\MEXP{m}

\section{Introduction -- deforming magnetic continua}

The magnetic materials which are not completely rigid represent interesting,
important, and difficult multi-physical concatenation of mere (thermo)continuum
mechanics and mere micromagnetism. Beside homogeneous visco-elastic magnets,
it may concern elastically rather soft materials filled with magnetic
particles, e.g.\ rocks (which can be considered soft on long time scales)
and polymers (i.e.\ so-called magneto-rheological elastomers or ferrogels),
which however needs to involve creep which is not consider in this paper
rather for not making the model too complicated. We will focus ourselves to 
general finite (also called large) strain mechanics in the Eulerian formulation.

This magneto-mechanical subject has been addressed in \cite[Ch.6]{DorOgd14NTEM}
or anisothermal but not with explicitly articulated
equations \cite{Brow66MI} and also, in a thermodynamic context,
\cite[Ch.6]{Maug88CMES}. Even in the purely mechanical isothermal cases,
and a-fortiori in anisothermal situations,
the visco-elastodynamics at finite strains has been articulated in
\cite{Ball02SOPE,Ball10PPNE}
as a difficult open problem as far as existence of weak solutions concerns. 
There is a certain agreement that, for analytical
reasons, a certain enough strong dissipation mechanism is to be
involved to make the dynamical problem parabolic, although some
hyperbolic models exist, as mentioned below. The simplest
variant is the {\it Kelvin-Voigt viscoelastic rheology}.

The mentioned {\it  Eulerian approach} is standardly believed
to be well fitted with fluids. It is particularly suitable in situations
when there is no natural reference configuration or where a reference
configuration becomes less and less relevant during long-time evolution,
which may however apply also for solids. A formulation of equations in
current deforming configuration needs rather
velocity/strain than displacement to be involved in the momentum
equation. The advantage is an easier possibility to involve
interaction with outer spatial fields (here magnetic and gravity) and
avoiding the pull-back and push-forward manipulation. On the other
hand, there is a necessity to involve convective derivative and transport
equations and also evolving the shape of the body is troublesome.
In isothermal situations, such model was formulated and analyzed as
incompressible in \cite{LeLiZh08GSIV,LiuWal01EDFC} and as compressible in
\cite{HuMas16GSRH,QiaZha10GWPC}. The mentioned higher gradients that
would allow for reasonable analysis can now be involved rather in the
dissipative than conservative part, so that their influence
manifests only in fast evolutions. In the isothermal situations
it was used in quasistatic case in \cite{Roub??VELS} and
in dynamical case in \cite{RouSte22VESS} when considering
the stored energy in the actual configuration, which then gives
an energy pressure in the stress tensor. In anisothermal situations,
such free-energy pressure would be directly added into stress tensor in an
non-integrable way and likely would cause technical
difficulties.

The main attributes of the devised model are:\\
\Item{{$\mbox{$\bm\oplus$}$}}{Concept of {\it hyperelastic materials} (whose conservative-stress
response comes from a free energy) combined with the {\it Kelvin-Voigt viscoelastic
rheology} and also evolution of magnetization is driven by this free energy.}
\Item{{$\mbox{$\bm\oplus$}$}}{Inertial effects in fully compressible context (in particular
with varying mass density) are considered.$_{_{}}$}
\noindent\Item{{$\mbox{$\bm\oplus$}$}}{The rate formulation in terms of velocity and deformation
gradient is used while the deformation itself does not explicitly occur.}
\Item{{$\mbox{$\bm\oplus$}$}}{Magnetic phenomena covered by the model includes:
ferro-to-para magnetic phase transition, hysteresis due to the 
pinning effects, exchange energy depending on deformation gradient (and
in particular on compression/expansion), and demagnetizing field.}
\Item{{$\mbox{$\bm\oplus$}$}}{Mechanical consistency in the sense that
{\it frame indifference} of the free energy (which is in particular
{\it nonconvex} in terms of deformation gradient and in magnetization) and 
its {\it singularity} under infinite compression in relation with
{\it local non-interpenetration} as well as objective corotational
time derivative for magnetization transport.}

\newpage

\noindent\Item{{$\mbox{$\bm\oplus$}$}}{Thermodynamic consistency of the
thermally coupled system in the sense that the {\it total energy is conserved}
in a closed system, the {\it Clausius-Duhem entropy inequality}
holds, and temperature stays non-negative.}
\noindent\Item{{$\mbox{$\bm\oplus$}$}}{The nonconservative part of the stress in the Kelvin-Voigt model
containing a higher-order component reflecting the concept of nonsimple
{\it multipolar media} is exploited.}
\Item{{$\mbox{$\bm\oplus$}$}}{The model allows for rigorous mathematical analysis as far as
existence and certain regularity of energy-conserving weak solutions
concerns.}

\vspace*{.1em}

\noindent
On the other hand, some simplifications are adopted:\\
\Item{{$\mbox{$\bm\ominus$}$}}{Relatively slow evolution is implicitly
assumed, which allows for
reducing the full Maxwell electromagnetodynamics to magneto-statics.}
\Item{{$\mbox{$\bm\ominus$}$}}{Electric conductivity (and in particular eddy currents) is not
considered.}

\vspace*{.3em}

\noindent
As far as the non-negativity of temperature, below we will be able
to prove only that at least some solutions enjoy this attribute, although
there is an intuitive belief that all possible solutions will make it and
a hope that more advanced analytical techniques would rigorously prove it.

The main notation used in this paper is summarized in the following table:

\begin{center}
\fbox{
\begin{minipage}[t]{16.5em}\small\smallskip

$\vv$ velocity (in m/s),

$\varrho$ mass density (in kg/m$^3$),

$\rhoR$ referential mass density, 

$\Fe$ deformation gradient,

$\mm$ magnetization (in A/m),

$\theta$ temperature (in K),

$\TT$ Cauchy stress (symmetric, in Pa),

$\KK$ Korteweg stress (symmetric, in Pa),

$\SS$ skew-symmetric stress (in Pa),

\vspace*{.0em}

$\DD$ dissipative stress (in Pa),

\vspace*{.0em}

$\mathscr{H}$ elastic hyperstress (in Pa\,m),

\vspace*{.0em}

$\jj$ heat flux (in W/m$^2$),

\vspace*{.0em}

\vspace*{.0em}

$\TRACTION$ traction load,

\vspace*{.0em}

$\hh$ (total) magnetic field (in A/m),

\vspace*{.0em}

$\hh_{\rm ext}$ external magnetic field (in A/m),

\vspace*{.0em}

$u$ demagnetizing-field potential (in A),

\vspace*{.0em}

$\bm\xi$ the return mapping,

\vspace*{.0em}

${\rm Cof}(\cdot)$ cofactor matrix,

\vspace*{.0em}

$\nu_\flat>0$ a boundary viscosity,

\vspace*{.0em}

$\det(\cdot)$ determinant of a matrix,

\vspace*{.0em}

$\R_{\rm sym}^{d\times d}=\{A\in\R^{d\times d};\ A^\top=A\}$,

\vspace*{.0em}

$\bbI\in\R_{\rm sym}^{d\times d}$ the unit matrix,

\smallskip \end{minipage}
\begin{minipage}[t]{22.8em}\small\smallskip

$\uppsi=\uppsi(\FF,\mm,\theta)$ referential free energy (in J/m$^3$=Pa),

$\upvarphi=\upvarphi(\Fe,\mm)$ referential stored energy (in J/m$^3$=Pa),

$\COUPLING=\COUPLING(\Fe,\mm,\theta)$ referential heat part of free energy,

$\ee(\vv)=\frac12\Nabla\vv^\top\!+\frac12\Nabla\vv$ small strain rate (in s$^{-1}$),


$\HC=\HC(\FF,\theta)$ magnetic coercive force (in A/m),

$\W$ heat part of internal energy (enthalpy, in J/m$^3$),

$(^{_{_{\bullet}}})\!\DT{^{}}=\pdt{}{^{_{_{\bullet}}}}\!+(\vv{\cdot}\Nabla)^{_{_{\bullet}}}$
convective time derivative,

$(^{_{_{\bullet}}})\!\ZJ{^{}}=(^{_{_{\bullet}}})\!\DT{^{}}-{\rm skw}(\nabla\vv)^{_{_{\bullet}}}$
corotational time derivative, 

${\rm skw}\,A=$ a skew-symmetric part, i.e. $\frac12A-\frac12A^\top$,

$\mathscr{S}$ skew-symmetric magnetic 
hyperstress (in Pa\,m),

$\cdot$ or $:$ scalar products of vectors or matrices, 

$\Vdots\ \ $ scalar products of 3rd-order tensors,

$\COND=\COND(\Fe,\theta)$ thermal conductivity (in W/m$^{-2}$K$^{-1}$),

$\kappa=\kappa(\Fe)$ exchange-energy coefficient (in kg${\cdot}$m/C$^2$),

$c=c(\Fe,\theta)$ heat capacity (in Pa/K),

$\tt$ magnetic ``driving force'' (in Jm$^{-2}$A$^{-1}$),

$\gamma=\gamma(\FF,\mm,\theta)$ gyromagnetic ratio (in C/kg),

$\tau$ ``viscous'' magnetic damping coefficient (in s),


$\DIS,\nu_2>0$ a bulk (hyper)viscosity coefficients,

$\hh_{\rm dem}{=}{-}\nabla u$ demagnetizing field (in A/m),

$\GRAVITY$ external bulk load (gravity acceleration in m/s$^{2}$),

$\mu_0$ vacuum permeability ($\sim1.257{\times}10^{-6}$\,H/m).
\smallskip \end{minipage}
}\end{center}

\vspace{-.5em}

\begin{center}
{\small\sl Table\,1.\ }
{\small
\,Summary of the basic notation used through Sections~\ref{sec-model}
and \ref{sec-anal}.\newline\vspace*{-.8em}\vskip -.8em
Convention: uprighted\,=\,referential, slanted\,=\,actual.\hspace*{4em}
}
\end{center}

\medskip

In comparison with \cite{Roub??VELS,RouSte22VESS},
the novelty of this paper is to apply the Eulerian approach
to solids in anisothermal situations, using the free energy
in a reference configuration,  which does not see the
energy-pressure in the stress tensor and which is also more fitted with
usually available experimental data. The analysis combines
$L^1$-theory for the heat equation adapted to the convective
time derivatives and the techniques from compressible fluid dynamics
adapted for solids.

For completeness, let us still mention a competitive, Lagrangian thermodynamic
formulation (including also diffusion) \cite{RouTom18TCMM} formulating
the equations in a certain fixed ``reference'' configuration. This approach
allows easily for deformation of the shape of the body and easier treatment
inertial forces but a frame-indifferent viscosity and interaction
with spatial gravity and magnetic forces is much more complicated.

The plan is as follows: formulation of the model in the actual Eulerian
configuration and its energetics and thermodynamics is presented in
Section~\ref{sec-model}, recalling first the micromagnetism and
Landau transition in rigid magnets in Sect.~\ref{sec-micro} and
finite-strain kinematics of deformable continua in Sect.~\ref{sec-mech}
before formulating the model in Section~\ref{sec-magnet} and showing
its energetics in Section~\ref{sec-engr}. Then, in Section~\ref{sec-anal},
the rigorous analysis by a suitable regularization and a (semi) Faedo-Galerkin
approximation is performed, combined with theory of transport by regular
velocity fields.

\section{The thermodynamic model and its energetics}\label{sec-model}

It is important to distinguish carefully the referential and
the actual time-evolving coordinates. Our aim is to formulate
the model eventually in actual configurations, i.e.\ the Eulerian
formulation, reflecting also the reality in many (or even most)
situations (and a certain general agreement) that a reference
configuration is only an artificial construction and, even
if relevant in some situations, becomes successively more and
more irrelevant during evolution at truly finite strains.
Typical materials involve magnetic gels or elastomers 
or magnetic rocks which are viscoelastic on geological timescales.
On the other hand, some experimental material data 
are related to some reference configuration -- typically it
concerns mass density and stored of free energies per mass (in J/kg)
or per referential volume (in J/m$^3$=Pa) as considered here.

We will present briefly the fundamental concepts and formulas which can
mostly be found in the monographs, as e.g.\ \cite[Part~XI]{GuFrAn10MTC}
or \cite[Sect.~7.2]{Mart19PCM}.

\subsection{Micromagnetism and ferro-parramagnetic transition}\label{sec-micro}
Let us briefly recall the micromagnetic model in rigid magnets and Landau's
phase-transition theory \cite{Land37TPT}, cf.\ also
\cite[Sec.39]{LanLif84ECM} or the monographs \cite{Bert98HM,Tremolet}. 
The basic ingredient governing static (and later also evolution) model is the free energy
$\uppsi=\uppsi(\mm,\theta)$ depending on magnetization $\mm$ and temperature $\theta$.

In the micromagnetism, the free energy $\uppsi$ is augmented by
the exchange energy $\frac{\kappa}2|\nabla\mm|^2$ with $\kappa$ a coefficient
determining an internal length-scale, responsible for a typical fine domain structure
in ferromagnets.

The magnetization itself induces a magnetic field, called a self-induced
{\it demagnetizing field} $\hh_{\rm dem}$. For many (or maybe most) applications,
full Maxwell electro-magnetic system is considered simplified to
{\it magnetostatics}, considering slow evolution and neglecting
in particular eddy currents and even confining on 
electrically nonconductive media. The Maxwell system then reduces to 
the Amp\`ere law ${\rm curl}\hh_{\rm dem}=\bm0$ and the 
Gauss law ${\rm div}\,\bb=0$ for the magnetic induction with is
given by $\bb=\mu_0\hh_{\rm dem}+\mu_0\mm$ where $\mu_0$ is the
physical constant (vacuum permeability). The Amp\`ere law
ensures existence of a scalar-valued potential $u$ such that
$\hh_{\rm dem}=-\nabla u$. These equations are considered on the whole
Universe $\R^d$ while, of course, the magnetization $\mm$ is only in the
body $\varOmega$ while outside it is considered zero, which is
articulated by introducing the characteristic function
$\chi_{\varOmega}^{}$ defined as $\chi_{\varOmega}^{}(\xx)=1$ if
$\xx\in\varOmega$ and $\chi_{\varOmega}^{}(\xx)=0$ if
$\xx\in\R^d{\setminus}\varOmega$. By substitution, we obtain the equation 
\begin{align}\label{u-eq}
{\rm div}(\nabla u-\chi_{\varOmega}^{}\mm)=0\ \ \ \text{ in }\ \R^d
\end{align}
to be considered in the sense of distributions. Under an external magnetic
field $\hh_{\rm ext}$, the overall effective
magnetic field $\hh$ is
\begin{align}
\hh=\hh_{\rm ext}-\hh_{\rm dem}=\hh_{\rm ext}+\nabla u\,.
\label{effective-field}\end{align}
Although not directly relevant in this paper, let us anyhow
remind that, for a fixed temperature $\theta$, the standard ferro-magnetostatic
theory is based on the free energy
$\uppsi(\mm,\theta,\nabla\mm)=\wt\uppsi(\mm,\theta)+\frac{\kappa}2|\nabla\mm|^2$
leading to the overall energy
\begin{align}
(\mm,u)&\mapsto\int_\varOmega\hspace*{-.5em}
\linesunder{\wt\uppsi(\mm,\theta)_{_{_{_{_{}}}}}}{free}{energy}\hspace*{-.7em}
+\hspace*{-.7em}\linesunder{\frac{\kappa}2|\nabla\mm|^2}{exchange}{energy}\hspace*{-.7em}
-\hspace*{-.7em}\linesunder{\mu_{0_{_{_{_{}}}}}(\hh_{\rm ext}{+}\nabla u){\cdot}\mm}{energy of $\mm$ in the}{magnetic field $\hh$}\hspace*{-.5em}\d\xx
-\int_{\R^d}\hspace*{-1.5em}\linesunder{\frac{\mu_0}2|\nabla u|^2}{energy of demag-}{netizing field}\hspace*{-1.5em}\d\xx\,.
\label{magnetostatic-functional}\end{align}
Notably, this functional is concave with respect to $u$ and has a saddle-point
character. The static configurations $(\mm,u)$ are standardly considered as
minimizing with respect to $\mm$ and maximizing with respect to $u$, i.e.\
a critical point or \eq{magnetostatic-functional}. The 1st-order
optimality conditions then gives the system
\begin{align}\label{magnetostatic-system}
\lineunder{\wt\uppsi_\mm'(\mm,\theta)-{\rm div}(\kappa\nabla\mm)}{$=\tt$ magnetic ``driving force''}\hspace*{-1em}=\mu_0\hh\ \ \ \text{ and }\ \ \ \eq{u-eq}\,.
\end{align}
The mentioned saddle-point character can be eliminated by executing maximization
with respect to $u$, i.e.\ in fact the partial Legendre transform. This gives,
when testing \eq{u-eq} by $\mu_0 u$, which gives
$\int_{\R^d}\mu_0|\nabla u|^2\,\d\xx=\int_\varOmega\mu_0\mm{\cdot}\nabla u\,\d\xx
=-\int_\varOmega\mu_0\mm{\cdot}\hh_{\rm dem}\,\d\xx$. Substituting it into
\eq{magnetostatic-functional}, the functional depending on $\mm$ which
should be minimized by static configurations is:
\begin{align}
\mm&\mapsto\int_\varOmega\hspace*{-.5em}
\linesunder{\wt\uppsi(\mm,\theta)_{_{_{_{_{}}}}}}{free}{energy}\hspace*{-.7em}
+\hspace*{-.7em}\linesunder{\frac{\kappa}2|\nabla\mm|^2}{exchange}{energy}
\hspace*{-.7em}-\hspace*{-.7em}\linesunder{\mu_{0_{_{_{_{}}}}}\hh_{\rm ext}{\cdot}\mm}{Zeeman}{energy}\hspace*{-.5em}\d\xx
+\int_{\R^d}\hspace*{-1.5em}\linesunder{\frac{\mu_0}2|\nabla u_\mm^{}|^2}{energy of demag-}{netizing field}\hspace*{-1em}\d\xx\,.
\label{magnetostatic-functional+}\end{align}
In the rest of this paper, we will couple it with mechanical effects and
a full thermodynamics, so that the minimization of energy will no longer
be relevant. 

In case of time-varying $\hh_{\rm ext}$, a dynamics of $\mm$ governed by the 
{\it Landau-Lifschitz-Gilbert equation} 
$\gamma^{-1}\pdt{}\mm=\mm{\times}\hh_{\rm eff}$ with $\hh_{\rm eff}=\hh-\tt-\tt_{\rm dis}$
an effective field composed from a conservative part $\tt$ arising from 
a free energy \eq{magnetostatic-functional}, cf.\ \eq{magnetostatic-system}
or also \eq{driving-field} below, while $\tt_{\rm dis}$ is a magnetic field
counting a  dissipative-processes phenomenology, and $\hh$ is from
\eq{effective-field}. Equivalently \cite{BePGVa01DDFG}, one can write it in
the Gilbert form $\gamma^{-1}\mm{\times}\pdt{}\mm=\hh_{\rm eff}$.
The basic choice of $\tt_{\rm dis}$ is the magnetic ``viscosity''
$\tau\pdt{}\mm$ with $\tau$ a phenomenological magnetic
damping coefficient. To cover the (temperature dependent) {\it hysteresis}
effects due to so-called pinning mechanism, we augment it by the
dry-friction term $\HC(\theta){\rm Dir}(\pdt{}\mm)$ where ``Dir'' denotes the
set-valued monotone ``direction'' mapping
\begin{align}
{\rm Dir}(\rr)=\begin{cases}\{r\in\R^d;\ |r|\le1\}&\text{if }\rr={\bm0}\,,\\
\quad\rr/|\rr|&\text{if }\rr\ne{\bm0}\,,\end{cases}
\end{align}
cf.\ \cite{RoToZa09GEDF} and Remark~\ref{rem-dry-friction} below. Note that
$\rr{\cdot}{\rm Dir}(\rr)=|\rr|$. Here, having in mind an isotropic
situation, $|\cdot|$ denotes the Euclidean norm, but in principle some
other anisotropic norms on $\R^d$ can be considered, too.
Altogether, we consider the specific Gilbert equation as
\begin{align}\label{LL-friction-}
\tau\pdt\mm+\HC(\theta){\rm Dir}\Big(\pdt\mm\Big)
-\frac{\mm}{\gamma(\theta)}{\times}\pdt\mm=\mu_0\hh-\tt\,.
\end{align}
The coercive force $\HC=\HC(\theta)$ determines the width of hysteresis
loops within slowly time-varying oscillatory external field $\hh_{\rm ext}$.
The gyromagnetic term should disappear under high temperatures, i.e.\ 
$1/\gamma(\cdot)$ going to 0 for temperatures around or above Curie
temperature, as articulated in \cite{Maug76NMHT}. Let us note that
\eq{LL-friction-} balances the terms in the physical units A/m, as standard.

\begin{example}[{\sl Ferro-to-para-magnetic transition}]\upshape
A simplest example of free energy in rigid isotropic magnetic materials is
\begin{align}\label{exa-phi}
\wt\uppsi(\mm,\theta)=
a_0(\theta{-}\theta_\text{\sc c})|\mm|^2+b_0|\mm|^4
+c_0\theta(1{-}{\rm ln}\theta)\,.
\end{align}
In static magnetically soft ferromagnetism, the magnetization minimizes the
energy. Here the minimum of $\wt\uppsi(\bulet,\theta)$ is attained on the orbit
$|\mm|=m_\text{\sc s}(\theta)$ with $m_\text{\sc s}(\theta)
=\sqrt{a_0(\theta_\text{\sc c}{-}\theta)/(2b_0)}$
if $0\le\theta\le\theta_\text{\sc c}$ \COMMENT{CHECK} and at $\mm=0$ if
$\theta\ge\theta_\text{\sc c}$, cf.\ the solid line in Figure~\ref{fig1}.
Under an applied magnetic field $\hh_{\rm ext}$, the minimum of
$\mm\mapsto\wt\uppsi(\mm,\theta)-\hh_{\rm ext}{\cdot}\mm$ is at some magnetization
whose magnitude is slightly bigger than $m_\text{\sc s}(\theta)$,
cf.\ the dashed line in Figure~\ref{fig1}.
This ansatz can be used for a ferro-para-magnetic transition
for a mechanically rigid magnets as formulated (and analyzed)
in \cite{PGRoTo10TCTF}. This may be quite equally interpreted
as ferri-antiferro-magnetic transition, too, cf.\ \cite{HBHC97PFTS}.
\begin{center}
\begin{my-picture}{12}{5}{fig1}
\psfrag{1}{\small$\theta_\text{\sc c}$}
\psfrag{t1}{\small$\theta$}
\psfrag{ferromagnetic}{\small ferromagnetic}
\psfrag{paramagnetic}{\small paramagnetic}
\psfrag{t2}{\small$m_\text{\sc s}$}
\psfrag{t3}{\small$\theta_\text{\sc c}$}
\psfrag{t4}{\small$0$}
\psfrag{t5}{{\small\begin{minipage}[t]{8em}\baselineskip=8pt\small
Curie temperature
\end{minipage}}}
\psfrag{t6}{\small$m_\text{\sc s}(0)$}
\includegraphics[width=25em]{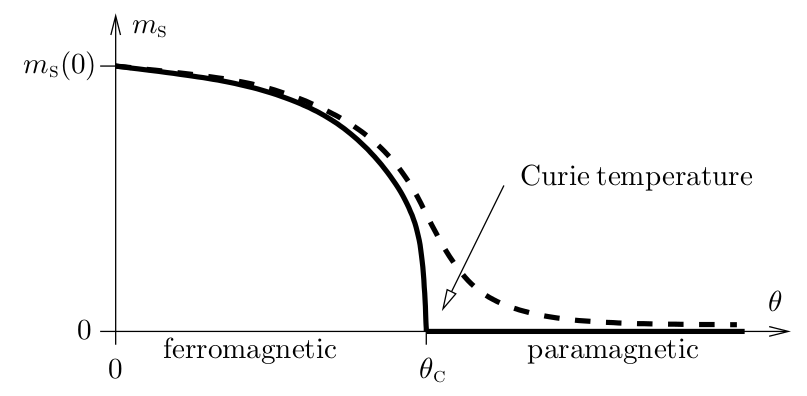}
\end{my-picture}
\nopagebreak
\\
{\small\sl\hspace*{-5em}Fig.~\ref{fig1}:~\begin{minipage}[t]{25em}
Typical dependence of saturation magnetization $m_\text{\sc s}$
on absolute temperature under zero applied field 
$\hh$ (solid line) and under some applied field (dashed line),
cf.~e.g.~\cite{Bert98HM}.
\end{minipage}
}
\end{center}
\end{example}

\subsection{Finite-strain kinematics and mass and momentum transport}\label{sec-mech}
In finite-strain continuum mechanics, the basic geometrical concept is the
time-evolving deformation $\upyy:\Omega\to\R^d$ as a mapping from a reference
configuration of the body $\Omega\subset\R^d$ into a physical space $\R^d$.
The ``Lagrangian'' space variable in the reference configuration will be
denoted as $\upxx\in\Omega$ while in the ``Eulerian'' physical-space
variable by $\xx\in\R^d$. The basic kinematic and geometrical objects are
the Lagrangian velocity $\upvv=\pdt{}\upyy$ and the Lagrangian
deformation gradient $\upFF=\Nabla_{\!\upxx}^{}\upyy$.

We will be interested in deformations $\xx=\upyy(t,\upxx)$ evolving in time,
which are sometimes called ``motions''. Further, assuming for a moment that
$\upyy(t,\cdot)$ is invertible, we define
the so-called {\it return} (sometimes called also a {\it reference})
{\it mapping} $\bm\xi:\xx\mapsto\upyy^{-1}(t,\upxx)$.
The important quantities are the
Eulerian velocity $\vv(t,\xx)=\upvv(t,\bm\xi(t,\xx))$ and the Eulerian
deformation gradient $\FF(t,\xx)=\upFF(t,\bm\xi(t,\xx))$.

Here and thorough the whole article, having the Eulerian velocity at disposal,
we use the dot-notation
$(\cdot)\!\DT{^{}}=\pdt{}+\vv{\cdot}\nabla_\xx$ for the {\it convective time
derivative} applied to scalars or, component-wise, to vectors or tensors.
Then the velocity gradient
$\Nabla\vv=\nabla_{\!\upxx}^{}\vv\nabla_{\!\xx}^{}\upxx=\DT\FF\FF^{-1}$,
where we used the chain-rule calculus  and
$\FF^{-1}=(\nabla_{\!\upxx}^{}\xx)^{-1}=\nabla_{\!\xx}^{}\upxx$. 
This gives the {\it transport equation-and-evolution  for the
deformation gradient} as
\begin{align}
\DT\FF=(\nabla\vv)\FF\,.
  \label{ultimate}\end{align}
From this, we also obtain the evolution-and-transport equation for Jacobian
$\det\FF$ as $\DT{\overline{\det\FF}}=(\det\FF){\rm div}\,\vv$ and its inverse as
\begin{align}\nonumber\\[-2.9em]
\DT{\overline{\!\!\!\bigg(\frac1{\det\FF}\bigg)\!\!\!}}\ =-\frac{{\rm div}\,\vv}{\det\FF}\,.
\label{DT-det-1}\end{align}
The return mapping $\bm\xi$ satisfies the transport equation
\begin{align}
\DT{\bm\xi}=\bm0\,;
\label{transport-xi}\end{align}
note that, since we confined on a spatially homogeneous material, 
actually $\bm\xi$ does not explicitly occur in the formulation of
the problem.

As $\FF$ depends on $\xx$, \eqref{ultimate}--\eqref{transport-xi} are equalities
which hold for a.a.\ $\xx$.
The same holds for \eqref{cont-eq+}--\eqref{LL-friction} below. Here we will benefit from the
boundary condition $\vv{\cdot}\nn=0$ below, which causes that the shape of
the actual domain $\varOmega$ does not evolve in time, i.e.\ $\varOmega=\Omega$.
The same convention concerns temperature $\theta$ and thus also $\TT$, $\eta$,
and $\DIS$ in \eq{stress-entropy} and \eq{Cauchy-dissip} below, which will
make the problem indeed fully Eulerian. Cf.\ the continuum-mechanics textbooks
as e.g.\ \cite{GuFrAn10MTC,Mart19PCM}.

The mass density (in kg/m$^3$) is an extensive variable, and its transport
(expressing that the conservation of mass) writes as the {\it continuity
equation} $\pdt{}\varrho+{\rm div}(\varrho\vv)=0$,
or, equivalently, the {\it mass evolution-and-transport equation}
\begin{align}
\DT\varrho=-\varrho\,{\rm div}\,\vv\,.
\label{cont-eq+}\end{align}
Alternatively to \eq{cont-eq+}, we will also use an
evolution-and-transport equation for the ``mass sparsity'' as the inverse
mass density $1/\varrho$:
\begin{align}
\DT{\overline{1/\varrho}}=(1/\varrho)\,{\rm div}\,\vv\,.
\label{cont-eq-inverse}\end{align}

The flow rule for the magnetization \eq{LL-friction-} is now to be
considered in deforming medium, and then the partial time derivative in
\eq{LL-friction-} should be replaced by an objective time derivative. Here
we use the Zaremba-Jaumann (corotational) time derivative $\ZJ\mm$, defined as
\begin{align}
\ZJ\mm=\DT\mm-{\rm skw}(\nabla\vv)\mm=\pdt\mm+(\vv{\cdot}\nabla)\mm-{\rm skw}(\nabla\vv)\mm\,,
\label{def-ZJ}\end{align}
where $\DT\mm=\pdt{}\mm+(\vv{\cdot}\nabla)\mm$ denotes the convective
derivative of $\mm$. 
Moreover, in deforming continuum, we can (and should) consider a more general
$\gamma=\gamma(\FF,\mm,\theta)$ and  $\HC=\HC(\FF,\theta)$.
Thus \eq{LL-friction-} turns into
\begin{align}\label{LL-friction}
\tau\ZJ\mm+\HC(\FF,\theta){\rm Dir}(\ZJ\mm)
-\frac{\mm{\times}\ZJ\mm}{\gamma(\FF,\mm,\theta)}=\mu_0\hh-\tt\,.
\end{align}
The convective derivative itself is not objective and would not be suitable
in our context, except perhaps some laminar-like deformation as implicitly
used in an incompressible isothermal variant
in \cite{BFLS18EWSE,KaKoSc21MAWS,SchZab18USMM,Zhao18LWPB} or in a
nanoparticle transport in fluids \cite{GruWei21FITM};
for usage of $\ZJ\mm$ in \eq{LL-friction} see Remark~\ref{rem-ZJ} below.

\subsection{Magneto-viscoelasticity and its thermodynamics}\label{sec-magnet}

The main ingredients of the model are the (volumetric) {\it free energy}
$\uppsi$ and the {\it dissipative stress}. The Helmholtz free energy
$\uppsi=\uppsi(\FF,\mm,\nabla\mm,\theta)$
is considered per the {\it referential volume}, while the free energy
per actual deformed volume is $\uppsi(\FF,\mm,\nabla\mm,\theta)/\!\det\FF$.
Considering the free energy per unit reference volume is more standard in
continuum physics \cite{GuFrAn10MTC,Mart19PCM} than the free energy per
actual evolving volume and well corresponds to experimentally available data.
Here also the anisotropy (which is typical in ferromagnets on microscopical
scale) in the stored energy needs rather large strains with referential stored
energy. This last benefit is related to the fact that the referential free energy
does not give an energy pressure contribution to the Cauchy stress (cf.\ the last
term in \eq{referential-stress} below or\ \cite[Rem.\,2]{Roub??VELS}) and allows
for more easy decoupling estimation strategy decoupling the magneto-mechanical
part and the thermal part of the coupled system.

We will select out the temperature independent stored energy $\upvarphi$ and
consider the split: 
\begin{align}
\uppsi(\Fe,\mm,\nabla\mm,\theta)=\upvarphi(\Fe,\mm)
+\COUPLING(\Fe,\mm,\theta)+\frac{\upkappa(\FF)}2|\nabla\mm|^2
\ \ \ \text{ with }\ \ \COUPLING(\Fe,\mm,0)=0\,.
\label{ansatz-}\end{align}
The free energy considered per actual (not referential) volume extended
by the Zeeman energy arising by an applied external actual (not referential)
magnetic field $\hh_{\rm ext}$, i.e.\ the Gibbs-type {\it actual free energy} is thus
\begin{align}
\psi_{\rm G}^{}(t;\Fe,\mm,\nabla\mm,\theta)=
\!\!\!\linesunder{\frac{\upvarphi(\Fe,\mm)}{\det\FF}}{stored}{energy}
\!\!\!\!+\!\!\!\!
\linesunder{\frac{\COUPLING(\Fe,\mm,\theta)}{\det\FF}}{coupling and}{heat energy}
\!\!\!\!-\!\!\!\!
\linesunder{\mu_0\hh_{\rm ext}(t){\cdot}\mm_{_{_{_{}}}}}{Zeeman}{energy}
\!\!\!\!+\!\!\!\!\linesunder{\frac{\upkappa(\FF)|\nabla\mm|^2}{2
\det\FF}}{exchange}{energy}\!\!\!.
\label{ansatz}\end{align}
Thus the stored energy contain, beside the elastic stored energy, also the
so-called anisotropy magnetic energy, which may distinguish directions of easy
magnetization on the microscopical single-crystal level.
Let us note that this rather general elasto-magnetic ansatz allows for modelling
the {\it magnetic shape-memory materials}
and {\it magnetostrictive effects}.

A particularly simple case 
$\upkappa(\FF)=k\det\FF$ with some constant $k$ 
would give the actual exchange energy $k|\nabla\mm|^2/2$ independent of
$\FF$ and the actual exchange coefficient as an intensive variable. Yet, as
there does not seem any physical reason for exchange coefficient to be an
intensive variable, we want to consider a general situation covering in
particular also the ``referential'' case $\upkappa(\FF)$ constant, which
corresponds the {\it actual exchange coefficient} $\upkappa(\FF)/\det\FF$
as an extensive variable.

From the free energy \eq{ansatz-}, we can read as partial (functional)
derivatives of $\uppsi$ with respect to $\FF$, $\nabla\mm$, $\mm$, and $\theta$
respectively the {\it conservative part of the actual Cauchy stress} $\TT$,
a {\it capillarity}-like {\it stress} $\KK$, the actual conservative
{\it magnetic driving force} $\tt$, and the {\it actual entropy} $\eta$ as:
\begin{subequations}\label{driving-fields}\begin{align}\label{driving-stress}
&\TT=\frac{\uppsi_\FF'(\FF,\mm,\nabla\mm,\theta)\FF^\top\!\!\!}{\det\FF}
=
\bigg(\!\upvarphi_\Fe'(\Fe,\mm)+\COUPLING_\Fe'(\Fe,\mm,\theta)
+\frac{\upkappa'(\FF)|\nabla\mm|^2\!}2\;\bigg)\frac{\Fe^\top}{\det\Fe}\,,
\\\label{capillarity-stress}
&\KK=-\frac{(\nabla\mm)^\top\uppsi_{\nabla\mm}'(\FF,\nabla\mm)}{\det\FF}=
-(\nabla\mm)^\top\frac{\mu_0\upkappa(\FF)\nabla\mm}{\det\FF}=
-\mu_0\upkappa(\FF)\frac{\nabla\mm{\otimes}\nabla\mm}{\det\FF}\,,
\\\nonumber
&\tt=\frac{\uppsi_\mm'(\FF,\mm,\theta)}{\det\FF}-{\rm div}\frac{\uppsi_{\nabla\mm}'(\FF,\nabla\mm)}{\det\FF}
\\&\ \ 
=\frac{\upvarphi_\mm'(\Fe,\mm)+\COUPLING_\mm'(\Fe,\mm,\theta)\!\!}{\det\Fe}
-{\rm div}\Big(\frac{\upkappa(\FF)\nabla\mm}{\det\Fe}\Big)
\,,\ \ \ \ \text{ and }\ \ \ \ 
\label{driving-field}\\
&\eta=-\frac{\uppsi_\theta'(\FF,\mm,\theta)}{\det\FF}
=-\,\frac{\!\COUPLING_\theta'(\Fe,\mm,\theta)}{\det\Fe}\,.
\label{stress-entropy}\end{align}\end{subequations}
The expected {\it symmetry} of such part  $\TT$ of the Cauchy stress is
granted by {\it frame indifference} of $\uppsi(\cdot,\cdot,\theta)$ and
or $\upkappa(\cdot)$. This means that
\begin{align}\nonumber
&\forall (\FF,\mm,\theta)\in{\rm GL}^+(d){\times}\R^d{\times}\R,\ \
Q\in{\rm SO}(d):\ \ \ \ 
\\&\quad\upvarphi(\FF,\mm)=\upvarphi(Q\FF,Q\mm)\,,\ \ 
\COUPLING(\FF,\mm,\theta)=\COUPLING(Q\FF,Q\mm,\theta)\,,
\ \text{ and }\ \upkappa(\FF)=\upkappa(Q\FF)\,,
\label{frame-indifference}
\end{align}
where $Q\in{\rm SO}(d)=\{Q\in\R^{d\times d};\ Q^\top Q=QQ^\top=\bbI\}$ is the
special orthogonal group and ${\rm GL}^+(d)=\{F\in\R^{d\times d};\ \det F>0\}$
denotes the orientation-preserving general linear group. This in particular
implies that the stress $\TT$ is symmetric.
The symmetry of the capillarity contribution $\KK$ to the Cauchy stress
is automatic; actually, this contribution as
$-(\nabla\mm)^\top[\psi_{\rm G}^{}]_{\nabla\mm}'(\FF,\nabla\mm)$ was devised
in \cite[Formula (2.27)]{BePGVa01DDFG} or \cite[Formula (5.16)]{DeSPoG96CTDF}.

Mainly for analytical reasons, we will use also a dissipative contribution
to the Cauchy stress which, together with the conservative part $\TT$,
will realize the {\it Kelvin-Voigt rheological model} and make the system
parabolic. To this goal, we consider a {\it dissipative} contribution to the
{\it Cauchy stress} involving the standard dissipative stress
depending (from the frame-invariancy reason) on the symmetric
velocity gradient and also a higher-order elastic {\it hyper-stress}
$\mathscr{H}$, both isotropic for simplicity:
\begin{align}\nonumber
\DD-{\rm div}\,\mathscr{H}
\ \ \ &\text{ with }\ \DD=\DD(\ee(\vv)) \ \text{ and }\ 
\mathscr{H}=\mathscr{H}(\nabla^2\vv)
\\
&\text{ for }\ \ \DD(\ee)=\DIS|\ee|^{p-2}\ee
\ \text{ and }\ \mathscr{H}({\bm E})=\nu_2|{\bm E}|^{p-2}{\bm E}\,.
\label{Cauchy-dissip}\end{align}
Actually, $\DIS$ and $\nu_2$ may depend on $\det\FF$ and $\theta$ without
causing any structural and analytical problems, but we ignore it
rather for notational simplicity.

The {\it momentum equilibrium} equation then balances the divergence of
the total Cauchy stress with the inertial and gravity force:
\begin{align}
\varrho\DT\vv-{\rm div}\big(\TT{+}\DD{+}\TT_{\!\rm mag}^{}
{-}{\rm div}(\mathscr{H}{+}\mathscr{S})\big)
=\varrho\GRAVITY+\ff_{\!\rm mag}^{}
\label{Euler-thermodynam1-}\end{align}
with $\TT$ from \eq{stress-entropy} and $\DD$ and $\mathscr{H}$ from
\eq{Cauchy-dissip}. Moreover, $\TT_{\!\rm mag}^{}$ and $\ff_{\!\rm mag}^{}$
are the magnetic stress and the magnetic force which balance the
energetics, cf.\ $\TT_{\!\rm mag}^{}:=\KK+\SS$ and
$\mu_0(\nabla\hh)^\top\mm-\mu_0\nabla(\hh{\cdot}\mm)
=:\ff_{\!\rm mag}^{}$ while
$\mathscr{S}:=\upkappa(\FF)
{\rm Skw}(\nabla\mm{\otimes}\mm)/\det\FF$ will be a ``magnetic exchange hyperstress''
in \eq{Euler-thermodynam1}.

The driving magnetic force \eq{driving-field} enters the
Landau-Lifschitz-Gilbert equation \eq{LL-friction} in the previous section.

The third ingredient, i.e.\ \eq{stress-entropy}, is subjected to 
the {\it entropy equation}:
\begin{align}
  \pdt\eta+{\rm div}\big(\vv\,\eta\big)
  =\frac{\xi-{\rm div}\,{\bm j}}\theta\ \ \ \ \text{ with }\ 
  \jj=-\COND(\Fe,\theta)\nabla\theta \,
\label{entropy-eq}\end{align}
and with $\xi=\xi(\FF,\theta;\ee(\vv),\nabla^2\vv,\ZJ\mm)$
denoting the heat production rate specified later in
\eq{Euler-thermodynam3}. The latter equality in \eq{entropy-eq} is the {\it Fourier law}
determining phenomenologically the heat flux $\jj$ proportional to the
negative gradient of temperature $\theta$ through the heat conduction
coefficient $\COND=\COND(\Fe,\theta)$. Assuming $\xi\ge0$ and
$\COND\ge0$ and integrating \eq{entropy-eq} over the domain $\varOmega$
while imposing the non-penetrability of the boundary in the sense that
the normal velocity $\vv{\cdot}\nn$ vanishes across the boundary
$\varGamma$ of $\varOmega$, we obtain the {\it Clausius-Duhem inequality}:
\begin{align}
\frac{\d}{\d t}\int_\varOmega\eta\,\d\xx
=\int_\varOmega\!\!\!\!\!\!\!\!\!\lineunder{\frac\xi\theta+\COND\frac{|\nabla\theta|^2}{\theta^2}}{entropy production rate}\!\!\!\!\!\!\!\!\!\!\!\d\xx
+\int_\varGamma\!\!\!\!\lineunder{\Big(\COND\frac{\nabla\theta}\theta-\eta\vv\Big)}{entropy flux}\!\!\!\!\!\!\cdot\nn\,\d S
\ge\int_\varGamma\!\COND\frac{\nabla\theta{\cdot}\nn}\theta\,\d S\,.
\label{entropy-ineq}\end{align}
If the system is thermally isolated in the sense that the normal heat flux
$\jj{\cdot}\nn$ vanishes across the boundary $\varGamma$, we recover the
{\it 2nd law of thermodynamics}, i.e.\ the total entropy in isolated
systems is nondecreasing in time.

Substituting $\eta$ from \eq{stress-entropy} into \eq{entropy-eq} written
in the form $\theta\DT\eta=\xi-{\rm div}\,\jj-\theta\eta{\rm div}\,\vv$,
we obtain 
\begin{align}\nonumber
&c(\FF,\mm,\theta)\DT\theta
=\xi\big(\FF,\theta;\ee(\vv),\nabla^2\vv,\ZJ\mm\big)
+\theta\Big(\frac{\COUPLING_{\theta}'(\FF,\mm,\theta)}{\det\FF}\Big)_{\!\FF}'
{:}\DT\FF
\\&\hspace{14.5em}\nonumber
+\theta\Big(\frac{\COUPLING_{\theta}'(\FF,\mm,\theta)}{\det\FF}\Big)_{\!\mm}'\!
{\cdot}\DT\mm
+\theta\,\frac{\!\COUPLING_\theta'(\FF,\mm,\theta)\!}{\det\FF}
\,{\rm div}\,\vv-{\rm div}\,{\bm j}
\\[.4em]&\nonumber\qquad\qquad\quad
=\xi\big(\FF,\theta;\ee(\vv),\nabla^2\vv,\ZJ\mm\big)
+\theta\frac{\COUPLING_{\FF\theta}''(\FF,\mm,\theta)}{\det\FF}{:}\DT\FF
+\theta\frac{\COUPLING_{\mm\theta}''(\FF,\mm,\theta)}{\det\FF}
{\cdot}\DT\mm-{\rm div}\,{\bm j}
\\[.1em]&\hspace{10em}\text{ with the heat capacity }\
c(\FF,\mm,\theta)=-\theta\,\frac{\!\COUPLING_{\theta\theta}''(\FF,\mm,\theta)}{\det\FF}\,,
\label{heat-eq+}\end{align}
which can be understood as the {\it heat equation} for the temperature
$\theta$ as an intensive variable.

The referential {\it internal energy} is given by the {\it Gibbs relation}
$\uppsi+\theta\eta$. In our Eulerian formulation, we will need rather the actual
internal energy, which, in view of \eq{stress-entropy}, equals here to 
\begin{align}
\morelinesunder{\frac{\uppsi-\theta\uppsi_\theta'}{\det\FF}}{actual}{internal}{energy}\!\!\!\!
=\!\!\!\!\linesunder{\frac{\upvarphi(\Fe,\mm)}{\det\FF}+\frac{\upkappa(\FF)|\nabla\mm|^2}{2
\det\FF}}
{actual stored and}{exchange energy}\!\!\!\!
+\!\!\!\!\!
\linesunder{\frac{\COUPLING(\Fe,\mm,\theta){-}\theta\COUPLING_\theta'(\Fe,\mm,\theta)}
{\det\FF}}{$=\W\:$ thermal part of}{the internal energy}\!\!\!.\ \ \ \ \ 
\end{align}
In terms of $\W$, the heat equation \eq{heat-eq+} can be written in the
so-called {\it enthalpy formulation}: 
\begin{align}\nonumber
\pdt\W+{\rm div}(\vv\W)&=\xi\big(\FF,\theta;\ee(\vv),\nabla^2\vv,\ZJ\mm\big)
+\frac{\COUPLING'_{\Fe}(\Fe,\mm,\theta)}{\det\Fe}{:}\DT\Fe
+\frac{\COUPLING'_{\mm}(\Fe,\mm,\theta)}{\det\Fe}{\cdot}\DT\mm
-{\rm div}{\bm j}
\\&
\hspace{12em}\text{ with }\ \ \ \W=
\frac{\COUPLING(\Fe,\mm,\theta)-\theta\COUPLING_\theta'(\Fe,\mm,\theta)}{\det\Fe}\,.
\label{Euler-thermodynam3-}\end{align}
Note that $\W$ is an extensive variable so that the left-hand side of
\eq{Euler-thermodynam3-} is not just a convective derivative $\DT\W$.
For the passage from \eq{heat-eq+} to \eq{Euler-thermodynam3-}, we use
the algebra $F^{-1}={\rm Cof}\,F^\top\!/\!\det F$ and the calculus
 $\det'(F)={\rm Cof}\,F$ and \eq{ultimate} so that
$(1/\det\FF)'{:}\DT\FF=-({\rm Cof}\FF/\det\FF^2){:}(\nabla\vv)\FF
=-(\FF^{-\top}\!/\det\FF)\FF^\top
{:}\nabla\vv=({\rm div}\,\vv)/\det\FF$, and thus we can calculate
\begin{align}\nonumber
\pdt\W&+{\rm div}(\vv\W)=\DT\W+\W\,{\rm div}\,\vv
 \\[-.6em]&\nonumber=\ \DT{\overline{\!\!\Big(\frac{\COUPLING(\FF,\mm,\theta)-\theta\COUPLING_\theta'(\FF,\mm,\theta)}{\det\FF}\Big)\!\!}}\
 +\frac{\COUPLING(\Fe,\mm,\theta)-\theta\COUPLING_\theta'(\Fe,\mm,\theta)\!}{\det\Fe}\,{\rm div}\,\vv
\\&\nonumber
=\bigg(\Big(\frac{\COUPLING(\FF,\mm,\theta)}{\det\FF}\Big)_{\!\FF}'
-\theta\Big(\frac{\COUPLING_\theta'(\FF,\mm,\theta)}{\det\FF}\Big)_{\!\FF}'\bigg)
{:}\DT\FF-\theta\frac{\COUPLING_{\theta\theta}''(\Fe,\mm,\theta)}{\det\FF}\DT\theta
\\&\nonumber\ \ 
+\bigg(\frac{\COUPLING_\mm'(\Fe,\mm,\theta)}{\det\FF}
-\theta\frac{\COUPLING_{\mm\theta}''(\Fe,\mm,\theta)}{\det\FF}\bigg)
{\cdot}\DT\mm
+\frac{\COUPLING(\Fe,\mm,\theta)-\theta\COUPLING_\theta'(\Fe,\mm,\theta)\!}{\det\Fe}\,{\rm div}\,\vv
\\&\nonumber
=\bigg(\frac{\COUPLING_{\!\FF}'(\FF,\mm,\theta)}{\det\FF}
-\theta\frac{\COUPLING_{\FF\theta}''(\FF,\mm,\theta)}{\det\FF}
+\big(\COUPLING_{\!\FF}'(\FF,\mm,\theta)-\theta\COUPLING_\theta'(\FF,\mm,\theta)\big)\Big(\frac1{\det\FF}\Big)'\bigg){:}\DT\FF
\\&\nonumber\ \ 
+\frac{\COUPLING_\mm'(\Fe,\mm,\theta){-}\theta\COUPLING_{\mm\theta}''(\Fe,\mm,\theta)}{\det\FF}
{\cdot}\DT\mm
+c(\Fe,\mm,\theta)\DT\theta
+\frac{\COUPLING(\Fe,\mm,\theta){-}\theta\COUPLING_\theta'(\Fe,\mm,\theta)\!}{\det\Fe}\,{\rm div}\,\vv
\\&\nonumber
=\frac{\COUPLING_{\!\FF}'(\FF,\mm,\theta){-}\theta\COUPLING_{\FF\theta}''(\FF,\mm,\theta)}{\det\FF}
{:}\DT\FF
+\frac{\COUPLING_\mm'(\Fe,\mm,\theta){-}\theta\COUPLING_{\mm\theta}''(\Fe,\mm,\theta)}{\det\FF}
{\cdot}\DT\mm
+c(\Fe,\mm,\theta)\DT\theta\,.
\end{align}
This shows that \eq{heat-eq+} are \eq{Euler-thermodynam3-} indeed
(formally) equivalent to each other. Alternatively in \eq{Euler-thermodynam3-},
we could use \eq{ultimate} for writing 
$\COUPLING'_{\Fe}(\Fe,\mm,\theta){:}\DT\Fe
=\COUPLING'_{\Fe}(\Fe,\mm,\theta)\Fe^\top\!{:}(\nabla\vv)$
$=\COUPLING'_{\Fe}(\Fe,\mm,\theta)\Fe^\top\!{:}\ee(\vv)$; here the assumed
frame indifference \eq{frame-indifference} of $\COUPLING(\cdot,\mm,\theta)$
itself, leading to symmetry of $\COUPLING'_{\Fe}(\Fe,\mm,\theta)\Fe^\top$, was
employed.

\begin{remark}[{\sl Gradient theories in rates}]\label{rem-grad}\upshape
So-called gradient theories in continuum mechanical models are nowadays very
standard, referred as {\it nonsimple materials}, determining some internal
length scales and often facilitating mathematical analysis. They can be applied to
the conservative stress through the free energy or to the dissipation stress.
Here, we have used the latter option in \eq{Cauchy-dissip} which is better
fitted to the rate formulation and which can make velocity field enough
regular, as vitally needed for the transport of $\varrho$ and $\FF$ in
the Eulerian models. The higher-gradient hyper-stress as used below in
\eq{Cauchy-dissip}
follows the theory by E.~Fried and M.~Gurtin \cite{FriGur06TBBC}, as already
articulated in the general nonlinear context of {\it multipolar fluids} by
J.~Ne\v cas at al. \ \cite{Neca94TMF,NeNoSi91GSCI,NecRuz92GSIV}
or of solids \cite{Ruzi92MPTM,Silh92MVMS}, inspired by 
R.A.\,Toupin \cite{Toup62EMCS} and R.D.\,Mindlin \cite{Mind64MSLE}.
\end{remark}

\begin{remark}[{\sl Dry friction in magnetization evolution}]\label{rem-dry-friction}
\upshape
Dry-friction-type rate-independent dissipation was proposed in
\cite{BalHel91DFM,Visi97MLLE}
as a device to model properly hysteresis in ferromagnets, modifying the Landau-Lifschitz
equation by augmenting suitably the effective magnetic field.
Although the original Gilbert's \cite{Gilb55LFGE} and Landau-Lifschitz' \cite{LanLif35TDMP}
equations are equivalent with each
other, the resulting augmented equations proposed in \cite{BalHel91DFM} and
\cite{Visi97MLLE} are no longer mutually equivalent. This has been pointed out
in \cite{PoGu01DMM}, where the conceptual differences between the Gilbert and
the Landau-Lifschitz formats have been elucidated.
\end{remark}

\begin{remark}[{\it Zaremba-Jaumann derivative $\ZJ\mm$}]\label{rem-ZJ}\upshape
In deformable (and deforming) magnetic medium, the Jaumann corotational
derivative for magnetization was suggested already by Maugin
\cite{Maug76CTDF} to model situations when the magnetization can
be ``frozen'' in hard-magnetic materials in their ferro- or ferri-magnetic
state. Later it was used in \cite{DeSPoG95ISIS,DeSPoG96CTDF} in the linear
viscosity (magnetic attenuation) term. It should be noted that the
gyromagnetic term $\mm{\otimes}\DT\mm/\gamma(\FF,\mm,\theta)$ can be seen in
literature; yet, it is the same as used in \eq{LL-friction} since
$\mm{\otimes}\DT\mm-\mm{\otimes}\ZJ\mm=\mm{\otimes}({\rm skw}(\nabla\vv)\mm)
={\rm skw}(\nabla\vv){:}(\mm{\otimes}\mm)=0$.
Instead of ${\rm skw}(\nabla\vv)$ in \eq{def-ZJ}, an angular velocity either
governed by a separate parabolic equation or approximated by as ${\rm curl}\,\vv/2$
was used in \cite{GruWei21FITM,NoSaTo16EFMN,Rose87MF,Scro19GWPC}. Usage of
another (Lie) derivative was proposed in \cite[Formula (74)]{VaPaEs21MTCM}.
Mere convective derivative for magnetization has been used in 
\cite{BFLS18EWSE,KaKoSc21MAWS,SchZab18USMM} to model rather (incompressible
isothermal) fluids containing magnetic particles.
\end{remark}

\begin{example}[{\sl Neo-Hookean elastic magnets}]\label{exa-neo-Hook} \upshape
Modifying slightly the ``rigid'' model \eq{exa-phi} 
and expanding it by standard neo-Hookean elastic ansatz, one obtains an example
for elastic magnetic material amenable for ferro-to-paramagnetic transition
and for complying with the (\ref{Euler-ass}b--e) below:
\begin{align}\nonumber
\uppsi(\FF,\mm,\nabla\mm,\theta)&=\frac12G\Big(\frac{{\rm tr}(\FF\FF^\top)}{(\det\FF)^{2/d}}-d\Big)
+v(\det\FF)
+\frac{\upkappa(\FF)}2|\nabla\mm|^2
\\&
+a_0(\det\FF)\Big(\frac\theta{1{+}\epsilon_1\theta}
{-}\frac{\theta_\text{\sc c}}{1{+}\epsilon_1\theta_\text{\sc c}}\Big)\frac{|\mm|^2}{1{+}\epsilon_2|\mm|^2}
+b_0(\det\FF)|\mm|^4+c_{0}\theta(1{-}{\rm ln}\theta)
\nonumber\end{align}
with some (referential) heat capacity $c_0>0$, some
$\epsilon_1,\epsilon_2>0$, some shear modulus $G>0$, and
the non-negative volumetric energy $v\in C^1(\R^+)$ 
and $a_0(J)\ge\delta J$ and $b_0(J)\ge\delta J$ for some $\delta>0$
so that $\upvarphi$ fulfills (up to an irrelevant constant) the coercivity
\eq{Euler-ass-phi} with $s=4$. Note
that $\uppsi(\FF,\mm,\nabla\mm,\cdot)$ is concave for $c_0,\varsigma>0$ and $\epsilon_1,\epsilon_2\ge0$.
Then, the split \eq{ansatz-} uses
\begin{align}
\COUPLING(\FF,\mm,\theta)=
\frac{a_0(\det\FF)\theta|\mm|^2}{(1{+}\epsilon_1\theta)(1{+}\epsilon_2|\mm|^2)}
+c_{0}\theta(1{-}{\rm ln}\theta)
\end{align}
so that
\begin{align}
\OMEGA(\FF,\mm,\theta)&=\frac{\COUPLING(\Fe,\mm,\theta)
-\theta\COUPLING_\theta'(\Fe,\mm,\theta)}{\det\FF}
=\frac{c_{0}\theta}{\det\FF}
+\frac{\epsilon_1a_0(\det\FF)\theta^2|\mm|^2}{(1{+}\epsilon_1\theta)^2(1{+}\epsilon_2|\mm|^2)\det\FF}
\end{align}
and the (actual)  heat capacity
$c(\FF,\mm,\theta)=-\theta\COUPLING_{\theta\theta}''(\Fe,\mm,\theta)/\!\det\FF$ is
\begin{align}
c(\FF,\mm,\theta)=\OMEGA_\theta'(\FF,\mm,\theta)
=\frac{c_{0}}{\det\FF}+
\frac{2\epsilon_1a_0(\det\FF)\theta|\mm|^2}{(1{+}\epsilon_1\theta)^3(1{+}\epsilon_2|\mm|^2)\det\FF}\,.
\end{align}
Note that ``thermo-coupling'' stress
$\COUPLING_\FF'(\FF,\mm,\theta)\FF^\top\!/\!\det\FF=
a_0'(\det\FF)\theta|\mm|^2\bbI/((1{+}\epsilon_1\theta)(1{+}\epsilon_2|\mm|^2))$
is bounded provided $a_0'$ is bounded on ${\rm GL}^+(d)$, so it
surely complies with \eq{Euler-ass-adiab} below.
Also $|\COUPLING_\mm'(\FF,\cdot,\cdot)/\!\det\FF|$ is bounded for $\FF$ ranging over compact
sets in ${\rm GL}^+(d)$, so it surely complies with \eq{Euler-ass-adiab};
here we use $\epsilon_1>0$ and $\epsilon_2>0$.
Also this ansatz satisfies \eq{Euler-ass-primitive-c}.
Moreover, $|\OMEGA_\FF'(\FF,\mm,\cdot)|$  has at
most linear growth while and $|\OMEGA_\mm'(\FF,\mm,\cdot)|$ is even bounded 
as well as $|\OMEGA_{\FF\theta}''|$ and  $|\OMEGA_{\mm\theta}''|$ for $\FF$ ranging over compact
sets in ${\rm GL}^+(d)$, so that \eq{Euler-ass-primitive-c} is satisfied, too.\COMMENT{CHECK}
\end{example}

\subsection{The thermo-magneto-mechanical system and its energetics}\label{sec-engr}
Let us summarize the thermodynamically coupled system composed of
six partial differential equations for $\varrho$, $\vv$, $\FF$, $\mm$, $u$,
and $\theta$. More specifically, it is composed from the mass continuity
equation for $\varrho$, the momentum equation written in terms of velocity
$\vv$, the evolution-and-transport of the deformation-gradient tensor $\FF$,
a flow rule (as an inclusion) for the magnetization $\mm$,
the Poisson equation for the demagnetizing-field potential $u$, and
the heat-transfer equation for temperature $\theta$.

Altogether, merging \eq{u-eq}, \eq{ultimate}, \eq{cont-eq+},
\eq{LL-friction}, \eq{Euler-thermodynam1-}, and \eq{Euler-thermodynam3-} with
\eq{driving-fields}, we obtain a system of six equations for
$(\varrho,\vv,\FF,\mm,u,\W)$, respectively also for $\theta$: 
\begin{subequations}\label{Euler-thermodynam}
\begin{align}\label{Euler-thermodynam0}
&\pdt\varrho=-\,{\rm div}(\varrho\vv)\,,
\\\nonumber
&\pdt{}(\varrho\vv)={\rm div}\Big(\TT{+}\KK{+}\SS{+}
\DD{-}{\rm div}(\mathscr{H}{+}\mathscr{S})-\varrho\vv{\otimes}\vv\Big)
  +\mu_0(\nabla\hh)^\top\mm-\mu_0\nabla(\hh{\cdot}\mm)+\varrho\GRAVITY
 \\[-.1em]\nonumber
    &\hspace*{1em}
\ \text{ with }\ \TT
=\Big(\frac{\ \upvarphi_{\Fe}'(\Fe,\mm)
+\COUPLING_{\Fe}'(\Fe,\mm,\theta)\!}{\det\Fe}+\frac{|\nabla\mm|^2\upkappa'(\FF)\!}{2\det\Fe}\:\Big)
\Fe^\top\,,\ \ \ \hh=\hh_{\rm ext}+\nabla u\,,
\\[.1em]\nonumber
    &\hspace*{4em}\KK=
\frac{\upkappa(\FF)}{\det\FF}\DELETE{\Big(}\nabla\mm{\otimes}\nabla\mm\DELETE{{-}\frac12|\nabla\mm|^2\bbI\Big)}
\,,\ \ \DD=\DIS|\EE(\vv)|^{p-2}\EE(\vv)\,,\ \ \ \  \mathscr{H}=\nu_2|\nabla^2\vv|^{p-2}\nabla^2\vv\,,
\\[.1em]
    &
    \hspace*{4em}\SS={\rm skw}\Big(\Big(\mu_0\hh {-}
 \frac{\uppsi_\mm'(\Fe,\mm,\theta)\!}{\det\Fe}\Big){\otimes}\mm\Big)\,,
\ \text{ and }\ \mathscr{S}=
\frac{\upkappa(\FF)}{\det\FF}{\rm Skw}\big(\mm{\otimes}\nabla\mm),\!\!\!\!
    \label{Euler-thermodynam1}
    \\[-.2em]
&\pdt\Fe=(\Nabla\vv)\Fe-(\vv{\cdot}\nabla)\Fe\,,
\label{Euler-thermodynam2}
\\[-.2em]&
\tau\ZJ\mm+\HC(\FF,\theta){\rm Dir}(\ZJ\mm)-\frac{\mm{\times}\ZJ\mm}{\gamma(\FF,\mm,\theta)}
\ni\mu_0\hh-\frac{\uppsi_\mm'(\Fe,\mm,\theta)\!}{\det\Fe}
+{\rm div}\Big(\frac{\upkappa(\FF)}{\det\FF}\nabla\mm\Big)\,,
\label{Euler-thermodynam4}\\&\label{Euler-thermodynam5}
\Delta u={\rm div}(\chi_\varOmega^{}\mm)\ \ \ \text{ on }\R^d\,,
\\&\nonumber
\pdt{\W}
=\xi(\FF,\theta;\EE(\vv),\nabla^2\vv,\ZJ\mm)
+\frac{\COUPLING'_{\Fe}(\Fe,\mm,\theta)\Fe^\top\!\!}{\det\Fe}{:}\ee(\vv)
+\frac{\COUPLING_\mm'(\FF,\mm,\theta)}{\det\FF}{\cdot}\DT\mm
-{\rm div}\big(\jj{+}\W\vv\big)
\\&\nonumber\hspace{8em}
\text{ with }\
\xi(\FF,\theta;\ee,\bm{G},\rr)=\DIS|\ee|^p+\nu_2|\bm{G}|^p+\mu_0\tau
|\rr|^2+\mu_0\HC(\FF,\theta)|\rr|
\\&\nonumber\hspace{8em}\text{ and }\ \jj=-\COND(\FF,\theta)\nabla\theta\,,
\\[-.4em]&\hspace{8em}\text{ and }\ \W=\OMEGA(\Fe,\mm,\theta)
=\frac{\COUPLING(\Fe,\mm,\theta)-\theta\COUPLING_\theta'(\Fe,\mm,\theta)}{\det\Fe}\,.
\label{Euler-thermodynam3}
\end{align}\end{subequations}
The product $\nabla\mm{\otimes}\nabla\mm$ in \eq{Euler-thermodynam1}
is to be understood componentwise, specifically  $[\nabla\mm{\otimes}\nabla\mm]_{ij}$
$=\sum_{k=1}^d\frac{\pl}{\pl x_i}m_k\frac{\pl}{\pl x_j}m_k$ with $\mm=(m_1,...,m_d)$
and $\xx=(x_1,...,x_d)$, while the skew-symmetric part ``Skw'' of the
3rd-order tensor is defined as
\begin{align}
\big[{\rm Skw}(\mm{\otimes}\nabla\mm)\big]_{ijk}
:=\frac12\Big(m_i\frac{\pl m_j}{\pl x_k}-m_j\frac{\pl m_i}{\pl x_k}\Big)\,.
\end{align}
The equations (\ref{Euler-thermodynam}a-d,f) are
considered on the domain $\varOmega$ while the static equation
\eqref{Euler-thermodynam5} is to hold on the
whole Universe in the distributional sense at all time instants.
The symmetric Korteweg-like stress
$\KK=\upkappa(\FF)\nabla\mm{\otimes}\nabla\mm/\!\det\FF$ in
\eq{Euler-thermodynam1} occur e.g.\ in \cite{BePGVa01DDFG,Eric02EETM}.
The skew-symmetric stress $\SS$ and couple-like hyperstress $\mathscr{S}$ come
from the calculation \eq{test-damage}; for a similar skew-symmetric stress $\SS$ see 
\cite[Formula (33)]{DeSPoG95ISIS} or \cite[Formula (5.37)]{DeSPoG96CTDF}
while the skew-symmetric hyperstress
$\mathscr{S}$ is like in the Cosserat theory in \cite{Toup64TECS}.
The pressures $\mu_0\hh{\cdot}\mm$ in the momentum equation
is related with that the Zeeman and the demagnetizing-field energies are actual
(not referential).
The magnetic, so-called Kelvin force $\mu_0(\nabla\hh)^\top\mm$ comes
from the magnetization $\mm$ exposed to the nonuniform magnetic field $\hh$,
and it arises from the calculus \eq{formula3} and \eq{calculus-PM4} below.

Denoting by $\nn$ the unit outward normal to the (fixed) boundary $\varGamma$
of the domain $\varOmega$, we complete this system by suitable boundary conditions: 
\begin{subequations}\label{Euler-thermodynam-BC}
\begin{align}\label{Euler-thermodynam-BC-1}
&\!\!\vv{\cdot}\nn=0,\ \ \ 
\big[(\TT{+}\KK{+}\SS
{+}\DD{-}{\rm div}(\mathscr{H}{+}\mathscr{S}))\nn{-}\divS\big(\mathscr{H}\nn{+}\mathscr{S}
\nn\big)\big]_\text{\sc t}^{}\!+\nu_\flat|\vv|^{p-2}\vv
=\TRACTION\,,\ \
\\&\label{Euler-thermodynam-BC-2}
\!\!\nabla^2\vv{:}(\nn{\otimes}\nn)={\bm0}\,,
\ \ \ (\nn{\cdot}\nabla)\mm=\bm{0}\,,\ \ \ u(\infty)=0\,,
\ \ \ \text{ and }\ \ \ 
\jj{\cdot}\nn=h(\theta)+\frac{\nu_\flat}2|\vv|^p
\end{align}\end{subequations}
with $\nu_\flat>0$ a boundary viscosity coefficient and with
$[\,\cdot\,]_\text{\sc t}^{}$ a tangential part of a vector
and with $\divS={\rm tr}(\nablaS)$ denoting the $(d{-}1)$-dimensional
surface divergence with ${\rm tr}(\cdot)$ being the trace of a
$(d{-}1){\times}(d{-}1)$-matrix and
$\nablaS v=\nabla v-\frac{\partial v}{\partial\nn}\nn$ 
being the surface gradient of $v$.
Naturally, $\TRACTION{\cdot}\nn=0$ is to be assumed
if we want to recover the boundary conditions \eqref{Euler-thermodynam-BC-1} in the
classical form, otherwise the weak form does not directly need it.

The first condition (i.e.\ normal velocity zero) expresses 
nonpenetrability of the boundary was used already for \eq{entropy-ineq}
and is most frequently adopted in literature for Eulerian formulation.
This simplifying assumption fixes the shape of $\varOmega$ in its
referential configuration allows also for considering fixed boundary even
for such time-evolving Eulerian description. The latter condition in
\eq{Euler-thermodynam-BC-1} involving a boundary viscosity comes from
the Navier boundary condition largely used in fluid dynamics and is here
connected with the technique used below, which is based on the total energy
balance as the departing point and which, unfortunately, does not allow to cope
with $\nu_\flat=0$ and simultaneously $\TRACTION\ne0$. This boundary viscosity
naturally may contribute to the heat production on the boundary as well as
to the outflow of the heat energy to the outer space. For notational
simplicity, we consider that it is just equally distributed, one part
remaining on the boundary of $\varOmega$ and the other part 
leaving outside, which is related with the coefficient $1/2$ in 
the last condition in \eq{Euler-thermodynam-BC-2}.
The condition $u(\infty)=0$ in \eq{Euler-thermodynam-BC-2} expresses shortly
that $\lim_{|\xx|\to0}u(\xx)=0$.

The magnetization flow rule \eqref{Euler-thermodynam4} with the corotational
derivative $\ZJ\mm$ in 
 see also \cite{BePGVa01DDFG,Maug76CTDF} where it is articulated that 
the magnetization is ``frozen'' in the deforming medium if $\ZJ\mm=0$ 
which then means that the magnetization is transported and rotates at
the same local rate as the deforming medium; this is the situation below
the blocking temperature $\theta_{\rm b}$ and when the total driving
magnetic field has small magnitude.

For the capillarity-like stress $\KK$ and and the skew-symmetric stress $\SS$
see also \cite{BePGVa01DDFG,DeSPoG95ISIS}.

The magneto-mechanical energy balance of the model can be seen when testing the
momentum equation \eq{Euler-thermodynam1} by $\vv$ while using the continuity
equation \eq{Euler-thermodynam0} tested by $|\vv|^2/2$ and the
evolution-and-transport equation \eq{Euler-thermodynam2} for $\FF$ tested
by the stress $[\upvarphi(\FF,\mm)/\!\det\FF]_\FF'$, the magnetic
Landau-Lifshitz-Gilberg equation \eq{Euler-thermodynam4}
by $\ZJ\mm$, and the (rest from the) Maxwell system \eq{Euler-thermodynam5}
by $\mu_0\pdt{u}$.

Let us first make the test of \eq{Euler-thermodynam1} by $\vv$. Using again
the algebra $F^{-1}={\rm Cof}\,F^\top\!/\!\det\,F$ and the calculus
$\det'(F)={\rm Cof}\,F$, we can write the part of the Cauchy stress arising
from the stored energy as 
\begin{align}\nonumber
&\frac{\upvarphi_\FF'(\FF,\mm)}{\det\FF}\FF^\top
=\frac{\upvarphi_\FF'(\FF,\mm)-\upvarphi(\FF,\mm)\FF^{-\top}\!\!\!}{\det\FF}\,\FF^\top\!+
\frac{\upvarphi(\FF,\mm)}{\det\FF}\bbI
\\&\ \ 
=\bigg(\frac{\upvarphi_\FF'(\FF,\mm)\!}{\det\FF}
-\frac{\!\upvarphi(\FF,\mm){\rm Cof}\FF}{(\det\FF)^2}\bigg)\FF^\top\!\!+
\frac{\upvarphi(\FF,\mm)}{\det\FF}\bbI
\Big[\frac{\upvarphi(\FF,\mm)}{\det\FF}\Big]_{\!\FF}'\FF^\top\!\!+\frac{\upvarphi(\FF,\mm)}{\det\FF}\bbI\,.
\label{referential-stress}\end{align}
Let us recall that $\upvarphi(\FF,\mm)/\!\det\FF$ in \eq{referential-stress}
is the stored energy per actual (not referential) volume. Using the calculus
\eq{referential-stress}, we obtain 
\begin{align}\nonumber
&\!\!\int_\varOmega{\rm div}\,\TT{\cdot}\vv\,\d\xx
=\!\int_\varGamma(\TT\nn){\cdot}\vv\,\d S-\!\int_\varOmega\!\TT{:}\ee(\vv)\,\d\xx
\ \ \ \ \text{ with}
\\&\!\!\nonumber\int_\varOmega\!\TT{:}\ee(\vv)\,\d\xx=
\!\int_\varOmega\!\Big(\frac{\upvarphi_\FF'(\FF,\mm)
+\COUPLING_\FF'(\FF,\mm,\theta)}{\det\FF}
+\frac{|\nabla\mm|^2\upkappa'(\FF)}{2\det\FF}\Big)\FF^\top{:}\ee(\vv)\,\d\xx
\\&\!\!
\nonumber
=\!\int_\varOmega\!\Big(\Big[\frac{\upvarphi(\FF,\mm)}{\det\FF}\Big]_{\!\FF}'\FF^\top
+\frac{\upvarphi(\FF,\mm)}{\det\FF}\bbI
+\frac{\COUPLING_\FF'(\FF,\mm,\theta)}{\det\FF}\FF^\top+\frac{|\nabla\mm|^2\upkappa'(\FF)}{2\det\FF}\FF^\top\Big){:}\ee(\vv)\,\d\xx
\\&\!\!\nonumber
=\!\int_\varOmega\!\Big[\frac{\upvarphi(\FF,\mm)}{\det\FF}\Big]_{\!\FF}'
{:}(\Nabla\vv)\FF
+\frac{\upvarphi(\FF,\mm)}{\det\FF\!}\,{\rm div}\,\vv
+\Big(\frac{\COUPLING_\FF'(\FF,\mm,\theta)\!\!}{\det\FF}
+\frac{|\nabla\mm|^2\upkappa'(\FF)}{2\det\FF}\Big)\FF^\top\!{:}\ee(\vv)\,\d\xx
\\[-.0em]&\!\!
=\frac{\d}{\d t}\int_\varOmega\!\frac{\upvarphi(\FF,\mm)\!\!}{\det\FF}\,\d\xx
-\!\int_\varOmega\!\frac{\upvarphi_\mm'(\FF,\mm)\!}{\det\FF}{\cdot}\DT\mm
-\Big(\frac{\COUPLING_\FF'(\FF,\mm,\theta)\!\!}{\det\FF}+\frac{|\nabla\mm|^2\upkappa'(\FF)}{2\det\FF}\Big)\FF^\top\!{:}\ee(\vv)
\,\d\xx.\!\!
\label{Euler-large-thermo}\end{align}
Here, we used the matrix algebra
$A{:}(BC)=(B^\top\!A){:}C=(AC^\top){:}B$ for any square matrices $A$, $B$, and $C$
and also we used \eq{Euler-thermodynam2} together with the Green formula
and the nonpenetrability boundary condition for
\begin{align}\nonumber
&\int_\varOmega\!\Big[\frac{\upvarphi(\FF,\mm)}{\det\FF}\Big]_{\!\FF}'
{:}(\Nabla\vv)\FF
+\frac{\upvarphi_\mm'(\FF,\mm)}{\det\FF}{\cdot}\DT\mm
+\frac{\upvarphi(\FF,\mm)}{\det\FF\!}\,{\rm div}\,\vv\,\d\xx
\\&\nonumber\hspace*{1em}
\!\!\!\!\stackrel{\eq{Euler-thermodynam2}}{=}\!\!\!
\int_\varOmega\!\Big[\frac{\upvarphi(\FF,\mm)}{\det\FF}\Big]_{\!\FF}'
{:}\Big(\pdt\FF{+}(\vv{\cdot}\nabla)\FF\Big)
+\frac{\upvarphi_\mm'(\FF,\mm)}{\det\FF}{\cdot}\Big(\pdt\mm
{+}(\vv{\cdot}\nabla)\mm\Big)
+\frac{\upvarphi(\FF,\mm)\!}{\det\FF}\,{\rm div}\,\vv\,\d\xx
\\&\nonumber\hspace*{1em}
=\frac{\d}{\d t}\int_\varOmega\frac{\upvarphi(\FF,\mm)}{\det\FF}\,\d\xx
+\!\int_\varOmega\!\nabla\Big(\frac{\upvarphi(\FF,\mm)}{\det\FF}\Big){\cdot}\vv+\frac{\upvarphi(\FF,\mm)\!}{\det\FF}\,{\rm div}\,\vv\,\d\xx
\\&\hspace*{1em}
=\frac{\d}{\d t}\int_\varOmega\frac{\upvarphi(\FF,\mm)}{\det\FF}\,\d\xx
+\!\int_\varGamma\frac{\upvarphi(\FF,\mm)}{\det\FF}(\hspace*{-.7em}\lineunder{\vv{\cdot}\nn}{$=0$}\hspace*{-.7em})\,\d S\,,
\label{formula4-}\end{align}
where we used
\begin{align}\nonumber
&\int_\varOmega\Big[\frac{\upvarphi(\FF,\mm)}{\det\FF}\Big]_{\!\FF}'{:}(\vv{\cdot}\Nabla)\FF
  +\frac{\upvarphi_{\mm}'(\FF,\mm)}{\det\FF}{\cdot}(\vv{\cdot}\Nabla)\mm\,\d\xx
\\[-.4em]&\qquad
=\int_\varOmega\Nabla\Big(\frac{\upvarphi(\FF,\mm)}{\det\FF}\Big){\cdot}\vv\,\d\xx
 =\int_\varGamma
 \frac{\upvarphi(\FF,\mm)}{\det\FF}\!\!\!\lineunder{\vv{\cdot}\nn}{=0}\!\!\!\!\!\,\d S-
\int_\varOmega\frac{\upvarphi(\FF,\mm)}{\det\FF}{\rm div}\,\vv\,\d\xx\,.
\label{formula4}\end{align}

The further contribution from the dissipative part of the Cauchy stress
uses Green's formula over $\varOmega$ twice and the surface Green formula
over $\varGamma$. We abbreviate the elastic hyperstress
$\mathscr{H}=\nu_2|\nabla^2\vv|^{p-2}\nabla^2\vv$. Then
\begin{align}\nonumber
&\int_\varOmega{\rm div}
\big(\DIS|\EE(\vv)|^{p-2}\EE(\vv)-{\rm div}\mathscr{H}\big){\cdot}\vv\,\d\xx
\\[-.3em]&\nonumber
=\int_\varGamma\vv{\cdot}
\big(\DIS|\EE(\vv)|^{p-2}\EE(\vv){-}{\rm div}\mathscr{H}\big)\nn\,\d S
-\int_\varOmega\!\big(\DIS|\EE(\vv)|^{p-2}\EE(\vv)-{\rm div}\mathscr{H}\big){:}\Nabla\vv\,\d\xx
\\&\nonumber=\int_\varGamma\vv{\cdot}
\big(\DIS|\EE(\vv)|^{p-2}\EE(\vv){-}{\rm div}\mathscr{H}\big)\nn
-\nn{\cdot}\mathscr{H}{:}\Nabla\vv\,\d S
-\!\int_\varOmega\!\DIS|\EE(\vv)|^p+\mathscr{H}\Vdots\Nabla^2\vv\,\d\xx
\\&\nonumber=\int_\varGamma\!
\mathscr{H}{\Vdots}(\pl_\nn\vv{\otimes}\nn{\otimes}\nn)
+\nn{\cdot}\mathscr{H}{:}\NablaS\vv+\vv{\cdot}
\big(\DIS|\EE(\vv)|^{p-2}\EE(\vv){-}{\rm div}\mathscr{H}\big)\nn\,\d S
\\[-.4em]&\nonumber\hspace{22em}
-\!\int_\varOmega\!\DIS|\EE(\vv)|^p\!+\nu_2|\nabla^2\vv|^p\,\d\xx
\\&\nonumber
=\int_\varGamma\mathscr{H}{\Vdots}(\pl_\nn\vv{\otimes}\nn{\otimes}\nn)-
\big(\divS(\nn{\cdot}\mathscr{H})+
\big(\DIS|\EE(\vv)|^{p-2}\EE(\vv){-}{\rm div}\mathscr{H}\big)\nn\big)
{\cdot}\vv\,\d S
\\[-.4em]&\hspace{22em}
-\!\int_\varOmega\!\DIS|\EE(\vv)|^p+\nu_2|\nabla^2\vv|^p\,\d\xx\,,
\label{Euler-test-momentum++}\end{align}
where we used the decomposition of $\Nabla\vv$ into its normal component
$\pl_\nn\vv$ and the tangential component, i.e.\ written componentwise
$\nabla\vv_i=(\nn{\cdot}\nabla\vv_i)\nn+\nablaS\vv_i$.

Furthermore, the inertial force 
$\pdt{}(\varrho\vv)+{\rm div}(\varrho\vv{\otimes}\vv)$ in
\eq{Euler-thermodynam1}
tested by $\vv$ gives the rate of kinetic energy $\varrho|\vv|^2/2$
integrating over $\varOmega$.
Here we use the continuity equation \eq{cont-eq+} tested by $|\vv|^2/2$
and the Green formula with the boundary condition $\vv{\cdot}\nn=0$: 
\begin{align}
&\int_\varOmega\Big(\pdt{}(\varrho\vv)+{\rm div}(\varrho\vv{\otimes}\vv)\Big){\cdot}\vv\,\d\xx=\int_\varOmega\varrho\DT\vv{\cdot}\vv\,\d\xx
=\frac{\d}{\d t}\int_\varOmega\frac\varrho2|\vv|^2\,\d\xx
+\!\int_\varGamma\varrho|\vv|^2\hspace{-.7em}\lineunder{\vv{\cdot}\nn}{$=0$}\hspace{-.7em}\d S\,.
\label{calculus-convective-in-F}
\end{align}

\def\MU{}

The test of \eqref{Euler-thermodynam4} by $\MU\ZJ\mm$ is quite technical. The
exchange-energy term ${\rm div}(\upkappa(\FF)\nabla\mm)/\!\det\FF)$ tested by
$\MU(\vv{\cdot}\Nabla)\mm$ is to be handled by using Green's formula twice. Namely,
\begin{align}\nonumber
&\int_\varOmega\!{\rm div}\Big(\frac{\upkappa(\FF)\nabla\mm}{\det\FF}\Big)
{\cdot}\MU(\vv{\cdot}\Nabla)\mm\,\d\xx
=\MU\!\int_\varGamma\frac{\upkappa(\FF)(\nn{\cdot}\Nabla)\mm}{\det\FF}
{\cdot}(\vv{\cdot}\Nabla)\mm\,\d S
\\[-.1em]&\nonumber\hspace*{9em}
-\MU\!\int_\varOmega\frac{\upkappa(\FF)\Nabla^2\mm}{\det\FF}\Vdots(\vv{\otimes}\Nabla\mm)
+\frac{\upkappa(\FF)(\Nabla\mm{\otimes}\Nabla\mm)}{\det\FF}{:}\EE(\vv)\,\d\xx
\\&\nonumber
=\MU\!\int_\varGamma
\frac{\upkappa(\FF)(\nn{\cdot}\Nabla)\mm}{\det\FF}{\cdot}\big((\vv{\cdot}\Nabla)\mm\big)-
\frac{\upkappa(\FF)|\Nabla\mm|^2}{2\det\FF}\vv{\cdot}\nn\,\d S
\\[-.1em]&\hspace*{1em}+\MU\!\int_\varOmega
\frac{|\Nabla\mm|^2}{2}\nabla\Big(\frac{\upkappa(\FF)}{\det\FF}\Big){\cdot}\vv
+\frac{\upkappa(\FF)|\Nabla\mm|^2}{2\det\FF}{\rm div}\,\vv
-\hspace*{-.7em}\lineunder{\frac{\upkappa(\FF)(\Nabla\mm{\otimes}\Nabla\mm)}{\det\FF}}{$=\KK$ from \eqref{Euler-thermodynam1}}\hspace*{-.7em}{:}\EE(\vv)
\,\d\xx\,,
\label{test-damage}\end{align}
where the boundary integral vanishes due to the boundary conditions
$(\nn{\cdot}\Nabla)\mm=\mathbf0$ and $\vv{\cdot}\nn=0$.
The latter equality in \eq{test-damage} follows by the
calculus and the Green formula:
\begin{align}\nonumber
&\int_\varOmega\frac{\upkappa(\FF)\Nabla^2\mm}{\det\FF}\Vdots(\vv{\otimes}\Nabla\mm)
\,\d\xx
=\int_\varGamma\frac{\upkappa(\FF)(\nn{\cdot}\Nabla)\mm}{\det\FF}
{\cdot}\big((\vv{\cdot}\Nabla)\mm\big)\,\d S
\\&\nonumber\ \ 
-\!\int_\varOmega
\frac{\upkappa(\FF)\nabla\mm{\otimes}\vv}{\det\FF}\Vdots\Nabla^2\mm
+|\nabla\mm|^{2}\Big(
\frac{\upkappa(\FF)}{\det\FF}{\rm div}\,\vv+
\nabla\Big(\frac{\upkappa(\FF)}{\det\FF}\Big){\cdot}\vv\Big)\,\d\xx
\\&\hspace{9em}=-\!\int_\varOmega\frac{\upkappa(\FF)|\nabla\mm|^{2}}{2\det\FF}{\rm div}\,\vv+
\frac{|\nabla\mm|^{2}}{2}\nabla\Big(\frac{\upkappa(\FF)}{\det\FF}\Big)
{\cdot}\vv\,\d\xx\,,
\end{align}
where again $(\nn{\cdot}\Nabla)\mm=\mathbf0$ was used.
For the last term, we can still use the calculus
\begin{align}
\!\!\nabla\Big(\frac{\upkappa(\FF)}{\det\FF}\Big)&=\Big(\frac{\upkappa(\FF)}{\det\FF}\Big)'{:}\nabla\FF
=\Big(\frac{\upkappa'(\FF)}{\det\FF}-\frac{\upkappa(\FF){\rm Cof}\FF}{\det\FF^2}\Big){:}\nabla\FF=\frac{\upkappa'(\FF)-\upkappa(\FF)\FF^{-\top}\!\!\!\!}{\det\FF}\,{:}\nabla\FF,\!\!
\label{calculus-flexo-force}
\end{align}
where we again used the algebra $F^{-\top}={\rm Cof}F/\!\det F$ and the
calculus $\det'={\rm Cof}$. Thus \eq{test-damage} can be written as:
\begin{align}\nonumber
&\int_\varOmega\!{\rm div}\Big(\frac{\upkappa(\FF)\nabla\mm}{\det\FF}\Big)
{\cdot}\MU(\vv{\cdot}\Nabla)\mm\,\d\xx
=\!\int_\varOmega\!\bigg(
\MU|\Nabla\mm|^2\frac{\upkappa'(\FF){-}\upkappa(\FF)\FF^{-\top}\!\!\!}{2\det\FF}{:}(\vv{\cdot}\nabla)\FF
\\[-.1em]&\hspace*{17em}
+\frac{\upkappa(\FF)|\Nabla\mm|^2}{2\det\FF}{\rm div}\,\vv-\KK{:}\EE(\vv)\bigg)
\,\d\xx
\,,
\label{test-damage+}\end{align}
Similarly, this exchange-energy term ${\rm div}(\upkappa(\FF)\nabla\mm/\!\det\FF)$
tested by $\MU\pdt{}\mm$ is to be handled by using Green's formula once:
\begin{align}\nonumber
&\int_\varOmega\!{\rm div}\Big(\frac{\upkappa(\FF)\nabla\mm}{\det\FF}\Big)
{\cdot}\MU\pdt\mm\,\d\xx
=\MU\!\int_\varGamma\frac{\upkappa(\FF)(\nn{\cdot}\Nabla)\mm}{\det\FF}
{\cdot}\pdt\mm\,\d S-\MU\!\int_\varOmega\frac{\upkappa(\FF)\nabla\mm}{\det\FF}{:}
\nabla\pdt\mm\,\d\xx
\\[-.1em]&\nonumber\hspace*{5em}
=-\frac{\d}{\d t}\!\int_\varOmega\!\frac{\MU\upkappa(\FF)|\Nabla\mm|^2\!\!}{2\det\FF}\,\d\xx
+\int_\varOmega\!\frac{\MU|\Nabla\mm|^2}2
\Big(\frac{\upkappa(\FF)}{\det\FF}\Big)'{:}\pdt\FF\,\d\xx
\\&\hspace*{5em}
=-\frac{\d}{\d t}\!\int_\varOmega\!\frac{\MU\upkappa(\FF)|\Nabla\mm|^2\!\!}{2\det\FF}\,\d\xx
+\int_\varOmega\!\MU|\Nabla\mm|^2\frac{\upkappa'(\FF){-}\upkappa(\FF)\FF^{-\top}\!\!\!}{2\det\FF}\,{:}\pdt\FF\,\d\xx\,,
\label{test-damage++}\end{align}
where we again used $(\nn{\cdot}\Nabla)\mm=\bm0$ on $\varGamma$
and, where we again used, as in \eq{calculus-flexo-force}, that
$(\upkappa(\FF)/\!\det\FF)'=(\upkappa'(\FF){-}\upkappa(\FF)\FF^{-\top})/\!\det\FF$.
To merge
\eq{test-damage+} and \eq{test-damage++}, we use \eq{Euler-thermodynam2} and
also the calculus
\begin{align}\nonumber
&\MU|\Nabla\mm|^2\frac{\upkappa'(\FF){-}\upkappa(\FF)\FF^{-\top}\!\!\!}{2\det\FF}\,{:}\pdt\FF+\MU|\Nabla\mm|^2\frac{\upkappa'(\FF){-}\upkappa(\FF)\FF^{-\top}\!\!\!\!\!}{2\det\FF}{:}(\vv{\cdot}\nabla)\FF
\\\nonumber
&=\MU|\Nabla\mm|^2\frac{\upkappa'(\FF){-}\upkappa(\FF)\FF^{-\top}\!\!\!\!\!}{2\det\FF}\,{:}
(\nabla\vv)\FF
=\MU|\Nabla\mm|^2\frac{\upkappa'(\FF)\FF^{\top}\!{-}\upkappa(\FF)\bbI}{2\det\FF}\,{:}\ee(\vv)
\\&=\MU|\Nabla\mm|^2\frac{\upkappa'(\FF)\FF^{\top}\!}{2\det\FF}\,{:}\ee(\vv)
-\MU\frac{|\Nabla\mm|^2\upkappa(\FF)}{2\det\FF}{\rm div}\,\vv
\label{exchange-stress}\end{align}
when the frame indifference of $\upkappa$ is assumed.
Noticing that the last term in \eq{exchange-stress} cancels with
the same pressure term in \eq{test-damage+}, we obtain
\begin{align}\label{test-by-DTm}
&\int_\varOmega\!{\rm div}\Big(\frac{\upkappa(\FF)\nabla\mm}{\det\FF}\Big)
{\cdot}\MU\DT\mm\,\d\xx
=\int_\varOmega\!
\Big(\MU\frac{\upkappa'(\FF)|\Nabla\mm|^2\FF^{\top}\!\!\!
}{2\det\FF}-\KK\Big){:}\ee(\vv)\,\d x
-\frac{\d}{\d t}\!\int_\varOmega\!\frac{\MU\upkappa(\FF)|\Nabla\mm|^2\!\!}{2\det\FF}\,\d\xx\,.
\end{align}

Moreover, we use the Green theorem also for the driving magnetic field
$\hh=\hh_\text{\rm ext}{-}\hh_\text{\rm dem}$ with the demagnetizing field
$\hh_\text{\rm dem}=-\nabla u$:
\begin{align}\nonumber
&\int_\varOmega\mu_0\hh{\cdot}\DT\mm\,\d\xx
=\int_\varOmega\mu_0\hh_{\rm ext}{\cdot}\pdt\mm
+\mu_0\hh{\cdot}(\vv{\cdot}\Nabla)\mm
-\mu_0\hh_{\rm dem}{\cdot}\pdt\mm\,\d\xx
\\[-.2em]&\ \
=\int_\varOmega\pdt{}\big(\mu_0\hh_{\rm ext}{\cdot}\mm\big)
-\mu_0\frac{\pl\hh_{\rm ext}}{\pl t}{\cdot}\mm
-\mu_0(\Nabla\hh)^\top\mm{\cdot}\vv
+\mu_0\nabla(\hh{\cdot}\mm)\,\vv
-\mu_0\hh_{\rm dem}{\cdot}\pdt\mm\,\d\xx\,.
\nonumber\end{align}

Altogether, this is used to handle the
right-hand side of \eqref{Euler-thermodynam4} tested by $\ZJ\mm$:
\begin{align}\nonumber
&\int_\varOmega\!\Big(\mu_0\hh
-\frac{\uppsi_\mm'(\FF,\mm,\theta)}{\det\FF}+{\rm div}\Big(\frac{\upkappa(\FF)\nabla\mm}{\det\FF}\Big)\Big){\cdot}\ZJ\mm\,\d\xx
\\[-.1em]&\nonumber\hspace*{.7em}=\nonumber\!\int_\varOmega
\!\Big(\mu_0\hh
-\frac{\uppsi_\mm'(\FF,\mm,\theta)}{\det\FF}+{\rm div}\Big(\frac{\upkappa(\FF)\nabla\mm}
{\det\FF}\Big)\Big){\cdot}\DT\mm-(\mu_0\hh{-}\tt){\cdot}{\rm skw}(\Nabla\vv){\mm}\,\d\xx
\\[-.0em]&\nonumber\hspace*{.7em}
=\nonumber
\frac{\d}{\d t}\int_\varOmega\!\mu_0\hh_\text{\rm ext}{\cdot}\mm
-\frac{\upkappa(\FF)|\Nabla\mm|^2}{2\det\FF}\,\d\xx
+\!\int_\varOmega\!\bigg(\!
\Big(\MU\frac{\upkappa'(\FF)|\Nabla\mm|^2\FF^{\top}\!\!\!}{2\det\FF}-\KK\Big){:}\EE(\vv)
-\mu_0\pdt{\hh_\text{\rm ext}}{\cdot}\mm
\\[-.0em]&\nonumber\hspace*{1.7em}
-\Big(\frac{\upvarphi_\mm'(\FF,\mm)}{\det\FF}
+\frac{\COUPLING_\mm'(\FF,\mm,\theta)}{\det\FF}\Big){\cdot}\DT\mm
-\big(\hspace{-.7em}\lineunder{(\mu_0(\Nabla\hh)^\top\mm
-\mu_0\nabla(\hh{\cdot}\mm)}{magnetic force in \eq{Euler-thermodynam1}}\hspace{-.7em}\big){\cdot}\vv
-\mu_0\hh_{\rm dem}{\cdot}\pdt\mm
\\[-.2em]&\hspace*{2em}
-\hspace*{-.7em}\lineunder{{\rm skw}\big((\mu_0\hh{-}\uppsi_\mm'(\FF,\mm,\theta)){\otimes}\mm\big)}{$=\SS$ from \eqref{Euler-thermodynam1}}\hspace*{-.7em}{:}\Nabla\vv
+\hspace*{-.7em}\lineunder{\upkappa(\FF)\frac{\mm{\otimes}\nabla\mm{-}\nabla\mm{\otimes}\mm}{2\det\FF}}{$=\mathscr{S}$
from \eq{Euler-thermodynam1}}\hspace*{-.7em}
\Vdots\nabla^2\vv\!\bigg)\,\d\xx\,
\label{formula3}
\end{align}
with $\tt$ and $\KK$ from \eq{Euler-thermodynam1}.
Here the algebra $(\mu_0\hh{-}\uppsi_\mm'){\cdot}{\rm skw}(\nabla\vv)\mm
={\rm skw}((\mu_0\hh{-}\uppsi_\mm'){\otimes}\mm){:}\nabla\vv$ has been used.
Beside, to see the magnetic hyperstress $\mathscr{S}$, we used the calculus
\begin{align}\nonumber
\int_\varOmega\!
&{\rm div}\Big(\frac{\upkappa(\FF)}{\det\FF}\nabla\mm\Big){\cdot}
{\rm skw}(\Nabla\vv){\mm}\,\d\xx
\\&\nonumber=\int_\varGamma\,
\frac{\upkappa(\FF)}{\det\FF}
\hspace*{-.7em}\lineunder{(\nn{\cdot}\nabla)\mm}{$=\bm0$}\hspace*{-.7em}
\cdot{\rm skw}(\Nabla\vv){\mm}\,\d S-\int_\varOmega\,\frac{\upkappa(\FF)}{\det\FF}\nabla\mm
{:}\nabla\big({\rm skw}(\Nabla\vv){\mm}\big)\,\d\xx
\\&=-\int_\varOmega\,\frac{\upkappa(\FF)}{\det\FF}
\hspace*{-.7em}\lineunder{(\nabla\mm{\otimes}\nabla\mm){:}{\rm skw}(\Nabla\vv){\mm}}{$=\bm0$}\hspace*{-.7em}
+\upkappa(\FF)\frac{\nabla\mm{\otimes}\mm{-}\mm{\otimes}\nabla\mm}{2\det\FF}
\Vdots\nabla^2\vv\,\d\xx\,.
\end{align}

It remains to exploit \eq{Euler-thermodynam5}. Testing it by $\mu_0\pdt{}u$, we
use the calculus, including the Green theorem for the convective term, to obtain
\begin{align}\nonumber
0&=\int_{\R^d}\mu_0{\rm div}\big(\chi_\varOmega^{}\mm-\Nabla u\big)\pdt{u}\,\d\xx=
\int_{\R^d}\mu_0\big(\Nabla u-\chi_\varOmega^{}\mm\big){\cdot}\Nabla\pdt{u}\,\d\xx
\\\nonumber
&=\frac{\d}{\d t}\int_{\R^d}\frac{\mu_0}2|\Nabla u|^2\,\d\xx
+\int_\varOmega\mu_0\mm{{\cdot}}\pdt{\hh_{\rm dem}}\,\d\xx
\\
&=\frac{\d}{\d t}\bigg(\int_{\R^d}\frac{\mu_0}2|\Nabla u|^2
+\int_\varOmega\mu_0\hh_{\rm dem}{\cdot}\mm\,\d\xx\bigg)
-\int_\varOmega\mu_0\pdt{\mm}{\cdot}\hh_{\rm dem}\,\d\xx\,.
\label{calculus-PM4}\end{align}
This is to be subtracted from \eq{formula3}, giving
cancellation of the terms $\pm\mu_0\pdt{}\mm{\cdot}\hh_{\rm dem}$
and merging $\mu_0(\hh_{\rm ext}{\cdot}\mm)-\mu_0(\hh_{\rm dem}{\cdot}\mm)
=\mu_0(\hh{\cdot}\mm)$.

Summarizing the above calculations, we formulated: 

\begin{proposition}[{\sl Magneto-mechanical dissipation energy balance}] Any
smooth solution of the system (\ref{Euler-thermodynam}a--e) with the boundary
conditions \eq{Euler-thermodynam-BC} satisfies the identity
\begin{align}\nonumber
  &\frac{\d}{\d t}\bigg(
  \int_\varOmega\!\!\!\!\linesunder{\frac\varrho2|\vv|^2}{kinetic}{energy}\!\!\!\!+\!\!\!\!\linesunder{\frac{\upvarphi(\FF,\mm)}{\det\FF}}{stored}{energy}\!\!\!\!+\!\!\!\!\linesunder{\frac{\upkappa(\FF)|\Nabla\mm|^2}{2\det\FF}}{exchange}{energy}\!\!\!\!
    -\!\!\!\!\linesunder{\mu_0\hh
    {\cdot}\mm_{_{_{_{_{_{}}}}}}\!\!\!}{Zeeman}{energy}\!\!\,\d\xx+\int_{\R^d}\!\!\!\!\!\!\!\!\!\!\linesunder{\frac{\mu_0}2|\Nabla u|^2_{_{_{_{_{_{}}}}}}\!\!\!}{energy of de-}{magnetizing field}\!\!\!\!\!\!\!\,\d\xx\bigg)
\\&\nonumber\hspace*{2em}
+\int_\varOmega\!\!\!\!\linesunder{\xi\big(\FF,\theta;\EE(\vv),\ZJ\mm\big)
    _{_{_{_{}}}}\!\!}{dissipation rate}{from \eqref{Euler-thermodynam3}}\!\!\!
    \d\xx
    +\!\int_\varGamma\!\!\!\!\!\!\!\!\!\!\linesunder{\nu_\flat|\vv|^p_{_{_{}}}\!}{boundary}{dissipation rate}\!\!\!\!\!\!\!\!\d S
=\int_\varOmega\!\bigg(\!\!\!\!\!\!\linesunder{\varrho\,\GRAVITY{\cdot}\vv_{_{_{_{}}}}\!\!}{power of}{gravity field}\!\!\!\!\!\!
-\!\!\!\!\linesunder{\mu_0\frac{\!\partial \hh_\text{\rm ext}\!}{\partial t}{\cdot}\mm
  }{power of}{external field}\!\!\!\!
  \\&
  \hspace*{7em}\lineunder{-\,
  \frac{\COUPLING_{\Fe}'(\Fe,\mm,\theta)\Fe^\top\!\!}{\det\Fe}{:}\ee(\vv)
-\frac{\COUPLING_\mm'(\FF,\mm,\theta)}{\det\FF}{\cdot}\DT\mm
\!}{adiabatic effects}\!\!\!\bigg)\d\xx+\!\int_\varGamma\!\!\!\!\!\!\linesunder{\!\!\TRACTION{\cdot}\vv_{_{_{_{}}}}\!\!}{power of}{traction}\!\!\!\!\!\d S\,.\hspace*{-1em}
\label{energy+}
\end{align}
\end{proposition}

When we add \eq{Euler-thermodynam3} tested by 1, the adiabatic and the
dissipative heat sources cancel with those in
\eq{energy+}. Thus we obtain (at least formally):

\begin{proposition}[{\sl Total energy balance}] Any
smooth solution of the evolution boundary-value problem
\eq{Euler-thermodynam}--\eq{Euler-thermodynam-BC} satisfies the identity
\begin{align}\nonumber
  &\frac{\d}{\d t}\bigg(\int_\varOmega\!\!\!\!
  \linesunder{\frac\varrho2|\vv|^2}{kinetic}{energy}\!\!\!\!+\!\!\!\!\!
  \linesunder{\frac{\upvarphi(\Fe,\mm)}{\det\Fe}}{stored}{energy}\!\!\!\!+\!\!\!\!\linesunder{\frac{\upkappa(\FF)|\Nabla\mm|^2}{2\det\FF}}{exchange}{energy}\!\!\!\!
  -\!\!\!\!\linesunder{\mu_0\hh
  {\cdot}\mm_{_{_{_{_{_{}}}}}}\!\!\!}{Zeeman}{energy}\!\!\,\d\xx+\!\!\!\!\!\linesunder{\OMEGA(\FF,\mm,\theta)}{heat}{energy}\!\!\d\xx
\\[-.2em]&\hspace{.1em}
+\int_{\R^d}\hspace{-1,9em}\linesunder{\frac{\mu_0}2|\Nabla u|^2_{_{_{_{_{_{}}}}}}\!\!\!}{energy of de-}{magnetizing field}\hspace{-1,3em}\d\xx\bigg)+\!\int_\varGamma\!\!\!\!\!\!\!\!\!\!\!\!\morelinesunder{\frac{\nu_\flat}{2}|\vv|^p_{_{_{}}}\!}{boundary}{heat production}{outflow}\!\!\!\!\!\!\!\!\!\!\d S
=\int_\varOmega\!\!\!\!\!\morelinesunder{\varrho\GRAVITY{\cdot}\vv_{_{_{_{}}}}}{power of}{gravity}{field}\!\!\!\!\!\!-\!\!\!\!\linesunder{\mu_0\frac{\!\partial \hh_\text{\rm ext}\!}{\partial t}{\cdot}\mm
  }{power of}{external field}\!\!\!\!
\d\xx
+\!\int_\varGamma\!\!\!\!\!\!\linesunder{\!\!\TRACTION{\cdot}\vv\!\!_{_{_{}}}}{power of}{traction\ }\!\!\!\!\!\!+\!\!\!\!\!\linesunder{h(\theta)}{heat}{flux}\!\!\!\d S\,.
\nonumber\\[-1.5em]
\label{thermodynamic-Euler-engr}\end{align}
\end{proposition}

Another aspect important both thermodynamically and also for mathematical
analysis is non-negativity of temperature, related with the 3rd law of
thermodynamics. This will be demonstrated later when we will exploit some
information about the quality of the velocity field extracted from
\eq{energy+}, cf.\ \eq{Euler-thermo-test-nonnegative} below.

\begin{remark}[{\sl Exchange hyper-stress}]\upshape
In principle, to balance the energetics, the magnetic exchange driving force
${\rm div}(\upkappa(\FF)\nabla\mm/\!\det\FF)$
in \eq{driving-field} may contribute either directly
to the skew-symmetric magnetic stress $\SS$ by
${\rm skw}({\rm div}(\frac{\upkappa(\FF)\nabla\mm}{\det\FF}){\otimes}\mm)$, as considered in
\cite{Roub22TMPE}, or to the  skew-symmetric hyperstress as $\mathscr{S}$.
Physically it is rather questionable which option is more relevant. The former
case would bring analytical troubles in the argumentation \eq{strong-hyper+}
below due to lack of compactness of $\nabla\vv$ as $\pdt{}\vv$ is not
estimated, in contrast to \cite{Roub22TMPE} where the inertial
force was handled in a simplified ``semi-compressible'' way. This
have led us to adopt the latter option here, which seems also more
physical and a similar skew-symmetric hyperstress
can be found in \cite{Toup64TECS}.
\end{remark}

\begin{remark}[{\sl Isotropic magnets}]\upshape
Let us note that, when $\uppsi(\FF,\cdot,\theta)$ is isotropic as in
Example~\ref{exa-neo-Hook}, the skew-symmetric magnetic stress
$\SS$ simplifies to $\mu_0{\rm skw}(\hh{\otimes}\mm)$ because
$\uppsi_\mm'(\FF,\mm,\theta)=k(\FF,\mm,\theta)\mm$ 
for some scalar-valued coefficient $k=k(\FF,\mm,\theta)$ so that
${\rm skw}(\uppsi_\mm'(\FF,\mm,\theta){\otimes}\mm)=
k(\FF,\mm,\theta){\rm skw}(\mm{\otimes}\mm)=\bm0$.
\end{remark}

\section{The analysis -- weak solutions of \eq{Euler-thermodynam}}
\label{sec-anal}

We will provide a proof of existence and certain regularity of weak
solutions. To this aim, the concept of multipolar viscosity is essential
but, anyhow, still quite nontrivial and carefully ordered arguments
will be needed. The peculiarities are that 
the inertial term in Eulerian setting involves varying mass density 
requiring sophisticated techniques from compressible fluid dynamics, 
the momentum equation is very geometrically nonlinear, and the heat equation 
has an $L^1$-structure with $\FF$-dependent heat capacity and
with the convective time derivative and ever-troubling adiabatic effects
due to necessarily general coupling of mechanical and thermal effect in
the deforming configuration in compressible media.

Usual analysis is made by some approximation, a-priori estimates,
and limit passage possibly in several steps. The mentioned strong
nonlinearity makes time discretization problematic. On the other hand,
the space discretization by a (conformal) Faedo-Galerkin method is
also not straightforward because of several ``nonlinear'' tests leading to
the basic energy balances in Section~\ref{sec-engr}, being confronted
in particular with the Lavrentiev phenomenon as occurring already in
static nonlinear elasticity \cite{Ball87SMSE,Ball02SOPE,FoHrMi03LGPN}.
Anyhow, careful suitably regularized ``semi-Galerkin'' discretization allowing estimation
of the magneto-mechanical part separately from the thermal part
and a successive limit passage  will work.

\subsection{Definition of weak solutions and the main results}
We will use the standard notation concerning the Lebesgue and the Sobolev
spaces, namely $L^p(\varOmega;\R^n)$ for Lebesgue measurable functions
$\varOmega\to\R^n$ whose Euclidean norm is integrable with $p$-power, and
$W^{k,p}(\varOmega;\R^n)$ for functions from $L^p(\varOmega;\R^n)$ whose
all derivative up to the order $k$ have their Euclidean norm integrable with
$p$-power. We also write briefly $H^k=W^{k,2}$. The notation
$p^*$ will denote the exponent from the embedding
$W^{1,p}(\varOmega)\subset L^{p^*}(\varOmega)$, i.e.\ $p^*=dp/(d{-}p)$
for $p<d$ while $p^*\ge1$ arbitrary for $p=d$ or $p^*=+\infty$ for $p>d$.
Moreover, for a Banach space
$X$ and for $I=[0,T]$, we will use the notation $L^p(I;X)$ for the Bochner
space of Bochner measurable functions $I\to X$ whose norm is in $L^p(I)$
while $W^{1,p}(I;X)$ stands for functions $I\to X$ whose distributional
derivative is in $L^p(I;X)$. Also, $C(\cdot)$ and $C^1(\cdot)$
will denote spaces of continuous and continuously differentiable functions.

Moreover, as usual, we will use $C$ for a generic constant which may vary
from estimate to estimate.

We will consider an initial-value problem, prescribing the initial conditions
\begin{align}\label{Euler-thermodynam-IC}
\varrho(0)=\varrho_0\,,
\ \ \ \ \vv(0)=\vv_0\,,\ \ \ \ 
\FF(0)=\FF_0\,,\ \ \ \ \mm(0)=\mm_0\,,\ \ 
\text{ and }\ \ \theta(0)=\theta_0\,;
\end{align}
here and in what follows, we will use the short-hand notation
as $[\varrho(t)](\xx)=\varrho(t,\xx)$.
Referring to the referential mass density $\rhoR$, the initial conditions
should satisfy $\varrho_0=\rhoR/\!\det\FF_0$.
To devise a weak formulation of the initial-boundary-value problem
\eq{Euler-thermodynam-BC} and \eq{Euler-thermodynam-IC} for
the system \eq{Euler-thermodynam}, we use the by-part integration in time and
the Green formula for the inertial force.
The nonsmoothness of ${\rm Dir}(\cdot)$ applied on $\ZJ\mm$ leads to a
variational inequality, arising by a standard definition of the
convex subdifferential of the convex potential of the monotone
set-valued mapping $\rr\mapsto\tau\rr+\HC(\FF,\theta){\rm Dir}(\rr)$,
let us denote it as $D(\FF,\theta;\rr)=\tau|\rr|^2/2+\HC(\FF,\theta)|\rr|$.
Then \eq{Euler-thermodynam4} has the form
$\pl_{\ZJ\mm} D(\FF,\theta;\ZJ\mm)\ni\mu_0\hh_{\rm ext}+\mu_0\nabla u
-\tt+\mm{\times}\ZJ\mm/\gamma(\FF,\mm,\theta)$ with $\tt$ from \eq{driving-field},
from which we obtain a variational inequality by taking into account the
standard definition of the (partial) convex subdifferential $\pl_{\ZJ\mm}$.
This involves $\tt{\cdot}\ZJ\mm$ which contains the product of 
${\rm div}(\upkappa_{\nabla\mm}'(\FF,\nabla\mm)/\!\det\FF)$ with
$\DT\mm$. This product would cause troubles in convergence of approximate
solutions, so we will better avoid it in the weak formulation by a 
substitution using \eq{test-by-DTm} integrated over $I$, i.e.
\begin{align}\nonumber
\int_0^T\!\!\!\int_\varOmega\!{\rm div}\Big(\frac{\upkappa(\FF)\nabla\mm\!}{\det\FF}\,\Big)
{\cdot}\DT\mm\,\d\xx\d t
&=\!\int_\varOmega\!
\frac{\upkappa(\FF_0)|\Nabla\mm_0|^2\!\!}{2\det\FF_0}-\frac{\upkappa(\FF(T))|\Nabla\mm(T)|^2}{2\det\FF(T)}
\,\d\xx
\\&\hspace*{1em}
+\int_0^T\!\!\!\int_\varOmega\!\Big(\MU\frac{\upkappa'(\FF)|\Nabla\mm|^2\FF^{\top}\!\!\!
}{2\det\FF}-\KK\Big){:}\ee(\vv)\,\d\xx\d t\,.\!
\end{align}
Also we use the orthogonality $(\mm{\times}\ZJ\mm){\cdot}\ZJ\mm=0$, which
eliminates this (otherwise not integrable) term and which altogether gives
the variational inequality \eq{Euler2-weak} below.

\begin{definition}[Weak solutions to \eq{Euler-thermodynam}]\label{def}
For $p\in[1,\infty)$, a six-tuple $(\varrho,\vv,\Fe,\mm,u,\theta)$
with $\varrho\in H^1(I{\times}\varOmega)$,
$\vv\in L^p(I;W^{2,p}(\varOmega;\R^d))$,
$\FF\in H^1(I{\times}\varOmega;\R^{d\times d})$,
$\mm\in H^1(I;L^2(\varOmega;\R^d))\cap L^\infty(I;H^1(\varOmega;\R^d))$,
$u\in L^\infty(I;H^1(\R^d))$, and
$\theta\in L^1(I;W^{1,1}(\varOmega))$ will be called a weak solution to the
system \eq{Euler-thermodynam} with the boundary conditions
\eq{Euler-thermodynam-BC} and the initial condition \eq{Euler-thermodynam-IC}
if 
\begin{subequations}\label{Euler-weak-}\begin{align}
&\frac{\uppsi_\Fe'(\Fe,\mm,\theta)\Fe^\top\!\!\!}{\det\Fe}\in L^1(I{\times}\varOmega;
\R_{\rm sym}^{d\times d})\,,
\\&\frac{\COUPLING'_{\Fe}(\Fe,\mm,\theta)\Fe^\top\!\!\!}{\det\Fe}\in
L^{q'}(I{\times}\varOmega;\R_{\rm sym}^{d\times d})\,,
\\&{\rm div}\Big(\frac{\upkappa(\FF)|\Nabla\mm|^2}{\det\FF}\Big)\in
L^2(I{\times}\varOmega;\R^d)
\end{align}\end{subequations}
with $\det\Fe>0$ a.e.\ such that the integral identities 
\begin{subequations}\label{Euler-weak}\begin{align}
&\nonumber\int_0^T\!\!\!\int_\varOmega\bigg(\Big(\frac{
\uppsi_\Fe'(\Fe,\mm,\theta)\Fe^\top\!\!\!}{\det\Fe}+
\frac{|\nabla\mm|^2\upkappa'(\FF)\Fe^\top\!\!\!}{2\det\Fe}
+\KK+\DIS|\ee(\vv)|^{p-2}\ee(\vv)-\varrho\vv{\otimes}\vv\Big){:}\ee(\widetilde\vv)
\\[-.1em]&\hspace{.1em}\nonumber
+\SS{:}{\rm skw}(\nabla\wt\vv)+\mu_0(\hh{\cdot}\mm){\rm div}\widetilde\vv
-\mu_0(\nabla\hh)^\top{:}(\mm{\otimes}\wt\vv)
+\big(\nu_2|\nabla^2\vv|^{p-2}\nabla^2\vv
{+}\mathscr{S}\big)\Vdots\nabla^2\widetilde\vv
\\[-.2em]&\hspace{.1em}
\label{Euler1-weak}
-\varrho\vv{\cdot}\pdt{\widetilde\vv}
\bigg)\,\d\xx\d t
=\!\int_0^T\!\!\!\!\int_\varOmega\varrho\GRAVITY{\cdot}\widetilde\vv
\,\d\xx\d t
+\!\int_0^T\!\!\!\!\int_\varGamma(\TRACTION{-}\nu_\flat|\vv|^{p-2}\vv){\cdot}\widetilde\vv
\,\d S\d t
+\!\int_\varOmega\!\varrho_0\vv_0{\cdot}\widetilde\vv(0)\,\d\xx
\intertext{with $\hh=\hh_{\rm ext}{+}\nabla u$, $\KK$, $\SS$, and $\mathscr{S}$
from \eq{Euler-thermodynam1} holds for any $\widetilde\vv$ smooth with
$\widetilde\vv{\cdot}\nn={\bm0}$ and $\widetilde\vv(T)=0$, and\COMMENT{CHECK SIGNS}}
\nonumber
&\int_0^T\!\!\!\int_\varOmega\!\bigg(\frac\tau2|\wt\rr|^2+\HC(\FF,\theta)|\wt\rr|
-{\rm div}\frac{\upkappa(\FF)\nabla\mm}{\MU\det\Fe}{\cdot}
{\rm skw}(\nabla\vv)\mm-\Big(\MU\frac{\upkappa'(\FF)|\Nabla\mm|^2\FF^{\top}\!\!\!
}{2\det\FF}-\KK\Big){:}\ee(\vv)
\\[-.3em]&\hspace*{3em}\nonumber
+\frac{\upkappa(\FF)\nabla\mm}{\det\Fe}{:}\nabla\wt\rr
-\Big(\mu_0\hh-\frac{\uppsi_\mm'(\Fe,\mm,\theta)\!}{\MU\det\Fe}\:\Big){\cdot}(\wt\rr-\ZJ\mm)
+\frac{\mm{\times}\ZJ\mm}{\gamma(\FF,\mm,\theta)}{\cdot}\wt\rr\bigg)\,\d\xx\d t
\\[-.1em]&\hspace*{1em}\label{Euler2-weak}\ge
\int_0^T\!\!\!\int_\varOmega\!\frac\tau2|\ZJ\mm|^2+\HC(\FF,\theta)|\ZJ\mm|\,\d\xx\d t
+\!\int_\varOmega\!
\frac{\upkappa(\FF(T))|\Nabla\mm(T)|^2}{2\MU\det\FF(T)}
-\frac{\upkappa(\FF_0)|\Nabla\mm_0|^2\!\!}{2\MU\det\FF_0}\,\d\xx
\intertext{for any $\wt\rr\in L^2(I;H^1(\varOmega;\R^d))$, and further}
&\int_{\R^3}\!\nabla u(t){\cdot}\nabla\wt u\,\d\xx
=\int_\varOmega\!\mm(t){\cdot}\nabla\wt u\,\d\xx
\intertext{holds for any $\wt u\in H^1(\R^d)$ and for a.a.\ $t\in I$, and}
\nonumber
&\int_0^T\!\!\!\int_\varOmega\!\bigg(\!\OMEGA(\FF,\mm,\theta)\pdt{\widetilde\theta}
+\big(\OMEGA(\FF,\mm,\theta)\vv{-}\COND(\FF,\theta)\nabla\theta\big){\cdot}\nabla\widetilde\theta
\\[-.4em]&\hspace{4em}\nonumber
+\Big(\xi(\FF,\theta;\EE(\vv),\nabla^2\vv,\ZJ\mm)
+\frac{\COUPLING'_{\Fe}(\Fe,\mm,\theta)\Fe^\top\!\!\!}{\det\Fe}\,
{:}\ee(\vv)
+\frac{\COUPLING_\mm'(\FF,\mm,\theta)}{\det\FF}{\cdot}\DT\mm\Big)
\,\widetilde\theta\bigg)\d\xx\d t
\\[-.4em]&\hspace{11em}
+\!\int_0^T\!\!\!\int_\varGamma\Big(h(\theta){+}\frac{\nu_\flat}2|\vv|^p\Big)\widetilde\theta\,\d S\d t
+\!\int_\varOmega\!\OMEGA(\FF_0,\theta_0)\widetilde\theta(0)\,\d\xx=0
\label{Euler3-weak}
\end{align}
\end{subequations}
with $\xi(\FF,\theta;\cdot,\cdot)$ and $\OMEGA(\cdot,\cdot,\cdot)$ from
\eq{Euler-thermodynam3}
holds for any $\widetilde\theta$ smooth with $\widetilde\theta(T)=0$, and
the equations \eq{Euler-thermodynam0} and \eq{Euler-thermodynam2} hold a.e.\
on $I{\times}\varOmega$ with $\vv(0)=\vv_0$ and $\FF(0)=\FF_0$
a.e.\ on $\varOmega$, and also $\mm(0)=\mm_0$ is to
hold a.e.\ on $\varOmega$.
\end{definition}

Before stating the main analytical result, let us summarize the data
qualification which will be fitted to the motivating Example~\ref{exa-neo-Hook}.
For some $\delta>0$ and $s>0$, we assume:
\begin{subequations}\label{Euler-ass}\begin{align}
&\varOmega\ \text{ a smooth bounded domain of $\R^d$, }\ d=2,3,
\\&\nonumber
\upvarphi\in C^1({\rm GL}^+(d){\times}\R^d),\ \forall\FF\,{\in}\,{\rm GL}^+(d){:}\
\ \ \ \upvarphi(\FF,\mm)\ge
\delta\big(1+|\mm|^s\det\FF\big),\COMMENT{\ CHECK}
\\&\hspace{1em}\exists\, C\in C({\rm GL}^+(d)),\ \forall\FF\,{\in}\,{\rm GL}^+(d),\ \mm\,{\in}\,\R^d:\ 
|\upvarphi_\mm'(\FF,\mm)|\le C(\FF)(1{+}|\mm|^{1+2^*/2}),
\label{Euler-ass-phi}
\\&\nonumber
\COUPLING\,{\in}\, C^2({\rm GL}^+(d){\times}\R^d{\times}\R^+),\ \forall(\FF,\mm,\theta)\,{\in}\,{\rm GL}^+(d){\times}\R^d{\times}\R^+:\ \
\COUPLING_{\theta\theta}''(\Fe,\mm,\theta)\le\frac{-\delta}{\det\FF}\,,
\\[-.2em]&\hspace{5em}
\Big|\frac{\COUPLING_{\Fe}'(\Fe,\mm,\theta)\Fe^\top\!\!}{\det\Fe}\;\Big|
+\Big|\frac{\COUPLING_\mm'(\FF,\mm,\theta)}{\det\FF}\Big|^2
\le C\Big(1{+}\frac{\upvarphi(\Fe,\mm){+}\theta}{\det\Fe}\Big),
\label{Euler-ass-adiab}
\\[-.0em]&\nonumber
\forall K{\subset}\,{\rm GL}^+(d)\text{ compact }\,\exists\,C_K<\infty
\ \forall(\FF,\mm,\theta)\,{\in}\,K{\times}\R^d{\times}\R^+{:}\ \  
|\OMEGA_\FF'(\FF,\mm,\theta)|\le  C_K(1{+}\theta),\\&\hspace{5em}
\OMEGA_\theta'(\FF,\mm,\theta)+|\OMEGA_\mm'(\FF,\mm,\theta)|+|\OMEGA_{\FF\theta}''(\FF,\mm,\theta)|+|\OMEGA_{\mm\theta}''(\FF,\mm,\theta)|\le C_K,
\label{Euler-ass-primitive-c}
\\[-.2em]&
\DIS>0,\ \nu_2>0,\ \nu_\flat>0,
\label{Euler-ass-xi}
\\&\nonumber
\gamma\,{\in}\, C({\rm GL}^+(d){\times}\R^d{\times}\R^+)\ \text{ positive},\ \
\forall K{\subset}\,{\rm GL}^+(d)\text{ compact }\
\\[-.7em]&\hspace{8em}\label{Euler-ass-gamma}
\exists\,C_K<\infty\ \ \forall(\FF,\mm,\theta)\,{\in}\,K{\times}\R^d{\times}\R^+{:}\ \ \  
\frac{|\mm|}{\gamma(\FF,\mm,\theta)}\le C_K\,,
\\[-.5em]&
\upkappa\in C^1({\rm GL}^+(d))\,,\ \ \
\mbox{$\inf_{\FF\in{\rm GL}^+(d)}^{}\upkappa(\FF)>0$}\,,
\label{Euler-ass-kappa}
\\&\COND\in C({\rm GL}^+(d){\times}\R^+)\ \text{ bounded},\ \ \
\mbox{$\inf_{\FF\in{\rm GL}^+(d),\theta\in\R^+}^{}\COND(\FF,\theta)>0$}\,,
\label{Euler-ass-cond}
\\&\nonumber
h:I{\times}\varGamma{\times}\R^+\to\R\ \text{ Carath\'eodory function},\ \ \
0\le\theta\,h(t,\xx,\theta)\le C(1{+}\theta^2)\ \ \text{ and}
\\[-.1em]&\hspace{11em}
h(t,\xx,\theta)\le h_{\max}(t,\xx)\ \ \text{ for some }\ h_{\max}\in L^1(I{\times}\varGamma)\,,
\label{Euler-ass-h}
\\&
{\bm g}\in L^1(I;L^\infty(\varOmega;\R^d))\,,
\ \ \hh_{\rm ext}\in W^{1,1}(I;L^{s'}(\varOmega;\R^d))\,,
\ \ \TRACTION\in L^2(I{\times}\varGamma;\R^d),\ \ \TRACTION{\cdot}\nn=0\,,
\label{Euler-ass-f-g}
\\&\vv_0\in L^2(\varOmega;\R^d)\,,\ \ \ 
\Fe_0\in W^{1,r}(\varOmega;\R^{d\times d})\,,\ \ r>d\,,\ \ \ \text{ with }\ \ \
{\rm min}_{\barOmega}^{}\det\Fe_0>0\,,
\label{Euler-ass-Fe0}
\\&\rhoR\in L^\infty(\varOmega)\cap W^{1,r}(\varOmega)\,,\ \ r>d\,,\ \
\text{ with }\ \ {\rm min}_{\barOmega}^{}\rhoR>0\,,
\label{Euler-ass-rhoR}
\\&\mm_0\in H^1(\varOmega;\R^d)\,,\ \ \ \ \theta_0\in L^1(\varOmega),\ \ \ \
\theta_0\ge0\ \text{ a.e.\ on }\ \varOmega\,,
\label{Euler-ass-theta0}
\end{align}\end{subequations}
where $\OMEGA$ in \eq{Euler-ass-primitive-c} is from \eq{Euler-thermodynam3}.
Let us note that the first condition in \eq{Euler-ass-adiab} is just a condition
on the heat capacity $c=c(\FF,\mm,\theta)=\OMEGA_\theta'(\Fe,\mm,\theta)$ and
implies the coercivity $\OMEGA(\Fe,\mm,\theta)\ge \delta\theta/\!\det\FF$ since
$\OMEGA(\Fe,\mm,0)=0$.
One should note that the referential stored energy (in contrast to the
actual stored energy) enters the model only through its derivatives and 
can be modified without loss of generality by adding a constant, so that
\eq{Euler-ass-phi} could be understood simply as coercivity
$\upvarphi(\FF,\mm)/\!\det\FF\ge\delta|\mm|^s$. Independently, the natural blow-up
under compression, i.e.\ $\upvarphi(\FF,\mm)\to\infty$ if $\det\FF\to0+$,
is allowed in \eq{Euler-ass-phi}. The condition \eq{Euler-ass-h} is well fitted with the
standard situation that the boundary flux is $h(\theta)=f(\theta_{\rm ext})-f(\theta)$
with an increasing function $f$ and with $\theta_{\rm ext}\ge0$ a prescribed
external temperature, so that
one can choose $h_{\max}=f(\theta_{\rm ext})$ provided we prove that
$\theta\ge0$. Also the condition $\theta\,h(t,\xx,\theta)\le C(1{+}\theta^2)$
is well compatible with this ansatz provided $f(0)=0$ and $f(\theta_{\rm ext})\in
L^2(I{\times}\varGamma)$. 

\begin{theorem}[Existence and regularity of weak solutions]\label{prop-Euler}
Let $p>d$ and $s\ge 2p/(p{-}2)$, and the assumptions \eq{frame-indifference} and
\eq{Euler-ass}
hold. Then:\\
\Item{(i)}
{there exist a weak solution $(\varrho,\vv,\FF,\mm,u,\theta)$ according
Definition~\ref{def} with a non-negative mass density
$\varrho\in L^\infty(I;W^{1,r}(\varOmega))$ such that
$\pdt{}\varrho\in L^\sigma(I;L^{r\sigma/(r+\sigma)}(\varOmega))$ with
$3\le\sigma<p(pd{+}4p{-}2d)/(4p{-}2d)$,\COMMENT{??? cf \eq{est-time-derivaitves}}
and a non-negative temperature $\theta\in L^\infty(I;L^1(\varOmega))
\,\cap\,L^\EXP(I;W^{1,\EXP}(\varOmega))$ with
$1\le\EXP<(d{+}2)/(d{+}1)$,
and further $\pdt{}{\Fe}\in L^{p}(I;L^r(\varOmega;\R^{d\times d}))$ and
$\Nabla\Fe\in L^\infty(I;L^r(\varOmega;\R^{d\times d\times d}))$,
and $\pdt{}{\mm}\in L^2(I{\times}\varOmega;\R^d)$ and
$\Nabla\mm\in L^\infty(I;L^2(\varOmega;\R^d))$,
and $\nabla^2\mm\in L^2(I{\times}\varOmega;\R^{d\times d\times d})$.}
\Item{(ii)}
{Moreover, this solution complies with energetics in the sense that
the energy dissipation balance \eq{energy+}
as well as the total energy balance \eq{thermodynamic-Euler-engr}
integrated over time interval
$[0,t]$ with the initial conditions \eq{Euler-thermodynam-IC} hold.}
\end{theorem}

\subsection{Some auxiliary results and formal a-priori estimates}

Let us first formulate two auxiliary assertions:

\begin{lemma}[See \cite{Roub22TVSE}]\label{lem1}
Given $\vv\in L^1(I;W^{2,p}(\varOmega;\R^d))$ with $p>d$ and $\vv{\cdot}\nn=0$ on
$I{\times}\varGamma$ and $\varrho_0\in W^{1,r}(\varOmega)$,
\eq{Euler-thermodynam0} has a unique weak solution
$\varrho\in C_{\rm w}(I;W^{1,r}(\varOmega))\cap W^{1,1}(I;L^r(\varOmega))$
which satisfies it a.e.\ on $I\times\varOmega$ and also the estimate holds:
\begin{align}\label{F-evol-est}
\|\varrho\|_{L^\infty(I;W^{1,r}(\varOmega))\,\cap\,W^{1,1}(I;L^r(\varOmega))}^{}\le 
\mathfrak{C}\Big(\|\Nabla\vv\|_{L^1(I;W^{1,p}(\varOmega;\R^{d\times d}))}^{}\,,\,
\|\varrho_0\|_{W^{1,r}(\varOmega)}^{}\Big)
\end{align}
holds with some $\mathfrak{C}\in C(\R^2)$. Moreover, the mapping
\begin{align}\label{v-mapsto-F}
\vv\mapsto\varrho:L^1(I;W^{2,p}(\varOmega;\R^d))\to
L^\infty(I;W^{1,r}(\varOmega))
\end{align}
is (weak,weak*)-continuous. The analogous assertion holds for \eq{DT-det-1}, assuming
$1/\det\FF_0\in W^{1,r}(\varOmega)$, and for \eq{cont-eq-inverse}, assuming
$1/\varrho_0\in W^{1,r}(\varOmega)$. Eventually, it holds $\R^{d\times d}$-valued also for
\eq{Euler-thermodynam2}, assuming $\FF_0\in W^{1,r}(\varOmega;\R^{d\times d})$.
\end{lemma}

For the approximation method in the proof below, 
we will still need a modification of Lemma~\ref{lem1} for a non-homogeneous
evolution-and-transport equation  \eq{Euler-thermo-reg5}, whose proof is
a straightforward modification (partly simplification) of \cite[Sect.4]{Roub22TVSE}:

\begin{lemma}\label{lem2}
Given $\vv\in L^p(I;W^{2,p}(\varOmega;\R^d))$ with $p>d$ and
$\rr\in L^1(I;L^2(\varOmega;\R^d))$, the equation
$\pdt{}\mm+(\vv{\cdot}\nabla)\mm-{\rm skw}(\nabla\vv)\mm=\rr$ with the initial
condition $\mm_0\in L^2(\varOmega;\R^d)$ and 
has a unique weak solution $\mm\in C_{\rm w}(I;L^2(\varOmega;\R^d))\cap W^{1,1}(I;H^1(\varOmega;\R^d)^*)$ and also the estimate holds:
\begin{align}\label{m-evol-est}
\|\mm\|_{L^\infty(I;L^2(\varOmega))\,\cap\,W^{1,1}(I;H^1(\varOmega;\R^d)^*)}^{}\le 
\mathfrak{C}\Big(\|\Nabla\vv\|_{L^1(I;W^{1,p}(\varOmega;\R^{d\times d}))}^{}\,,\,
\|\mm_0\|_{L^2(\varOmega;\R^d)}^{}\,,\,\|\rr\|_{L^1(I;L^2(\varOmega;\R^d))}^{}\Big)
\end{align}
holds with some $\mathfrak{C}\in C(\R^3)$. Moreover, the mapping
\begin{align}\label{v-mapsto-m}
(\vv,\rr)\mapsto\mm:L^1(I;W^{2,p}(\varOmega;\R^d))\times
 L^1(I;L^2(\varOmega;\R^d))\to L^\infty(I;L^2(\varOmega;\R^d))
\end{align}
is (weak,weak*)-continuous.
\end{lemma}

Formally, the assumptions \eq{Euler-ass} yield some a-priori bounds which can
be obtained from the total energy balance \eq{thermodynamic-Euler-engr} and the 
mechanical energy-dissipation balance \eq{energy+} for any sufficiently regular
solution $(\varrho,\vv,\FF,\mm,u,\theta)$ with $\theta\ge0$ a.e.\ in
$I{\times}\varOmega$.
Later, we will prove existence of such solutions, but unfortunately we are not
able to claim that every weak solution has $\theta$ non-negative. For the
approximation method used in the proof below, we assume the data $\uppsi$,
$\COND$, and $h$ to be defined also for the negative temperature by extending
them as
\begin{align}\nonumber
&\uppsi(\FF,\mm,\nabla\mm,\theta):=\upvarphi(\FF,\mm)
+\theta\big({\rm ln}(-\theta){-}1\big)+\upkappa(\FF)|\nabla\mm|^2/2\,,
\\&\COND(\FF,\theta):=\COND(\FF,-\theta)\,,
\ \ \ \text{ and }\ \ \ h(t,\xx,\theta):=h(t,\xx,-\theta)
\ \ \ \text{ for }\ \ \theta<0
\label{extension-negative}\end{align}
with $\upvarphi$ and $\COUPLING$ from the split \eq{ansatz}. This definition
makes $\uppsi:{\rm GL}^+(d){\times}\R^d{\times}\R^{d\times d}{\times}\R\to\R$
continuous and implies that $\OMEGA(\FF,\mm,\cdot)$ as well as
$\zeta_\FF(\FF,\mm,\cdot)$ and $\zeta_\mm(\FF,\mm,\cdot)$ continuous;
note that $\OMEGA(\FF,\mm,\theta)=\theta/\!\det\FF$,
$\zeta_\FF'(\FF,\mm,\theta)=\bm0$, and $\zeta_\mm'(\FF,\mm,\theta)=\bm0$ for $\theta$ negative.

First, we use  the total energy balance \eq{thermodynamic-Euler-engr}
integrated over a time interval $[0,t]$.
At this point, we must now assume (while being later proved at least for some
solution) that $\theta\ge0$, and similarly we now assume $\det\FF>0$.
In particular, we
have also $\OMEGA(\FF,\mm,\theta)\ge0$ and thus we are ``only'' to estimate the
right-hand side in \eq{thermodynamic-Euler-engr} together with the
Zeeman energy. For the bulk term  $\varrho\GRAVITY{\cdot}\vv$
and the boundary terms $\TRACTION{\cdot}\vv+h(\theta)$ we refer to 
\cite{Roub22TVSE}. The gravity force $\varrho\GRAVITY$ tested by
the velocity $\vv$ can be estimated by the H\"older/Young inequality as
\begin{align}\nonumber
\int_\varOmega\varrho\GRAVITY{\cdot}\vv\,\d \bm{x}
&=\int_\varOmega\sqrt{\frac{\rhoR}{\det\FF}}\sqrt{\varrho}\vv{\cdot}\GRAVITY\,\d \bm{x}
\le\Big\|\sqrt{\frac{\rhoR}{\det\FF}}\Big\|_{L^{2}(\varOmega)}
\big\|\sqrt{\varrho}\vv\big\|_{L^2(\varOmega;\R^d)}^{}
\big\|\GRAVITY\big\|_{L^\infty(\varOmega;\R^d)}^{}
\\&\nonumber\le
\frac12\bigg(\Big\|\sqrt{\frac{\rhoR}{\det\FF}}\Big\|_{L^2(\varOmega)}^2
\!+\big\|\sqrt{\varrho}\vv\big\|_{L^2(\varOmega;\R^d)}^2\bigg)
\,\big\|\GRAVITY\big\|_{L^\infty(\varOmega;\R^d)}^{}
\\&\nonumber=\big\|\GRAVITY\big\|_{L^\infty(\varOmega;\R^d)}^{}\int_\varOmega
\frac{\rhoR}{2\det\FF\!}+\frac\varrho2|\vv|^2\,\d\xx
\\&
\le\big\|\GRAVITY\big\|_{L^\infty(\varOmega;\R^d)}^{}
\bigg(\frac{\max\rhoR(\barOmega)\!\!}{2
\inf\varphi({\rm GL}^+(d){\times}\R^d)}\!\int_\varOmega\frac{\upvarphi(\FF,\mm)}{\det\FF\!}\,\d\xx
+\int_\varOmega\frac\varrho2|\vv|^2\,\d\xx\bigg)\,.
\label{Euler-est-of-rhs}\end{align}
The integral on the right-hand side of \eq{Euler-est-of-rhs} can then be
treated by the Gronwall lemma. In order to apply the Gronwall lemma one needs
the qualification \eq{Euler-ass-f-g} for $\GRAVITY$.
The boundary terms in \eq{thermodynamic-Euler-engr} can be
estimated, at current time instant $t\in I$, as
\begin{align}\label{h.v-est}
\int_\varGamma\!\TRACTION{\cdot}\vv+h(\theta)\,\d S\le
\Big(\frac2{\nu_\flat}\Big)^{1/(p-1)}\|\TRACTION\|_{L^{p'}(\varGamma;\R^d)}^{p'}
+\frac{\nu_\flat}{2p}\|\vv\|_{L^p(\varGamma;\R^d)}^p+\|h_{\rm max}\|_{L^1(\varGamma)}^{}
\end{align}
with the term  $\nu_\flat\|\vv\|_{L^p(\varGamma;\R^d)}^p/(2p)$ to be absorbed
in the left hand side of \eq{thermodynamic-Euler-engr}; here we used
the modelling assumption that part of the heat produced by the
boundary viscosity leaves the system, otherwise we would have to
confine ourselves to $\ff=0$. Then, integrating 
\eq{thermodynamic-Euler-engr} in time,
one can use the qualification of $\TRACTION$ in \eq{Euler-ass-f-g}.
The Zeeman energy $\mu_0\hh_{\rm ext}{\cdot}\mm$ and the power of external
field $\mu_0\frac{\pl}{\pl t}\hh_{\rm ext}{\cdot}\mm$ are to be estimated as
\begin{align}\nonumber
&\int_\varOmega\mu_0\hh_{\rm ext}(t){\cdot}\mm(t)-\hh_{\rm ext}(0){\cdot}\mm_0\,\d\xx
-\int_0^t\!\!\int_\varOmega\mu_0\frac{\pl}{\pl t}\hh_{\rm ext}{\cdot}\mm\,\d\xx\d t
\\&\qquad\nonumber\le C_{\mu_0,\delta}\|\hh_{\rm ext}(t)\|_{L^{s'}(\varOmega;\R^d)}^{s'}+
\frac\delta2\|\mm(t)\|_{L^s(\varOmega;\R^d)}^s
\\&\qquad\ \ \
+\mu_0\!\int_0^t\Big\|\pdt{\hh_{\rm ext}}\Big\|_{L^{s'}(\varOmega;\R^d)}
\big(1+\|\mm\|_{L^s(\varOmega;\R^d)}^s\big)\,\d t
+\mu_0\|\hh_{\rm ext}(0){\cdot}\mm_0\|_{L^1(\varOmega)}^{}
\label{Zeeman-est}\end{align}
with some $C_{\mu_0,\delta}$ depending on $\mu_0$ and $\delta$ chosen according to
the assumption \eq{Euler-ass-phi}. This assumption is 
then to be exploited for the stored energy on the left-hand side of
\eq{thermodynamic-Euler-engr} and, together with the
qualification \eq{Euler-ass-f-g} of $\hh_{\rm ext}$, used
for the Gronwall inequality.
As a result, since $\det\FF>0$, we obtain the (formal) a-priori estimates
\begin{subequations}\label{Euler-est}\begin{align}
&\big\|\sqrt\varrho\vv\big\|_{L^\infty(I;L^2(\varOmega;\R^d))}\le C\,,\ \ \
\label{est-rv2}
\\&\label{est-phi}
\Big\|\frac{\upvarphi(\Fe,\mm)}{\det\Fe}\Big\|_{L^\infty(I;L^1(\varOmega))}^{}\le C\,,
\\&\label{est-nabla-m}
\Big\|\frac{\upkappa(\FF)|\nabla\mm|^2}{\det\Fe}\Big\|_{L^\infty(I;L^1(\varOmega))}^{}\le C\,,
\\&\label{est-u}
\|\nabla u\|_{L^\infty(I;L^2(\R^d;\R^d))}^{}\le C\,,\ \ \text{ and}
\\&
\Big\|\frac\theta{\det\FF}\Big\|_{L^\infty(I;L^1(\varOmega))}^{}\le C\,.
\label{est-w/detF}
\intertext{From \eq{est-phi}, using \eq{Euler-ass-phi}, we also obtain}
&\label{est-m}
\big\|\mm\big\|_{L^\infty(I;L^s(\varOmega;\R^d))}^{}\le C\,.
\end{align}\end{subequations}

Now we come to \eq{energy+}; here we used the assumed  frame indifference of
$\COUPLING(\cdot,\mm,\theta)$ so that
$\COUPLING_{\Fe}'(\Fe,\mm,\theta)\Fe^\top$ is symmetric and thus
$\COUPLING_{\Fe}'(\Fe,\mm,\theta)\Fe^\top\!{:}\nabla\vv=
\COUPLING_{\Fe}'(\Fe,\mm,\theta)\Fe^\top\!{:}\ee(\vv)$. 
The issue is now estimation of the adiabatic term in
\eq{energy+}. Then, while relying on $p>d$, using the embedding 
\begin{align}\label{Korn+}
\|\nabla\vv\|_{L^\infty(\varOmega;\R^{d\times d})}^p\le
K\Big(\big\|\nabla^2\vv\big\|_{L^p(\varOmega;\R^{d\times d\times d})}^p
+\big\|\vv|_\varGamma\big\|_{L^p(\varGamma;\R^d)}^p\Big)\,,
\end{align}
we estimate 
\begin{align}\nonumber
&\!\!\int_\varOmega
\Big|\frac{\!\COUPLING_{\Fe}'(\Fe,\mm,\theta)\Fe^\top\!\!\!}{\det\Fe}{:}\ee(\vv)
+\frac{\!\COUPLING_\mm'(\Fe,\mm,\theta)\!}{\det\Fe}{\cdot}\DT\mm
\Big|\,\d\xx
\\&\nonumber\le
\Big\|\frac{\!\COUPLING_{\Fe}'(\Fe,\mm,\theta)\Fe^\top\!\!\!}{\det\Fe}\
\Big\|_{L^1(\varOmega;\R^{d\times d})}\|\ee(\vv)\|_{L^\infty(\varOmega;\R^{d\times d})}
+\Big\|\frac{\!\COUPLING_\mm'(\Fe,\mm,\theta)}{\det\Fe}
\Big\|_{L^2(\varOmega;\R^d)}\|\DT\mm\|_{L^2(\varOmega;\R^d)}
\\&\nonumber
\!\!\!\stackrel{\eq{Euler-ass-adiab}}\le\! C\Big\|1{+}\frac{\upvarphi(\Fe,\mm){+}\theta}{\det\Fe}
\Big\|_{L^1(\varOmega)}
\|\ee(\vv)\|_{L^\infty(\varOmega;\R^{d\times d})}
+C\Big\|1{+}\frac{\upvarphi(\Fe,\mm){+}\theta}{\det\Fe}
\Big\|_{L^1(\varOmega)}^{1/2}\|\DT\mm\|_{L^2(\varOmega;\R^d)}
\\&\nonumber
\!\!\!\!\stackrel{(\ref{Euler-est}b,e)}{\le}\!
\frac{C^{p'}K^{p'}}{\delta^{1/(p-1)}}
\Big\|1{+}\frac{\upvarphi(\Fe,\mm){+}\theta}{\det\Fe}\Big\|_{L^1(\varOmega)}^{p'}\!\!
+\delta\|\nabla^2\vv\|_{L^p(\varOmega;\R^{d\times d\times d})}^p\!
+\delta\big\|\vv|_\varGamma\big\|_{L^p(\varGamma;\R^d)}^p
\\&\nonumber\hspace*{4em}
+\frac{C^2}{\delta}
\Big\|1{+}\frac{\upvarphi(\Fe,\mm){+}\theta}{\det\Fe}
\Big\|_{L^1(\varOmega)}\!+\delta\|\ZJ\mm\|_{L^2(\varOmega;\R^d)}^2\!
+\delta\|{\rm skw}(\nabla\vv)\mm\|_{L^2(\varOmega;\R^d)}^2
\\&\nonumber
\le
\Big(\frac{C^{p'}K^{p'}}{\delta^{1/(p-1)}}{+}\frac{C^{2p}}{\delta^p}\Big)
\Big\|1{+}\frac{\upvarphi(\Fe,\mm){+}\theta}{\det\Fe}\Big\|_{L^1(\varOmega)}^{p'}\!\!
\!+C_{p,s,K}\big(1{+}\|\mm\|_{L^s(\varOmega;\R^d)}^s\big)
\\&\hspace*{10em}
+\delta\|\ZJ\mm\|_{L^2(\varOmega;\R^d)}^2\!
+2\delta\|\nabla^2\vv\|_{L^p(\varOmega;\R^{d\times d\times d})}^p\!
+2\delta\big\|\vv|_\varGamma\big\|_{L^p(\varGamma;\R^d)}^p\,.
\label{Euler-est-of-rhs+}\end{align}
The last inequality is due to the estimate, using \eq{Korn+}, that 
\begin{align}\nonumber
\|{\rm skw}(\nabla\vv)\mm\|_{L^2(\varOmega;\R^d)}^2
&\le C_{p,s,K}\big(1{+}\|\mm\|_{L^s(\varOmega;\R^d)}^s\big)+ 
K^{-1}\|\nabla\vv\|_{L^\infty(\varOmega;\R^{d\times d})}^p
\\
&\le C_{p,s,K}\big(1{+}\|\mm\|_{L^s(\varOmega;\R^d)}^s\big)
+\big\|\nabla^2\vv\big\|_{L^p(\varOmega;\R^{d\times d\times d})}^p\!
+\big\|\vv|_\varGamma\big\|_{L^p(\varGamma;\R^d)}^p
\label{Euler-est-of-rhs++}\end{align}
which holds if $1/p+1/s\le1/2$, as implied by the assumption $s\ge 2p/(p{-}2)$. 
Choosing $\delta>0$ sufficiently small, the last three terms in
\eq{Euler-est-of-rhs++} can be absorbed in the left-hand side of \eq{energy+}
while the other terms in \eq{Euler-est-of-rhs+} are already estimated by (\ref{Euler-est}b,e,f).
Thus we can exploit the dissipation rates in the magneto-mechanical energy
dissipation balance \eq{energy+} to obtain estimates 
\begin{align}
\|\ee(\vv)\|_{L^p(I{\times}\varOmega;\R^{d\times d})}^{}\le C\,,\ \ \ \ \
\|\nabla^2\vv\|_{L^p(I{\times}\varOmega;\R^{d\times d\times d})}^{}\le C\,,
\ \ \text{ and }\ \ \|\ZJ\mm\|_{L^2(I\times\varOmega;\R^d)}^{}\le C\,.
\label{est-e(v)}\end{align}

When assuming $p>d$, the estimate \eq{est-e(v)} is essential by preventing
evolution of singularities of the quantities transported by such a smooth
velocity field. Here, due to qualification of $\FF_0$ and 
$\varrho_0=\rhoR/\!\det\FF_0$ in \eq{Euler-ass-Fe0}  and \eq{Euler-ass-rhoR}
which implies also the qualification of the initial conditions
\begin{align}&\nonumber
\nabla\Big(\frac1{\det\FF_0}\Big)=-\frac{\det'(\FF_0){:}\nabla\FF_0}
{(\det\FF_0)^2}=-\frac{\FF_0^{-\top}{:}\nabla\FF_0}{\det\FF_0}
\in L^r(\varOmega;\R^d)\,,
\\&\nonumber
\nabla\varrho_0=\frac{\nabla\rhoR\ }{\det\FF_0}-\rhoR
\frac{\FF_0^{-\top}{:}\nabla\FF_0}{\det\FF_0}\in L^r(\varOmega;\R^d)\,,\ \ \text{ and }\ \ 
\\&\nonumber
\nabla\Big(\frac1{\varrho_0}\Big)=\nabla\Big(\frac{\det\FF_0}{\rhoR}\Big)=
\frac{{\rm Cof}\FF_0{:}\nabla\FF_0}{\rhoR}
-\frac{\det\FF_0\nabla\rhoR}{\rhoR^2}\in L^r(\varOmega;\R^d)\,,
\end{align}
we can see that Lemma~\ref{lem1} yields the estimates
\begin{subequations}\label{est+}
\begin{align}
\label{est+Fes}&\|\FF\|_{L^\infty(I;W^{1,r}(\varOmega;\R^{d\times d}))}\le C_r\,,
\ \ \ \:\Big\|\frac1{\det\FF}\Big\|_{L^\infty(I;W^{1,r}(\varOmega))}\le C_r\,,
\\&\label{est+rho}\|\varrho\|_{L^\infty(I;W^{1,r}(\varOmega))}^{}\le C_r\,,
\ \ \text{ and }\ \ \Big\|\frac1\varrho\Big\|_{L^\infty(I;W^{1,r}(\varOmega))}\!\le C_r
\ \ \ \text{ for any $1\le r<+\infty$}\,.
\intertext{From \eq{est-rv2} and \eq{est+rho} with $r>d$, we then have also}
&\|\vv\|_{L^\infty(I;L^2(\varOmega;\R^d))}^{}\le
\|\sqrt\varrho\vv\|_{L^\infty(I;L^2(\varOmega;\R^d))}^{}\Big\|\frac1{\sqrt\varrho}\Big\|_{L^\infty(I\times\varOmega)}^{}\le C\,,
\label{basic-est-of-v}
\intertext{and from \eq{est-nabla-m}, using \eq{Euler-ass-kappa}, we obtain}
&\|\nabla\mm\|_{L^\infty(I;L^2(\varOmega;\R^{d\times d}))}^{}\le
\Big\|\frac{\det\Fe}{\upkappa(\Fe)}\Big\|_{L^\infty(I\times\varOmega)}^{1/2}
\Big\|\frac{\upkappa(\Fe)|\nabla\mm|^2}
{\det\Fe}\Big\|_{L^\infty(I;L^1(\varOmega))}^{1/2}\!\le C\,.
\label{basic-est-of-m}
\end{align}\end{subequations}

Furthermore, having $\ZJ\mm$ estimated in \eq{est-e(v)} and 
by using the calculus ${\rm div}(\upkappa(\FF)\nabla\mm/\det\FF)$ $=
\upkappa(\FF)\Delta\mm/\det\FF)+
(\upkappa'(\FF)/\!\det\Fe-\upkappa(\FF){\rm Cof}\FF/\!\det\Fe^2){\Vdots}(\nabla\FF{\otimes}\nabla\mm)$, we can exploit \eq{Euler-thermodynam4} in the form
\begin{align}\nonumber
&\hspace{0em}\Delta\mm\in\frac{\det\FF}{\upkappa(\FF)}\bigg(
\tau\ZJ\mm+\HC(\FF,\theta){\rm Dir}(\ZJ\mm)
-\frac{\mm{\times}\ZJ\mm}{\gamma(\FF,\mm,\theta)}-\hh_{\rm ext}-\nabla
u
\\&\hspace{3em}
+\frac{\upvarphi_\mm'(\Fe,\mm)+\COUPLING_\mm'(\Fe,\mm,\theta)\!}{\det\Fe}
-\Big(\frac{\upkappa'(\FF)}{\det\Fe}-\frac{\upkappa(\FF){\rm Cof}\FF}{\det\Fe^2}\Big){\Vdots}(\nabla\FF{\otimes}\nabla\mm)\bigg)
\label{Delta-m-}
\end{align}
to estimate $\nabla^2\mm$ by the $H^2$-regularity of the
Laplacean with the homogeneous Neumann boundary conditions, as available
on smooth or convex domains. Here we use \eq{est+Fes} with $r>d$ and
\eq{basic-est-of-m}, we have $\nabla\FF{\otimes}\nabla\mm$ bounded in
$L^\infty(I;L^{2r/(r+2)}(\varOmega;\R^{d\times d\times d\times d}))$.
By \eq{Euler-ass-gamma} and by \eq{est-e(v)}, we can still see that
$\mm{\times}\ZJ\mm/\gamma(\FF,\mm,\theta)\in L^2(I{\times}\varOmega;\R^{d\times d})$.
 By \eq{est-u}, $\nabla u|_\varOmega^{} \in L^2(I{\times}\varOmega;\R^d)$. 
Moreover, by \eq{Euler-ass-phi}, $\upvarphi_\mm'(\FF,\mm)\in L^\infty(I;L^{2^*2/(2^*+2)}(\varOmega;\R^d))$
and by \eq{Euler-ass-adiab} $\COUPLING_\mm'(\FF,\mm,\theta)\in L^\infty(I;L^2(\varOmega;\R^d))$.
By comparison and the mentioned $H^2$-regularity, from \eq{Delta-m-}
we obtain
\begin{align}\label{Delta-m+}
&\|\nabla^2\mm\|_{L^2(I\times\varOmega;\R^{d\times d\times d})}^{}\le C\,.
\end{align}
Noting that the embedding $H^1(\varOmega)\subset L^2(\varOmega)$
is compact for $r>d$, \eq{Delta-m+} gives always a certain additional
information about $\nabla\mm$ in comparison with \eq{basic-est-of-m}.

\subsection{Proof of Theorem~\ref{prop-Euler}}

For clarity, we will divide the proof into ten steps. The inertial
term and the continuity equation \eq{Euler-thermodynam0} are treated
as in \cite{Roub22TVSE} and we thus sketch the proof in these aspects.\COMMENT{SOME MORE??}

Let us outline main technical difficulties: The time discretization
(Rothe's method) standardly needs convexity of $\upvarphi$ (which
is not a realistic assumption in finite-strain mechanics) possibly
weakened if there is some viscosity in $\FF$ (which is not directly
considered here, however). Also the conformal space discretization
(i.e.\ the Faedo-Galerkin method) is difficult since it cannot
directly copy the energetics because the ``nonlinear'' test
of \eqref{Euler-thermodynam2} by $[\upvarphi/\det]_\FF'(\FF)$ needed in
\eq{Euler-large-thermo} is problematic in this approximation
as $[\upvarphi/\!\det]_\FF'(\FF)$ is not in the respective finite-dimensional
space in general and similarly also the tests of \eq{Euler-thermodynam0}
by $|\vv|^2$ and of \eq{Euler-thermodynam4} by
$\ZJ\mm=\pdt{}\mm+(\vv{\cdot}\nabla)\mm-{\rm skw}(\nabla\vv)\mm$ are problematic.

\medskip\noindent{\it Step 1: a regularization}.
Referring to the formal estimates \eq{est+Fes}, we can choose $\LAM>0$ so
small that, for any possible sufficiently regular solution, it holds 
\begin{align}
&\det\FF>\LAM\ \ \ \ \text{ and }\ \ \ \ |\FF|<\frac1\LAM
\ \ \text{ a.e.\ on }\ I{\times}\varOmega\,.
\label{Euler-quasistatic-est-formal4}
\end{align}
We first regularize the stress $\TT$ and the other nonlinearities in
\eq{Euler-thermodynam} by considering a smooth cut-off
$\pi_\LAM\in C^1(\R^{d\times d})$ defined as
\begin{align}&\pi_\LAM(\FF)
:=\begin{cases}
\qquad\qquad1&\hspace{-8em}
\text{for $\det \FF\ge\LAM$ and $|\FF|\le1/\LAM$,}
\\
\qquad\qquad0&\hspace{-8em}\text{for $\det \FF\le\LAM/2$ or $|\FF|\ge2/\LAM$,}
\\
\displaystyle{\Big(\frac{3}{\LAM^2}\big(2\det\FF-\LAM\big)^2
-\frac{2}{\LAM^3}\big(2\det\FF-\LAM\big)^3\,\Big)\,\times}\!\!&
\\[.2em]
\qquad\qquad\displaystyle{\times\,\big(3(\LAM|\FF|-1)^2
-2(\LAM|\FF|-1)^3\,\big)}\!\!&\text{otherwise}.
\end{cases}
\label{cut-off-general}
\end{align}
Here $|\cdot|$ stands for the Frobenius norm $|\FF|=(\sum_{i,j=1}^d
F_{ij}^2)^{1/2}$ for $\FF=[F_{ij}]$, which makes $\pi_\LAM$ frame indifferent.
Thus we can regularize in a smooth way the singular nonlinearity
$1/\!\det(\cdot)$ and also extend $\upkappa(\cdot)/\!\det(\cdot)$ and $\COND$ and
also $\OMEGA$:
\begin{subequations}\begin{align}\label{cut-off-det}
&{\det}_\LAM(\FF):=\pi_\LAM(\FF)\det\FF+1-\pi_\LAM(\FF)\,,
\\&\label{cut-off-kappa}
\upkappa_\LAM(\FF):=\pi_\LAM(\FF)\upkappa(\FF)+(1{-}\pi_\LAM(\FF))\det\FF\,,
\ \ \text{ and }
\\&\label{cut-off-K}
\COND_\LAM(\FF,\theta):=\pi_\LAM(\FF)\COND(\FF,\theta)+1-\pi_\LAM(\FF)\,.
\end{align}\end{subequations}

Using the operator $\Delta^{-1}{\rm div}:L^2(\varOmega;\R^d)\to H^1(\R^d)$
defined by $u=[\Delta^{-1}{\rm div}](\mm)$ as a unique weak
solution to \eq{u-eq} with the ``boundary'' condition $u(\infty)=0$,
we can eliminate $u$; actually, this is rather the scenario for $d=3$ otherwise
it is formally possible because only $\nabla u$ but not $u$ itself occurs
in the system and its energetics. Altogether, for the above chosen $\LAM$ and
for any $\EPS>0$, we consider the regularized system 
\begin{subequations}\label{Euler-thermo-reg}
\begin{align}\label{Euler-thermo-reg0}
&\pdt\varrho=-\,{\rm div}(\varrho\vv)\,,
\\\nonumber
&\pdt{}(\varrho\vv)={\rm div}\Big(\TT_{\LAM,\EPS}(\FF,\mm,\theta)
{+}\KK_\LAM(\FF,\nabla\mm){+}\SS_{\LAM,\EPS}(\FF,\mm,\theta)
{+}\DIS|\ee(\vv)|^{p-2}\ee(\vv)
-\varrho\vv{\otimes}\vv
\\[-.6em]&\nonumber\hspace*{6em}
{-}{\rm div}\big(\mathscr{H}{+}\mathscr{S}_\LAM(\FF,\mm,\nabla\mm)\big)\Big)
+\mu_0(\nabla\hh)^\top\mm-\mu_0\nabla(\hh\cdot\mm)
  +\sqrt{\frac{\rhoR\varrho}{\!\det_\LAM(\FF)\!}\!}\ \GRAVITY
\\[.1em]\nonumber
    &\hspace*{1em}\text{ with }\
\TT_{\LAM,\EPS}(\FF,\mm,\theta)
=\Big(\frac{[\pi_\LAM\upvarphi]_\FF'(\FF,\mm)}{\det\FF}
+\frac{\pi_\LAM(\FF)\COUPLING_\FF'(\FF,\mm,\theta)}{(1{+}\EPS|\theta|
)\det\FF}
+\frac{|\nabla\mm|^2\upkappa_\LAM'(\FF)\!}{2\det\Fe}\:\Big)\FF^\top\,,
\\[-.1em]\nonumber
&\hspace*{4em}\hh=\hh_{\rm ext}{+}\nabla\Delta^{-1}{\rm div}(\mm)\,,\ \ \ \
\KK_\LAM(\FF,\nabla\mm)=
\frac{\upkappa_\LAM(\FF)}{\det\FF}\nabla\mm{\otimes}\nabla\mm\,,
  \\[.1em]\nonumber&\hspace*{4em}\SS_{\LAM,\EPS}(\FF,\mm,\theta)=
 {\rm skw}\big(\big(\mu_0\hh
 {-}\wh\tt_{\LAM,\EPS}(\FF,\mm,\theta)\big){\otimes}\mm\big)\,,\ \ \ 
  \mathscr{H}=\nu_2|\nabla^2\vv|^{p-2}\nabla^2\vv\,,
\\[.1em]
&\hspace*{4em}\mathscr{S}_\LAM(\FF,\mm,\nabla\mm)=
\frac{\upkappa_\LAM(\FF)}{\det\FF}{\rm Skw}(\nabla\mm{\otimes}\mm)\,,
    \label{Euler-thermo-reg1}
    \\[-.2em]
&\pdt\Fe=(\Nabla\vv)\Fe-(\vv{\cdot}\nabla)\Fe\,,
\label{Euler-thermo-reg2}
\\[-.2em]&\nonumber
\tau\rr+\HC(\FF,\theta){\rm Dir}(\rr)-\frac{\mm{\times}\rr}{\gamma(\FF,\mm,\theta)}
\ni\mu_0\hh-\wh\tt_{\LAM,\EPS}(\FF,\mm,\theta)
-{\rm div}\Big(\frac{\upkappa_\LAM(\FF)\nabla\mm}{\det\Fe}\Big)
\ \
\\[-.2em]&\hspace{7.3em}\label{Euler-thermo-reg4}
\text{ with }\ \ 
\wh\tt_{\LAM,\EPS}(\FF,\mm,\theta)=
\frac{\pi_\LAM(\FF)\upvarphi_\mm'(\Fe,\mm)}{\det\Fe}
+\frac{\pi_\LAM(\FF)\COUPLING_\mm'(\Fe,\mm,\theta)}{(1{+}\EPS|\theta|^{1/2})\det\Fe}\,,
\\[-.5em]&\pdt\mm={\rm skw}(\nabla\vv)\mm-(\vv{\cdot}\nabla)\mm+\rr\,,
\label{Euler-thermo-reg5}\\&\nonumber
\pdt{\W}=\xi_\EPS(\FF,\theta;\ee(\vv),\nabla^2\vv,\rr)
+{\rm div}\big(\COND(\FF,\theta)\nabla\theta-\W\vv\big)
+\frac{\pi_\LAM(\Fe)\COUPLING'_{\Fe}(\Fe,\mm,\theta)\Fe^\top\!\!}{(1{+}\EPS|\theta|
)\det\Fe}{:}\ee(\vv)
\\[-.3em]&\nonumber\hspace{5em}
+\frac{\pi_\LAM(\Fe)
\COUPLING_\mm'(\FF,\mm,\theta)}{(1{+}\EPS|\theta|^{1/2})\det\FF}{\cdot}
\big(\rr{-}{\rm skw}(\nabla\vv)\mm\big)\ \ \ \ \ \text{ with }\ \ \
\W=\OMEGA(\FF,\mm,\theta)
\\&
\hspace{9em}
\text{ and }\
\xi_\EPS(\FF,\theta;\ee,\bm{G},\rr):=
\frac{\DIS|\ee|^p{+}\tau|\rr|^2{+}\HC(\FF,\theta)|\rr|
{+}\nu_2|\bm{G}|^p}{1{+}\EPS|\EE|^p{+}\EPS|\bm{G}|^p{+}\EPS|\rr|^2}\,,
\label{Euler-thermo-reg3}
\end{align}\end{subequations}
where $\OMEGA(\cdot,\cdot)$ is from \eq{Euler-thermodynam3}.
We complete this system with the correspondingly regularized boundary
conditions on $I{\times}\varGamma$:
\begin{subequations}\label{Euler-thermo-reg-BC-IC}
\begin{align}\nonumber
&
\Big[\big(\TT_{\LAM,\EPS}(\FF,\mm,\theta)
{+}\KK_\LAM(\FF,\nabla\mm){+}\SS_{\LAM,\EPS}(\FF,\mm,\theta)
{+}\DIS|\ee(\vv)|^{p-2}\ee(\vv)
\\[-.2em]&\hspace{.5em}
{-}{\rm div}(\mathscr{H}{+}\mathscr{S_\LAM(\FF,\mm,\nabla\mm)})\big)\nn{-}\divS\big(\mathscr{H}
\nn{+}\mathscr{S}_\LAM(\FF,\mm,\nabla\mm)\nn\big)\Big]_\text{\sc t}^{}\!\!
+\nu_\flat|\vv|^{p-2}\vv=\TRACTION\,,
\label{Euler-thermo-reg-BC-IC-1}
\\&
\vv{\cdot}\nn=0,\ \ \ \ \nabla^2\vv{:}(\nn{\otimes}\nn)={\bm0}\,,\ \ \
\frac{\upkappa_\LAM(\FF)(\nn{\cdot}\nabla)\mm\!\!}{\det\Fe}=\bm0
\,,\ \ \ \ u(\infty)=0\,,\ \ \ \text{ and }\ \ \ 
\\&
\COND(\FF,\theta)\nabla\theta{\cdot}\nn-\frac{\nu_\flat|\vv|^p}{2{+}\EPS|\vv|^p}=h_\EPS(\theta):=\frac{h(\theta)}{1{+}\varepsilon|h(\theta)|}
\label{Euler-thermo-reg-BC-IC-2}
\intertext{and initial conditions}
&\varrho(0)=\varrho_0\,,
\ \ \ \ \vv(0)=\vv_0\,,\ \ \ \ 
\FF(0)=\FF_0\,,\ \ \ \ \mm(0)=\mm_0\,,\ \ \ \ \theta(0)=\theta_{0,\varepsilon}:=
\frac{\theta_{0}}{1{+}\varepsilon\theta_{0}}\,.
\end{align}\end{subequations}

Note that $\pi_\LAM\upvarphi\in C^1(\R^{d\times d})$ if
$\upvarphi\in C^1({\rm GL}^+(d))$ and that 
$[\pi_\LAM\upvarphi]'$ together with the regularized Cauchy stress $\TT_{\LAM,\EPS}$
are bounded, continuous, and vanish if $\FF$ ``substantially'' violates
the constraints \eq{Euler-quasistatic-est-formal4}.
Altogether, also from \eq{cut-off-kappa} and \eq{cut-off-K} we have:
$$
\bigg(\!\!\det\FF\le\frac\LAM2\ \ \text{ or }\ \ |\FF|\ge\frac2\LAM\bigg)
\ \ \ \Rightarrow\ \ \ \TT_{\LAM,\EPS}={\bm0},\ \ \ \frac{\upkappa_\LAM(\FF)}{\det\FF}=1,
\ \text{ and }\ \COND(\FF,\theta)=1,
$$
where a continuous (and smooth) extension of $[\upkappa_\LAM/\!\det](\cdot)$
at $0$ is considered. Also, due to \eq{Euler-ass-adiab},
$\TT_{\LAM,\EPS}:\R^{d\times d}\times\R^d\times\R\to\R_{\rm sym}^{d\times d}$ is
bounded. Thus the part (\ref{Euler-thermo-reg}a,b,d,e) of the
system allows for bounds for $\varrho_\EPS\vv_\EPS$, $\nabla\vv_\EPS$,
and $\mm_\EPS$ independently of 
$(\varrho_\EPS,\FFeps,\theta_\EPS)$ which can be estimated subsequently.
Moreover, recalling the extension \eq{extension-negative}, note that
it is defined also for negative temperatures.

The corresponding weak formulation of
\eq{Euler-thermo-reg}--\eq{Euler-thermo-reg-BC-IC} a'la
Definition~\ref{def} is quite straightforward and we will not 
explicitly write it, also because it will be obvious from its Galerkin
version \eq{Euler-weak-Galerkin} below. The philosophy of the
regularization \eq{Euler-thermo-reg} is that the estimation of the
magneto-mechanical part (\ref{Euler-thermo-reg}a-e) and of the thermal part 
\eq{Euler-thermo-reg3} decouples since $\TT_{\LAM,\EPS}$ and $\SS_{\LAM,\EPS}$ are
bounded and that heat sources in the heat equation are bounded for fixed
$\EPS>0$. Simultaneously, the heat equation has a non-negative solution and,
for $\EPS\to0$, the physical a~priori estimates
are the same as the formal estimates \eq{Euler-est}--\eq{est+} and,
when taking $\LAM>0$ small to comply with \eq{Euler-quasistatic-est-formal4},
the $\LAM$-regularization becomes eventually inactive, cf.\ Step~9 below.

\medskip\noindent{\it Step 2: a semi-discretization}.
For $\varepsilon>0$ fixed, we use a spatial semi-discretization, keeping the
transport equations \eq{Euler-thermo-reg0}, \eq{Euler-thermo-reg2}, and
\eq{Euler-thermo-reg5} continuous
(i.e.\ non-discretised) when exploiting Lemmas~\ref{lem1} and \ref{lem2}.
More specifically, we make a conformal Galerkin approximation of
\eq{Euler-thermo-reg1} by using a collection of nested finite-dimensional
subspaces $\{V_k\}_{k\in\N}$ whose union is dense in $W^{2,p}(\varOmega;\R^d)$
and a conformal Galerkin approximation of \eq{Euler-thermo-reg4} and of
\eq{Euler-thermo-reg3} by using  a collection of nested finite-dimensional
subspaces $\{Z_k\}_{k\in\N}$ whose union is dense in $H^1(\varOmega)$.
Without loss of generality, we can assume $\vv_0\in V_1$ and
$\theta_{0,\varepsilon}\in Z_1$.

The approximate solution of the regularized system \eq{Euler-thermo-reg}
will be denoted by
$$
(\varrho_{\EPS k},\vv_{\EPS k},\FFepsk,\rr_{\EPS k},\mm_{\EPS k},
\theta_{\EPS k}):I\to W^{1,r}(\varOmega){\times}V_k{\times}
W^{1,r}(\varOmega;\R^{d\times d}){\times}Z_k^d{\times}H^1(\varOmega;\R^d){\times}Z_k\,.
$$
Specifically, such a six-tuple should satisfy
\begin{subequations}\label{Euler-weak-Galerkin}\begin{align}
&\pdt{\varrho_{\EPS k}}=-{\rm div}(\varrho_{\EPS k}\vv_{\EPS k}) 
\ \ \ \text{ in the $L^1(I{\times}\varOmega)$-sense,}
\label{Euler-weak-Galerkin=1}
\\[-.0em]&\pdt{\FFepsk}=(\nabla\vv_{\EPS k})\FFepsk
-(\vv_{\EPS k}{\cdot}\nabla)\FFepsk\ \ \ \text{ in the $L^1(I{\times}\varOmega;\R^{d\times d})$-sense, and}
\label{Euler-weak-Galerkin-F}
\\[-.0em]&\label{Euler-weak-Galerkin-m}
\pdt{\mm_{\EPS k}}={\rm skw}(\nabla\vv_{\EPS k})\mm_{\EPS k}
-(\vv_{\EPS k}{\cdot}\nabla)\mm_{\EPS k}+\rr_{\EPS k}
\ \ \ \text{ in the weak sense}, 
\intertext{relying on $\varrho_{\EPS k}\in W^{1,1}(I{\times}\varOmega)$
and $\FFepsk\in W^{1,1}(I{\times}\varOmega;\R^{d\times d})$ 
which will be indeed proved later, together with the following integral identities}
&\nonumber
\int_0^T\!\!\!\int_\varOmega\!
\bigg(\Big(\TT_{\LAM,\EPS}(\FFepsk,\mm_{\EPS k},\theta_{\EPS k})
{+}\KK_\LAM(\FFepsk,\nabla\mm_{\EPS k}){-}\varrho_{\EPS k}\vv_{\EPS k}{\otimes}\vv_{\EPS k}
{+}\DIS|\ee(\vv_{\EPS k})|^{p-2}\ee(\vv_{\EPS k})
\Big){:}\ee(\widetilde\vv)
\\[-.1em]&\hspace{2em}\nonumber
+\SS_{\LAM,\EPS}(\FFepsk,\mm_{\EPS k},\theta_{\EPS k}){:}{\rm skw}(\widetilde\vv)
-\mu_0\big(\nabla\hh_{\EPS k}\big){:}(\mm_{\EPS k}{\otimes}\wt\vv)
-\mu_0\hh_{\EPS k}{\cdot}\mm_{\EPS k}({\rm div}\,\widetilde\vv)
-\varrho_{\EPS k}\vv_{\EPS k}{\cdot}\pdt{\widetilde\vv}
\\[-.2em]&\hspace{2em}\nonumber
+\big(\nu_2|\nabla^2\vv_{\EPS k}|^{p-2}\nabla^2\vv_{\EPS k}{+}
\mathscr{S}_\LAM(\FFepsk,\mm_{\EPS k},\nabla\mm_{\EPS k})\big)
\Vdots\Nabla^2\widetilde\vv\bigg)\,\d\xx\d t
=\!\int_\varOmega\varrho_0\vv_0{\cdot}\widetilde\vv(0)\,\d\xx
\\[-.2em]&\hspace{2em}
+\int_0^T\!\!\!\int_\varOmega\sqrt{\frac{\rhoR\varrho_{\EPS k}}
{\det_\LAM(\FFepsk)}}\GRAVITY{\cdot}\widetilde\vv
\,\d\xx\d t
+\!\int_0^T\!\!\!\int_\varGamma(\TRACTION-\nu_\flat|\vv_{\EPS k}|^{p-2}\vv_{\EPS k}){\cdot}\widetilde\vv\,\d S\d t
\label{Euler1-weak-Galerkin}
\intertext{with $\TT_{\LAM,\EPS}$, $\KK_\LAM$, $\SS_{\LAM,\EPS}$, and $\mathscr{S}_\LAM$
from \eq{Euler-thermo-reg1} and
$\hh_{\EPS k}=\hh_{\rm ext}{+}\nabla^2\Delta^{-1}{\rm div}(\mm_{\EPS k})$
for any $\widetilde\vv\in L^\infty(I;V_k)$ with
$\widetilde\vv{\cdot}\nn=0$ on $I{\times}\varGamma$ and $\widetilde\vv(T)={\bm0}$, and }
\nonumber
&\int_0^T\!\!\!\int_\varOmega\!\bigg(\frac\tau2|\wt\rr|^2+\HC(\FFepsk,\theta_{\EPS k})|\wt\rr|
-{\rm div}\frac{\upkappa_\LAM(\FFepsk)\nabla\mm_{\EPS k}}{\MU\det\FFepsk}{\cdot}
{\rm skw}(\nabla\vv_{\EPS k})\mm_{\EPS k}
+\frac{\upkappa_\LAM(\FFepsk)\nabla\mm_{\EPS k}}{\det\FFepsk}{:}\nabla\wt\rr
\\[-.3em]&\hspace*{3em}\nonumber
-\big(\mu_0\hh-\wh\tt_{\LAM,\EPS}(\FFepsk,\mm_{\EPS k},\theta_{\EPS k})\big)
{\cdot}(\wt\rr-\rr_{\EPS k})
+\frac{\mm_{\EPS k}{\times}\rr_{\EPS k}}{\gamma(\FFepsk,\mm_{\EPS k},\theta_{\EPS k})}{\cdot}\wt\rr+\Big(\KK_\LAM(\FFepsk,\nabla\mm_{\EPS k})
\\[-.1em]&\nonumber\hspace*{3em}
-\MU\frac{\upkappa_\LAM'(\FFepsk)|\Nabla\mm_{\EPS k}|^2\FFepsk^{\top}\!}{2\det\FFepsk}
\Big){:}\ee(\vv_{\EPS k})
\bigg)\,\d\xx\d t
\ge\int_0^T\!\!\!\int_\varOmega\!\frac\tau2|\rr_{\EPS k}|^2+\HC(\FF,\theta)|\rr_{\EPS k}|\,\d\xx\d t
\\[.1em]&\hspace*{14em}+\!\int_\varOmega\!
\frac{\upkappa_\LAM(\FFepsk(T))|\Nabla\mm_{\EPS k}(T)|^2}{2\MU\det\FFepsk(T)}-\frac{\upkappa_\LAM(\FF_0)|\Nabla\mm_0|^2\!\!}{2\MU\det\FF_0}\,\d\xx
\label{Euler2-weak-Galerkin}
\intertext{holding for any $\wt\rr\in Z_k^d$ and for a.a.\ $t\in I$, and further}
\nonumber
&\!\!\!\int_0^T\!\!\!\int_\varOmega\bigg(\W_{\EPS k}\pdt{\widetilde\theta}
+\big(\W_{\EPS k}\vv_{\EPS k}
{-}\COND(\FFepsk,\theta_{\EPS k})\nabla\theta_{\EPS k}\big)
{\cdot}\nabla\widetilde\theta
+\xi_\EPS(\FFepsk,\theta_{\EPS k};\ee(\vv_{\EPS k}),\nabla^2\vv_{\EPS k},\rr_{\EPS k})
\widetilde\theta
\\[-.2em]&\hspace{-.3em}\nonumber
+\pi_\LAM(\FFepsk)\Big(
\frac{\COUPLING'_{\Fe}(\FFepsk,\mm_{\EPS k},\theta_{\EPS k})\FFepsk^\top}
{(1{+}\EPS|\theta_{\EPS k}|)\det\FFepsk}
{:}\ee(\vv_{\EPS k})
+\frac{\!\!\COUPLING_\mm'(\FFepsk,\mm_{\EPS k},\theta_{\EPS k}){\cdot}
\big(\rr_{\EPS k}{-}{\rm skw}(\nabla\vv_{\EPS k})\mm_{\EPS k}\big)}
{(1{+}\EPS|\theta_{\EPS k}|^{1/2})\det\FFepsk}\Big)
\widetilde\theta\,\bigg)\d\xx\d t
\\[-.2em]&\hspace{-.3em}\nonumber
+\!\int_\varOmega\!\OMEGA(\FF_0,\theta_{0,\varepsilon})\widetilde\theta(0)\,\d\xx
+\!\int_0^T\!\!\!\int_\varGamma\!\Big(h_\EPS(\theta_{\EPS k}){+}
\frac{\nu_\flat|\vv_{\EPS k}|^p}{2{+}\EPS|\vv_{\EPS k}|^p}\Big)
\widetilde\theta\,\d S\d t=0
\\[.1em]&\hspace{20em}
\ \ \ \ \text{ with }\ \ \ \W_{\EPS k}=
\OMEGA(\FFepsk,\mm_{\EPS k},\theta_{\EPS k})
\label{Euler3-weak-Galerkin}
\end{align}
\end{subequations}
holds for any $\widetilde\theta\in C^1(I;Z_k)$ with $\widetilde\theta(T)=0$.

Existence of this solution is based on the standard theory of systems of
ordinary differential equations first locally in time combined here
with the abstract $W^{1,r}(\varOmega)$- and $L^2(\varOmega)$
valued differential equations based on Lemmas~\ref{lem1} and \ref{lem2}
for the scalar, the vector, and the tensor transport equations
(\ref{Euler-weak-Galerkin}a--c)
and then by successive prolongation on the whole time interval based on the
$L^\infty$-estimates below. Usage of Lemmas~\ref{lem1} and \ref{lem2}
with the fixed initial conditions $\varrho_0$, $\FF_0$, and $\mm_0$ 
defines the nonlinear operators $\mathfrak{R}:I\times L^p(I;W^{2,p}(\varOmega;\R^d))
\to W^{1,r}(\varOmega)$, $\mathfrak{F}:I\times L^p(I;W^{2,p}(\varOmega;\R^d))
\to W^{1,r}(\varOmega;\R^{\d\times d})$, and
$\mathfrak{M}:I\times L^p(I;W^{2,p}(\varOmega;\R^d))\times L^2(I{\times}\varOmega;\R^d)
\to L^2(\varOmega;\R^d)$ by 
\begin{align}
\varrho_{\EPS k}(t)=\mathfrak{R}\big(t,\vv_{\EPS k}\big)\,,
\ \ \ \ \ 
\FFepsk(t)=\mathfrak{F}\big(t,\vv_{\EPS k}\big)\,,\ \text{ and }\ 
\mm_{\EPS k}(t)=\mathfrak{M}\big(t,\vv_{\EPS k},\rr_{\EPS k}\big)\,.
\end{align}

\medskip\noindent{\it Step 3: first a~priori estimates}.
In the Galerkin approximation, it is legitimate to use
$\widetilde\vv=\vv_{\EPS k}$ for \eq{Euler1-weak-Galerkin} and
$\widetilde\theta=\theta_{\EPS k}$ for \eq{Euler3-weak-Galerkin}.
We take the benefit from having the transport equations
(\ref{Euler-weak-Galerkin}a,b) non-discretized and thus we can test them
by the nonlinearities $|\vv_{\EPS k}|^2/2$ and
$[\pi_\LAM\upvarphi(\FFepsk,\mm_{\EPS k})/\!\det\FFepsk]_\FF'$, 
respectively. In particular, we can use the calculus
\eq{Euler-large-thermo} which holds also for 
$\pi_\LAM\upvarphi$ instead of $\upvarphi$ and the calculus
\eq{calculus-convective-in-F} also for the semi-Galerkin approximate solution.
Also we can use the calculus \eq{test-damage}--\eq{formula3} with
$\pi_\LAM\COUPLING_\FF'/(1{+}|\theta|)$ and $\pi_\LAM\COUPLING_\mm'/(1{+}|\theta|^{1/2})$ 
and $\upkappa_\LAM$ instead of $\COUPLING_\FF'$ and $\COUPLING_\mm'$ and $\upkappa$,
respectively, using also that we have the nondiscretized equation
\eq{Euler-weak-Galerkin-m} at disposal. The philosophy of the regularization
\eq{Euler-thermo-reg} is that, for this estimation procedure, the system
decouples to the magneto-mechanical part and the thermal part which allows
for basic estimates independent of  $\vv_{\EPS k}$, $\mm_{\EPS k}$, and $\rr_{\EPS k}$.

Specifically,  from \eq{Euler1-weak-Galerkin} tested by $\vv_{\EPS k}$
and \eq{Euler2-weak-Galerkin} tested by $\bm0$, like \eq{energy+} we obtain
the inequality 
\begin{align}\nonumber
  &\hspace*{0em}\frac{\d}{\d t}
  \int_\varOmega\!\frac{\varrho_{\EPS k}}2|\vv_{\EPS k}|^2+
  \frac{\pi_\LAM(\FFepsk)\upvarphi(\FFepsk,\mm_{\EPS k})\!}{\det\FFepsk}
  +\frac{\upkappa_\LAM(\FFepsk)}{2\det\FFepsk}|\Nabla\mm_{\EPS k}|^2
  -\mu_0\hh_\text{\rm ext}{\cdot}\mm_{\EPS k}\,\d\xx
\\[-.1em]&\hspace{.0em}\nonumber
\qquad+\!\int_\varOmega\!
\xi_\EPS(\FFepsk,\theta_{\EPS k};\ee(\vv_{\EPS k}),\nabla^2\vv_{\EPS k},\rr_{\EPS k})
\,\d\xx
+\!\int_\varGamma\!\nu_\flat|\vv_{\EPS k}|^p\,\d S
+\frac{\d}{\d t}\int_{\R^d}\frac{\mu_0}2|\Nabla u_{\EPS k}|^2\,\d\xx
\\[-.1em]&\nonumber\hspace{.0em}
\le\int_\varGamma\!\TRACTION{\cdot}\vv_{\EPS k}\,\d S
+\int_\varOmega\bigg(\sqrt{\frac{\rhoR\varrho_{\EPS k}}{\det_\LAM(\FFepsk)}}\GRAVITY{\cdot}\vv_{\EPS k}
-\frac{\pi_\LAM(\FFepsk)\COUPLING_{\Fe}'(\FFepsk,\mm_{\EPS k},\theta_{\EPS k})
\FFepsk^\top}{(1{+}\EPS|\theta_{\EPS k}|)\det\FFepsk}{:}\ee(\vv_{\EPS k})
\\[-.1em]&\hspace{.0em}
\qquad\qquad-\frac{\pi_\LAM(\FFepsk)
\COUPLING_\mm'(\FFepsk,\mm_{\EPS k},\theta_{\EPS k})\!}
{(1{+}\EPS|\theta_{\EPS k}|^{1/2})\det\FFepsk}
{\cdot}\big(\rr_{\EPS k}{-}{\rm skw}(\nabla\vv_{\EPS k})\mm_{\EPS k}\big)
-\mu_0\frac{\!\partial \hh_\text{\rm ext}\!}{\partial t}{\cdot}\mm_{\EPS k}\!\bigg)\,\d\xx
\label{thermodynamic-Euler-mech-disc}
\end{align}
with $u_{\EPS k}=\Delta^{-1}{\rm div}(\mm_{\EPS k})$.

The bulk-force term can be estimated similarly as \eq{Euler-est-of-rhs}: 
\begin{align}\nonumber
\int_\varOmega\sqrt{\frac{\rhoR}{\det_\LAM(\FFepsk)}}
&\sqrt{\varrho_{\EPS k}}\vv_{\EPS k}{\cdot}\GRAVITY\,\d \bm{x}
\le
\Big\|\sqrt{\frac{\rhoR}{\det_\LAM(\FFepsk)}}\Big\|_{L^{2}(\varOmega)}\big\|\sqrt{\varrho_{\EPS k}}\vv_{\EPS k}\big\|_{L^2(\varOmega;\R^d)}^{}\big\|\GRAVITY\big\|_{L^\infty(\varOmega;\R^d)}^{}
\\&\nonumber\le
\big\|\GRAVITY\big\|_{L^\infty(\varOmega;\R^d)}^{}\int_\varOmega
\frac{\rhoR}{2\det_\LAM(\FFepsk)\!}+\frac{\varrho_{\EPS k}}2|\vv_{\EPS k}|^2\,\d\xx
\\&
=\big\|\GRAVITY\big\|_{L^\infty(\varOmega;\R^d)}^{}
\Big(\frac{{\rm meas}(\varOmega)(1{+}\max_{\varOmega}\rhoR)\!}{2\LAM}\,
+\int_\varOmega\frac{\varrho_{\EPS k}}2|\vv_{\EPS k}|^2\,\d\xx\Big)\,
\label{Euler-est-bulk-force}\end{align}
while the other terms  in \eq{thermodynamic-Euler-mech-disc} can be estimated as done
in \eq{h.v-est}, \eq{Zeeman-est},  \eq{Euler-est-of-rhs+}, and \eq{Euler-est-of-rhs++}.

By the Gronwall inequality, we obtain the estimates
\begin{subequations}\label{Euler-quasistatic-est1}
\begin{align}\label{Euler-quasistatic-est1-1}
&\|\vv_{\EPS k}\|_{L^p(I;W^{2,p}(\varOmega;\R^d))}^{}\le C_\EPS\ \ \ \text{ and }\ \ \
\big\|\sqrt{\varrho_{\EPS k}}\vv_{\EPS k}\big\|_{L^\infty(I;L^2(\varOmega;\R^d))}\le C_\EPS\,,
\\&
\|\rr_{\EPS k}\|_{L^2(I\times\varOmega;\R^d)}^{}\le C_\EPS\,,\ \ \ \text{ and }\ \ \
\|u_{\EPS k}\|_{L^\infty(I;H^1(\R^d))}^{}\le C_\EPS
\intertext{with $C_\EPS$ depending on $\EPS>0$ considered fixed in this Step.
The former estimate in \eq{Euler-quasistatic-est1-1} relies also on the Navier-type
boundary conditions and allows us to use Lemma~\ref{lem1} to obtain the estimate}
&\label{Euler-est1-2}
\|\FFepsk\|_{L^\infty(I;W^{1,r}(\varOmega;\R^{d\times d}))}^{}\le C_{r,\EPS}\ \ \text{ with }\ \
\Big\|\frac1{\det\FFepsk}\Big\|_{L^\infty(I;W^{1,r}(\varOmega))}^{}\!\le C_{r,\EPS}\,,
\\&\label{Euler-est1-3}
\|\varrho_{\EPS k}\|_{L^\infty(I;W^{1,r}(\varOmega))}^{}\le C_{r,\EPS}
\ \ \ \text{ with }\ \ \ \Big\|\frac1{\varrho_{\EPS k}}\Big\|_{L^\infty(I;W^{1,r}(\varOmega))}^{}\!\le C_{r,\EPS}\,,
\\&\label{Euler-est1-4}
\|\vv_{\EPS k}\|_{L^\infty(I;L^2(\varOmega;\R^d))}^{}\le C_\EPS\,,\ \ \ \text{ and }
\\\label{Euler-est1-5}
&\|\mm_{\EPS k}\|_{L^\infty(I;H^1(\varOmega;\R^d))}^{}\le C_\EPS\,;
\end{align}\end{subequations}
for \eq{Euler-est1-4} and \eq{Euler-est1-5} we used argumentation like in
\eq{basic-est-of-v} and \eq{Euler-est1-5}, respectively.
For $\LAM$ and $\EPS$ fixed, it is important
that these estimates can be made independently of $\theta$ since $\TT_{\LAM,\EPS}$ is
a-priori bounded. It is also important that, due to the latter estimate in
\eq{Euler-est1-2}, the singularity in $\COUPLING(\cdot,\mm,\theta)$ is
not active and $\OMEGA(\FFepsk,\mm_{\EPS k},\theta_{\EPS k})$ is well defined.

Also, we have the time derivatives estimated by comparison from
\eq{Euler-weak-Galerkin} as
\begin{align}
\Big\|\pdt{\varrho_{\EPS k}}\Big\|_{L^p(I;L^r(\varOmega))}\!\le C_\EPS,\ \ \ \
\Big\|\pdt{\FFepsk}\Big\|_{L^p(I;L^r(\varOmega;\R^{d\times d}))}\!\le C_\EPS,\ \ \ \
\Big\|\pdt{\mm_{\EPS k}}\Big\|_{L^2(I{\times}\varOmega;\R^d)}\!\le C_\EPS.
\label{est-time-derivaitves}\end{align}

The further estimates can be obtained by testing the Galerkin approximation
of \eq{Euler3-weak-Galerkin} by $\widetilde\theta=\theta_{\EPS k}$, This
is to be made carefully not to see terms as
$\theta_{\EPS k}\OMEGA_\Fe'(\FFepsk,\mm_{\EPS k},\theta_{\EPS k}){:}(\vv_{\EPS k}{\cdot}\nabla)\FFepsk$
which is not integrable. To this goal, we consider the convective-derivative
form $\pdt{}\W+{\rm div}(\W\vv)=\DT\W+\W\,{\rm div}\vv$ in
\eq{Euler-thermo-reg3}. We denote by $\widehat\OMEGA(\FF,\mm,\theta)$
a primitive function to
$\theta\mapsto\theta\OMEGA_\theta'(\Fe,\mm,\theta)$ depending smoothly
on $\Fe$ and on $\mm$, specifically
\begin{align}
\widehat\OMEGA_\LAM(\Fe,\mm,\theta)=\int_0^1\!\!r\theta^2
\OMEGA_\theta'(\Fe,\mm,r\theta)\,\d r\,.
\label{primitive}\end{align}
For $\W_{\EPS k}=\OMEGA(\FFepsk,\mm_{\EPS k},\theta_{\EPS k})$ and
using \eq{Euler-weak-Galerkin-F}, the mentioned test by $\theta_{\EPS k}$ then gives
\begin{align}\nonumber
\theta_{\EPS k}\DT\W_{\EPS k}&=
\theta_{\EPS k}\OMEGA_\Fe'(\FFepsk,\mm_{\EPS k},\theta_{\EPS k}){:}\DT\FFepsk
\\&\nonumber
\ \ +\theta_{\EPS k}\OMEGA_\mm'(\FFepsk,\mm_{\EPS k},\theta_{\EPS k}){\cdot}
\DT\mm_{\EPS k}+
\theta_{\EPS k}\OMEGA_\theta'(\FFepsk,\mm_{\EPS k},\theta_{\EPS k})\DT\theta_{\EPS k}
\\\nonumber&
=\Big(\theta_{\EPS k}\OMEGA_\Fe'(\FFepsk,\mm_{\EPS k},\theta_{\EPS k})
-\widehat\OMEGA_\Fe'(\FFepsk,\mm_{\EPS k},\theta_{\EPS k})\Big){:}
(\nabla\vv_{\EPS k})\FFepsk
\\&\ \ +\Big(\theta_{\EPS k}\OMEGA_\mm'(\FFepsk,\mm_{\EPS k},\theta_{\EPS k})
{-}\widehat\OMEGA_\mm'(\FFepsk,\mm_{\EPS k},\theta_{\EPS k})\Big){\cdot}
\DT\mm_{\EPS k}\!
+\DT{\overline{\widehat\OMEGA(\FFepsk,\mm_{\EPS k},\theta_{\EPS k})}}\,.
\label{Euler-thermodynam3-test++}\end{align}
Integrating the last term over $\varOmega$ gives, by the Green formula,
$\int_\varOmega\DT{\overline{\widehat\OMEGA(\FFepsk,\mm_{\EPS k},\theta_{\EPS k})}}\,\d\xx
=\frac{\d}{\d t}\int_\varOmega\widehat\OMEGA(\FFepsk,\mm_{\EPS k},\theta_{\EPS k})\,\d\xx
-\int_\varOmega\widehat\OMEGA(\FFepsk,\mm_{\EPS k},\theta_{\EPS k}){\rm div}\vv_{\EPS k}\,\d\xx+\int_\varGamma\widehat\OMEGA(\FFepsk,\mm_{\EPS k},\theta_{\EPS k})\vv_{\EPS k}{\cdot}\nn\,\d S$. Thus, we obtain:
\begin{align}\nonumber
&\frac{\d}{\d t}\int_\varOmega\widehat\OMEGA(\FFepsk,\mm_{\EPS k},\theta_{\EPS k})\,\d\xx
+\int_\varOmega\COND(\FFepsk,\theta_{\EPS k})|\nabla\theta_{\EPS k}|^2\,\d\xx
\\[-.2em]&\hspace{.5em}\nonumber
=\int_\varOmega\!\bigg(
\Big(\xi_\EPS(\FFepsk,\theta_{\EPS k};\ee(\vv_{\EPS k}),\nabla^2\vv_{\EPS k},\rr_{\EPS k})
+\frac{\pi_\LAM(\FFepsk)
\COUPLING'_{\Fe}(\FFepsk,\mm_{\EPS k},\theta_{\EPS k})
\FFepsk^\top}{(1{+}\EPS|\theta_{\EPS k}|)\det\FFepsk}{:}\ee(\vv_{\EPS k})
\\[.1em]&\nonumber\hspace{3.5em}
+\frac{\pi_\LAM(\FFepsk)
\COUPLING'_\mm(\FFepsk,\mm_{\EPS k},\theta_{\EPS k}){\cdot}
\big(\rr_{\EPS k}{-}{\rm skw}(\nabla\vv_{\EPS k})\mm_{\EPS k}\big)}
{(1{+}\EPS|\theta_{\EPS k}|^{1/2})\det\FFepsk}-\OMEGA(\FFepsk,\mm_{\EPS k},\theta_{\EPS k})
{\rm div}\,\vv_{\EPS k}\Big)\theta_{\EPS k}
\\[-.3em]&\nonumber\hspace{3.5em}
-\Big(\theta_{\EPS k}\OMEGA_\Fe'(\FFepsk,\mm_{\EPS k},\theta_{\EPS k})
-\widehat\OMEGA_\Fe'(\FFepsk,\mm_{\EPS k},\theta_{\EPS k})\Big){:}
(\nabla\vv_{\EPS k})\FFepsk
\\[-.0em]&\nonumber\hspace{3.5em}
-\Big(\theta_{\EPS k}\OMEGA_\mm'(\FFepsk,\mm_{\EPS k},\theta_{\EPS k})
{-}\widehat\OMEGA_\mm'(\FFepsk,\mm_{\EPS k},\theta_{\EPS k})\Big){\cdot}
\big(\rr_{\EPS k}{-}{\rm skw}(\nabla\vv_{\EPS k})\mm_{\EPS k}\big)
\\[-.3em]&\hspace{3.5em}
+\widehat\OMEGA(\FFepsk,\mm_{\EPS k},\theta_{\EPS k}){\rm div}\,\vv_{\EPS k}
\bigg)\,\d\xx
+\int_\varGamma\Big(h_\EPS(\theta_{\EPS k})-\frac{\nu_\flat|\vv_{\EPS k}|^p}
{2{+}\EPS|\vv_{\EPS k}|^p}\Big)\theta_{\EPS k}\,\d S\,.
\label{Euler3-Galerkin-est}\end{align}

We integrate \eq{Euler3-Galerkin-est} in time over an interval
$[0,t]$ with $t\in I$. For the left-hand side, let us realize
that $\widehat\OMEGA(\Fe,\mm,\theta)\ge c_K^{}\theta^2$
due to \eq{Euler-ass-adiab} with $c_K^{}>0$ depending
on the fixed $\LAM>0$ used in \eq{Euler-quasistatic-est-formal4}.
Then the integrated right-hand side of \eq{Euler3-Galerkin-est} is to
be estimated from above, in particular relying on \eq{Euler-ass-adiab}
and on \eq{Euler-ass-primitive-c}. Let us discuss the difficult terms.
In view of \eq{primitive}, it holds that
$\wh\OMEGA_\Fe'(\Fe,\mm,\theta)=\int_0^1r\theta^2
\OMEGA_{\Fe\theta}''(\Fe,\mm,r\theta)\,\d r$ and
$\OMEGA_\mm'(\Fe,\mm,\theta)=$ $\int_0^1r\theta^2
\OMEGA_{\mm\theta}''(\Fe,\mm,r\theta)\,\d r$.
Recalling \eq{Euler-ass-primitive-c}, we have 
$|\wh\OMEGA_\Fe'(\Fe,\mm,\theta)
-\theta\OMEGA_\Fe'(\Fe,\mm,\theta)|\le C(1{+}|\theta|^2)$
and $|\wh\OMEGA_\mm'(\Fe,\mm,\theta)
-\theta\OMEGA_\mm'(\Fe,\mm,\theta)|\le C(1{+}|\theta|)$. 
It allows for estimation 
\begin{align}\nonumber
&\bigg|\int_\varOmega\Big(\wh\OMEGA_\Fe'(\FFepsk,\mm_{\EPS k},\theta_{\EPS k})-\theta_{\EPS k}
\OMEGA_\Fe'(\FFepsk,\mm_{\EPS k},\theta_{\EPS k})\Big){:}
(\nabla\vv_{\EPS k})\FFepsk\,\d\xx\bigg|
\\&\hspace{10em}\le
C\big(|\varOmega|+\|\theta_{\EPS k}\|_{L^{2}(\varOmega)}^{2}\big)
\big\|(\nabla\vv_{\EPS k})\FFepsk\big\|_{L^\infty(\varOmega;\R^{d\times d})}
\label{adiab-est-F}\ \ \ \ \text{ and}
\\\nonumber
&\bigg|\int_\varOmega\!\Big(\wh\OMEGA_\mm'(\FFepsk,\mm_{\EPS k},\theta_{\EPS k})
-\theta_{\EPS k}\OMEGA_\mm'(\FFepsk,\mm_{\EPS k},\theta_{\EPS k})\Big){\cdot}
\big(\rr_{\EPS k}{-}{\rm skw}(\nabla\vv_{\EPS k})\mm_{\EPS k}\big)\,\d\xx\bigg|
\\&\hspace{10em}\le
2C^2\big(|\varOmega|+\|\theta_{\EPS k}\|_{L^{2}(\varOmega)}^{2}\big)
\!+\big\|(\rr_{\EPS k}{-}{\rm skw}(\nabla\vv_{\EPS k})\mm_{\EPS k}\big\|_{L^2(\varOmega;\R^d)}^2\,.
\label{adiab-est-m}\end{align}
Using \eq{Euler-ass-primitive-c} together with $\OMEGA(\FF,\mm,0)=0$
so that $|\OMEGA(\FF,\mm,\theta|\le C_K|\theta|$,
the convective terms
$\OMEGA(\FFepsk,\mm_{\EPS k},\theta_{\EPS k})({\rm div}\vv_{\EPS k})\theta_{\EPS k}$
and $\wh\OMEGA(\FFepsk,\mm_{\EPS k},\theta_{\EPS k}){\rm div}\,\vv_{\EPS k}$
in \eq{Euler3-Galerkin-est} can be estimated as
\begin{align}
\!\!\!\!\!\bigg|\int_\varOmega\!\!\Big(\wh\OMEGA(\FFepsk,\mm_{\EPS k},\theta_{\EPS k})
{-}\theta_{\EPS k}\OMEGA(\FFepsk,\mm_{\EPS k},\theta_{\EPS k})\Big){\rm div}\vv_{\EPS k}
\d\xx\bigg|\le
\frac{\!C_K}2\big\|\theta_{\EPS k}\big\|_{L^2(\varOmega)}^2
\big\|{\rm div}\,\vv_{\EPS k}\big\|_{L^\infty(\varOmega)}.\!\!
\label{Euler3-Galerkin-est++}\end{align}
The terms $\|\theta_{\EPS k}\|_{L^{2}(\varOmega)}^{2}$ in \eq{adiab-est-F}
and in \eq{adiab-est-m} and in \eq{Euler3-Galerkin-est++}are to be treated by the
Gronwall inequality. The boundary term
in \eq{Euler3-Galerkin-est} can be estimated by \eq{Euler-ass-h}, taking also
into the account the extension \eq{extension-negative}, as 
\begin{align}\nonumber
\int_\varGamma\Big(h_\EPS(\theta_{\EPS k})-\frac{\nu_\flat|\vv_{\EPS k}|^p}
{2{+}\EPS|\vv_{\EPS k}|^p}\Big)\theta_{\EPS k}\,\d S
&\le C_{\EPS,\nu_\flat,a}+a\|\theta_{\EPS k}\|_{L^2(\varGamma)}^2
\\&\le C_{\EPS,\nu_\flat,a}+a N^2\big(\|\theta_{\EPS k}\|_{L^2(\varOmega)}^2
+\|\nabla\theta_{\EPS k}\|_{L^2(\varOmega;\R^d)}^2\big)\,,
\label{boundary-heat-est}\end{align}
where $C_{\EPS,\nu_\flat,\delta}$ depends also on $C$ from \eq{Euler-ass-h}, 
$N$ is the norm of the trace operator $H^1(\varOmega)\to L^2(\varGamma)$.
For $a>0$ in \eq{boundary-heat-est} sufficiently small, the last term
can be absorbed in the left-hand side of \eq{Euler3-Galerkin-est}.
Exploiting again the bound (\ref{Euler-quasistatic-est1}a,b), we eventually
obtain the estimate
\begin{subequations}\label{Euler-quasistatic-est2}
\begin{align}\label{Euler-quasistatic-est2-theta}
&\|\theta_{\EPS k}\|_{L^\infty(I;L^2(\varOmega))\,\cap\,L^2(I;H^1(\varOmega))}^{}\le C
\ \ \text{ and also}
\\&\|\W_{\EPS k}\|_{L^\infty(I;L^2(\varOmega))\,\cap\,L^2(I;H^1
(\varOmega))}^{}\le C.
\label{Euler-weak-sln-w}
\end{align}\end{subequations}
For \eq{Euler-weak-sln-w}, we used the calculus 
\begin{align}\nonumber
\nabla\W_{\EPS k}
=[\OMEGA_\LAM]_\theta'(\FFepsk,\mm_{\EPS k},\theta_{\EPS k})\nabla\theta_{\EPS k}
&+[\OMEGA_\LAM]_\FF'(\FFepsk,\mm_{\EPS k},\theta_{\EPS k})\nabla\FFepsk
\\&+[\OMEGA_\LAM]_\mm'(\FFepsk,\mm_{\EPS k},\theta_{\EPS k})\nabla\mm_{\EPS k}\in
L^2(I{\times}\varOmega;\R^d)
\label{w=...}\end{align}
together with the already proved information that
$|[\OMEGA_\LAM]_\FF'(\FFepsk,\mm_{\EPS k},\theta_{\EPS k})|$ is bounded in
$L^2(I;L^{2^*}(\varOmega))$ 
and $|[\OMEGA_\LAM]_\mm'(\FFepsk,\mm_{\EPS k},\theta_{\EPS k})|$ is bounded in
$L^\infty(I{\times}\varOmega)$ due to 
\eq{Euler-ass-primitive-c}, while $|\nabla\FFepsk|$ is bounded in
$L^\infty(I;L^r(\varOmega))$ and $|\nabla\mm_{\EPS k}|$ is
bounded in $L^\infty(I;L^2(\varOmega))$.

\medskip\noindent{\it Step 4: Limit passage in the magneto-mechanical part for $k\to\infty$}.
Using the Banach selection principle, we can extract some subsequence of
$\{(\varrho_{\EPS k},\vv_{\EPS k},\FFepsk,\mm_{\EPS k},\rr_{\EPS k},\W_{\EPS k})\}_{k\in\N}$
and its limit
$(\varrho_\EPS,\vv_\EPS,\FFeps,\mm_\EPS,\rr_\EPS,\W_\EPS):I\to W^{1,r}(\varOmega)\times
L^2(\varOmega;\R^d)\times W^{1,r}(\varOmega;\R^{d\times d})\times H^1(\varOmega;\R^d)
\times L^2(\varOmega;\R^d)\times L^2(\varOmega)$ such that
\begin{subequations}\label{Euler-weak-sln}
\begin{align}\label{Euler-weak-sln-rho}
&\!\!\varrho_{\EPS k}\to\varrho_\EPS&&\text{weakly* in $\
L^\infty(I;W^{1,r}(\varOmega))$}\,,
\\\label{Euler-weak-sln-v}
&\!\!\vv_{\EPS k}\to\vv_\EPS&&\text{weakly* in $\
L^\infty(I;L^2(\varOmega;\R^d))\cap
L^p(I;W^{2,p}(\varOmega;\R^d))$,}\!\!&&
\\\label{Euler-weak-sln-F}
&\!\!\FFepsk\to\FFeps
\!\!\!&&\text{weakly* in $\ L^\infty(I;W^{1,r}(\varOmega;\R^{d\times d}))$,}\!\!
\\\label{Euler-weak-sln-m}
&\!\!\mm_{\EPS k}\to\mm_\EPS
\!\!\!\!\!&&\text{weakly* in $\ L^\infty(I;H^1(\varOmega;\R^d))$,}\!\!
\\\label{Euler-weak-sln-r}
&\!\!\rr_{\EPS k}\to\rr_\EPS
\!\!\!&&\text{weakly \,\,in $\ \,L^2(I{\times}\varOmega;\R^d)$,}\!\!
\\
&\!\!\W_{\EPS k}\to\W_\EPS&&\text{weakly* in $\ 
L^\infty(I;L^2(\varOmega))\,\cap\,L^2(I;H^1(\varOmega))$.}
\label{Euler-weak-sln-w++}
\end{align}\end{subequations}

Relying on the assumption $r>d$ and on estimates \eq{est-time-derivaitves}
on $\pdt{}\varrho_{\EPS k}$, $\pdt{}\FFepsk$, and $\pdt{}\mm_{\EPS k}$, by the Aubin-Lions
lemma we also have that
\begin{subequations}\begin{align}\label{rho-conv}
&\varrho_{\EPS k}\to\varrho_\EPS\ \text{ strongly in }\ C(I{\times}\barOmega)\,,
\\&\FFepsk\to\FFeps\text{ strongly in
$C(I{\times}\barOmega;\R^{d\times d})$,}\ \ \text{ and }\ \ 
\\&\mm_{\EPS k}\to\mm_\EPS\ \text{ strongly in }\ C(I{\times}\barOmega;\R^d)\,.
\end{align}\end{subequations}
This already allows for the limit passage in the evolution-and-transport equations
\eq{Euler-weak-Galerkin}, cf.\ \eq{v-mapsto-F} and \eq{v-mapsto-m}.

Further, by comparison in the equation \eq{Euler-thermo-reg3} with the boundary condition
\eq{Euler-thermo-reg-BC-IC-2} in its Galerkin approximation,
we obtain a bound on \EEE $\pdt{}\W_{\EPS k}$ in seminorms $|\cdot|_l$ on $L^2(I;H^1(\varOmega)^*)$
arising from this Galerkin approximation, defined as
$|f|_l^{}:=\sup\{\int_0^T\!\int_\varOmega f\widetilde\theta\,\d\xx\d t;\
\|\widetilde\theta\|_{L^2(I;H^1(\varOmega))}^{}\le1, \
\widetilde\theta(t)\in Z_l\ \text{for }t\in I\}$. More specifically, for any $k\ge l$,
we can estimate
\begin{align}\nonumber&
\Big|\pdt{\W_{\EPS k}}\Big|_l^{}=\!\!
\sup\limits_{\stackrel{{\scriptstyle{\widetilde\theta(t)\in Z_l\ \text{for }t\in I}}}{{\scriptstyle{\|\widetilde\theta\|_{L^2(I;H^1(\varOmega))}^{}\le1}}}}
\int_0^T\!\!\!\int_\varOmega\bigg(
{-}\COND(\FFepsk,\theta_{\EPS k})\nabla\theta_{\EPS k}
{\cdot}\nabla\widetilde\theta
+\Big(\xi_\EPS(\FFepsk,\theta_{\EPS k};\ee(\vv_{\EPS k}),\nabla^2\vv_{\EPS k},\rr_{\EPS k})
\\[-2em]&\hspace{17em}\nonumber
+\frac{\pi_\LAM(\FFepsk)\COUPLING'_{\Fe}(\FFepsk,\mm_{\EPS k},\theta_{\EPS k})\FFepsk^\top
}{(1{+}\EPS|\theta_{\EPS k}|)\det\FFepsk}{:}\ee(\vv_{\EPS k})
\\[-.2em]&\hspace{8em}\nonumber
+\frac{\pi_\LAM(\FFepsk)\COUPLING'_\mm(\FFepsk,\mm_{\EPS k},\theta_{\EPS k})}
{(1{+}\EPS|\theta_{\EPS k}|^{1/2})\det\FFepsk}{\cdot}
\big(\rr_{\EPS k}{-}{\rm skw}(\nabla\vv_{\EPS k})\mm_{\EPS k}\big)
\Big)\widetilde\theta\bigg)\d\xx\d t
\\[-.1em]&\hspace{17em}
+\!\int_0^T\!\!\!\int_\varGamma\!
\Big(h_\EPS(\theta_{\EPS k}){+}
\frac{\nu_\flat|\vv_{\EPS k}|^p}{2{+}\EPS|\vv_{\EPS k}|^p}\Big)
\widetilde\theta\,\d S\d t\le C
\label{Euler3-weak-Galerkin+}\end{align}
with some $C$ depending on the estimates (\ref{Euler-quasistatic-est1}a,b)
and \eq{Euler-quasistatic-est2-theta} but independent on $l\in N$. 
Thus, by \eq{Euler-weak-sln-w++} and by a generalized Aubin-Lions theorem
\cite[Ch.8]{Roub13NPDE}, we obtain 
\begin{subequations}\label{Euler-weak+}
\begin{align}
&\label{w-conv}
\W_{\EPS k}\to \W_\EPS\hspace*{-0em}&&
\hspace*{-3em}\text{strongly in $L^s(I{\times}\varOmega)$ for $1\le s<2+4/d$}.
\intertext{Since $\OMEGA_\LAM(\FFepsk,\mm_{\EPS k},\cdot)$ is increasing, we can
write $\theta_{\EPS k}=[\OMEGA_\LAM(\FFepsk,\mm_{\EPS k},\cdot)]^{-1}(\W_{\EPS k})$. 
Thanks to the continuity of
$(\FF,\mm,\W)\mapsto[\OMEGA_\LAM(\FF,\mm,\cdot)]^{-1}(\W):\R^{d\times d}\times\R^d\times\R\to\R$
and the at most linear growth with respect to $\W$ uniformly
 with respect to $\FF$ from any compact $K\subset{\rm GL}^+(d)$,
cf.\ \eq{Euler-ass-adiab}, we have also}
&\label{z-conv}
\theta_{\EPS k}\to \theta_\EPS=[\OMEGA_\LAM(\FFeps,\mm_{\EPS k},\cdot)]^{-1}(\W_\EPS)
\hspace*{-0em}&&\hspace*{-1em}\text{strongly in $L^s(I{\times}\varOmega)$
for $1\le s<2+4/d$};
\intertext{actually, (\ref{Euler-weak+}a,b) results from interpolation of
\eq{Euler-quasistatic-est2}. Note that we do not have any direct information about
$\pdt{}\theta_{\EPS k}$
so that we could not use the Aubin-Lions arguments straight for
$\{ \theta_{\EPS k}\}_{k\in\N}$. Thus, by the continuity of the corresponding
Nemytski\u{\i} (or here simply superposition) mappings, also the 
conservative part of the regularized Cauchy stress as well as the heat part
of the internal energy, namely}
&
\TT_{\LAM,\EPS}(\FFepsk,\mm_{\EPS k},\theta_{\EPS k})\to
\TT_{\LAM,\EPS}(\FFeps,\mm_\EPS,\theta_\EPS)
\hspace*{-.5em}&&\hspace*{0em}\text{strongly in
$L^c(I{\times}\varOmega;\R_{\rm sym}^{d\times d}),\ \ 1\le c<\infty$,}
\label{Euler-T-strong-conv}
\\&\nonumber
\frac{\pi_\LAM(\FFepsk)\COUPLING'_{\Fe}(\FFepsk,\mm_{\EPS k},\theta_{\EPS k})\FFepsk^\top
}{(1{+}\EPS|\theta_{\EPS k}|)\det\FFepsk}
\\&\hspace{2em}
\to\frac{\pi_\LAM(\FFeps)\COUPLING'_{\Fe}(\FFeps,\mm_\EPS,\theta_\EPS)\FFeps^\top}
{(1{+}\EPS\theta_\EPS)\det\FFepsk}
\hspace*{0em}&&\hspace*{0em}\text{strongly in
$L^c(I{\times}\varOmega;\R_{\rm sym}^{d\times d}),\ \ 1\le c<\infty$,}
\\&\SS_{\LAM,\EPS}(\FFepsk,\mm_{\EPS k},\theta_{\EPS k})\to
\SS_{\LAM,\EPS}(\FF_{\EPS},\mm_{\EPS},\theta_{\EPS})\hspace*{-.5em}&&
\text{strongly in $L^2(I;L^2(\varOmega;\R^{d\times d}_{\rm skw})),\ 1\le c<\infty$},\!
\\\nonumber
&\!\!\mathscr{S}_{\LAM}(\FFepsk,\mm_{\EPS k},\!\nabla\mm_{\EPS k})
\\&\qquad\qquad\qquad
\to\!\mathscr{S}_{\LAM}(\FF_{\EPS},\mm_{\EPS},\!\nabla\mm_{\EPS})\!\!\!
&&\text{weakly* in $L^\infty(I;L^{2^*2/(2^*+2)}(\varOmega;\R^{d\times d\times d}))$},\!\!
\label{hyper-S-weak}
\\&\widehat\tt_\LAM(\FFepsk,\mm_{\EPS k},\theta_{\EPS k})\to\,
\widehat\tt_\LAM(\FFeps,\mm_\EPS,\theta_\EPS)\hspace*{-0em}&&\hspace*{0em}\text{strongly in
$L^c(I{\times}\varOmega;\R^d),\ \ 1\le c<2$},
\\&\OMEGA_\LAM(\FFepsk,\mm_{\EPS k},\theta_{\EPS k})\to
\OMEGA_\LAM(\FFeps,\mm_\EPS,\theta_\EPS)\hspace*{-0em}&&\hspace*{0em}\text{strongly in $L^c(I{\times}\varOmega),\ \ 1\le c<2+4/d$.}\label{Euler-weak+w}
\end{align}\end{subequations} 
It is important to notice that $\nabla(\varrho_{\EPS k}\vv_{\EPS k})=
\nabla\varrho_{\EPS k}{\otimes}\vv_{\EPS k}+\varrho_{\EPS k}\nabla\vv_{\EPS k}$
is bounded in $L^{\UUU \infty}(I;L^r(\varOmega;\R^{d\times d}))$ due to the already
obtained bounds (\ref{Euler-quasistatic-est1}a,c,d).
Therefore, $\varrho_{\EPS k}\vv_{\EPS k}$ converges weakly* in
$L^{\infty}(I;W^{1,r}(\varOmega;\R^d))$. The limit of
$\varrho_{\EPS k}\vv_{\EPS k}$ can be identified as $\varrho_\EPS\vv_\EPS$ because
we already showed that $\varrho_{\EPS k}$ converges strongly in \eq{rho-conv}
and $\vv_{\EPS k}$ converges weakly due to \eqref{Euler-weak-sln-v}. 

By comparison \eq{Euler1-weak-Galerkin}, we also obtain an information about
$\pdt{}(\varrho_{\EPS k}\vv_{\EPS k})$. Specifically, for any
$\widetilde\vv\in L^\infty(I;L^2(\varOmega;\R^d))\,\cap\, L^p(I;W^{2,p}(\varOmega;\R^d))$ with
$\widetilde\vv(t)\in V_k$ for a.a.\ $t\in I$, we have 
\begin{align}\nonumber
&\int_0^T\!\!\!\int_\varOmega
\pdt{}(\varrho_{\EPS k}\vv_{\EPS k}){\cdot}\widetilde\vv\,\d\xx\d t
=\!\int_0^T\!\!\!\int_\varGamma\!\!
\big(\TRACTION{+}\nu_\flat|\vv_{\EPS k}|^{p-2}\vv_{\EPS k}\big){\cdot}\widetilde\vv
\,\d S\d t+\int_0^T\!\!\!\int_\varOmega\bigg(
\sqrt{\frac{\rhoR\varrho_{\EPS k}}{\det_\LAM(\FFepsk)}}
\GRAVITY{\cdot}\widetilde\vv
\\&\hspace{0em}\nonumber
+\Big(\varrho_{\EPS k}\vv_{\EPS k}{\otimes}\vv_{\EPS k}
-\TT_{\LAM,\EPS}(\FFepsk,\mm_{\EPS k},\theta_{\EPS k})
-\KK_\LAM(\FFepsk,\nabla\mm_{\EPS k})
-\SS_{\LAM,\EPS}(\FFepsk,\mm_{\EPS k},\theta_{\EPS k})
\\[-.2em]&\hspace{0em}\nonumber
-\DIS|\ee(\vv_{\EPS k})|^{p-2}\ee(\vv_{\EPS k})\Big){:}\EE(\widetilde\vv)
-\Big(\nu_2|\nabla^2\vv_{\EPS k}|^{p-2}\nabla^2\vv_{\EPS k}{+}
\mathscr{S}_\LAM(\FFepsk,\mm_{\EPS k},\nabla\mm_{\EPS k})\!\Big)
\Vdots\nabla\EE(\widetilde\vv)
\\[-.2em]&\hspace{0em}
+\mu_0\nabla\hh_{\EPS k}{:}(\mm_{\EPS k}{\otimes}\wt\vv)
+\mu_0\hh_{\EPS k}{\cdot}\mm_{\EPS k}({\rm div}\,\widetilde\vv)\!
\bigg)\,\d\xx\d t\le C\|\widetilde\vv\|_{L^\infty(I;L^2(\varOmega;\R^d))\,\cap\,L^p(I;W^{2,p}(\varOmega;\R^d))}^{}
\label{est-of-DT-rho.v}
\end{align}
with $C$ independent of $k$. This yields a bound for
$\pdt{}(\varrho_{\EPS k}\vv_{\EPS k})$ in a seminorm on
$L^1(I;L^2(\varOmega;\R^d))$ $+L^{p'}(I;W^{2,p}(\varOmega;\R^d)^*)$ induced by
the Galerkin discretization by $V_k$, and by any $V_l$ with $l\le k$ with
$C$ in \eq{est-of-DT-rho.v} independent of $k$.
Here we used in particular that $\KK_\LAM(\FFepsk,\nabla\mm_{\EPS k})$
and $\SS_{\LAM,\EPS}(\FFepsk,\mm_{\EPS k},\theta_{\EPS k})$ are bounded in
$L^\infty(I;L^1(\varOmega;\R^{d\times d}))$  which is surely in duality
with $L^p(I;W^{1,p}(\varOmega;\R^{d\times d}))$ and that 
$\mathscr{S}_\LAM(\FFepsk,\mm_{\EPS k},\nabla\mm_{\EPS k})$ is bounded in
$L^\infty(I;L^{2^*2/(2^+2)}(\varOmega;\R^{d\times d\times d}))$ which is in duality
with $L^p(I{\times}\varOmega;\R^{d\times d\times d})$ if $p>d$.
By a generalization of the Aubin-Lions compact-embedding theorem,
cf.\ \cite[Lemma 7.7]{Roub13NPDE}, we then obtain 
\begin{subequations}\label{rho-v-conv}\begin{align}\label{rho.v-conv}
&\varrho_{\EPS k}\,\vv_{\EPS k}\to\varrho_{k}\vv_\EPS
&&\hspace*{-1em}\text{strongly in }L^{c}(I{\times}\varOmega;\R^d)\ \ \text{ for any $1\le c<4$}\,.
\intertext{Since obviously
$\vv_{\EPS k}=(\varrho_{\EPS k}\vv_{\EPS k})(1/\varrho_{\EPS k})$,
thanks to \eq{rho-conv} and \eq{rho.v-conv}, we also have that}
&\vv_{\EPS k}\to\vv_\EPS
&&\hspace*{-1em}\text{strongly in }L^c(I{\times}\varOmega;\R^d)\ \ \text{ with
any $1\le c<4$}\,.
\end{align}\end{subequations}

For the limit passage in the momentum equation, one uses the monotonicity
of the dissipative stress, i.e., the monotonicity of the quasilinear operator
$\vv\mapsto{\rm div}\big({\rm div}(\nu_2|\Nabla^2\vv|^{p-2}\nabla^2\vv)
-\DIS|\ee(\vv)|^{p-2}\ee(\vv)\big)$,
as well as of the time-derivative operator. One could use the already obtained
weak convergences and the so-called Minty trick but, later, we will need a
strong convergence of $\ee(\vv_{\EPS k})$ to pass to the limit in the heat
equation. Thus we first prove this strong convergence, which then allows for
the limit passage in the momentum equation directly.
We will use the weak convergence of the inertial force
\begin{align}
&
\int_0^T\!\!\!\int_\varOmega
\Big(\pdt{(\varrho_{\EPS k}\vv_{\EPS k})}+{\rm div}(\varrho_{\EPS k}\vv_{\EPS k}{\otimes}\vv_{\EPS k})
\Big){\cdot}\widetilde\vv\,\d \xx\d t
{\buildrel{k\to\infty}\over{\longrightarrow}}\!
\int_0^T\!\!\!\int_\varOmega
\Big(\pdt{(\varrho_\EPS\vv_\EPS)}+{\rm div}(\varrho_\EPS\vv_\EPS{\otimes}\vv_\EPS)\Big){\cdot}\widetilde\vv\,\d \xx\d t\,;
\label{conv-of-inirtia.v}\end{align}
cf.\ \cite{Roub22TVSE,RouSte22VESS}.
Further, relying on the calculus \eq{calculus-convective-in-F}, we will used the identity
\begin{align}\nonumber
\!\!\int_\varOmega\!\frac{\varrho_{\EPS k}(T)}2\big|\vv_{\EPS k}(T){-}\vv_\EPS(T)\big|^2\d \xx
&=\int_0^T\!\!\!\int_\varOmega
\Big(\pdt{(\varrho_{\EPS k}\vv_{\EPS k})}
+{\rm div}(\varrho_{\EPS k}\vv_{\EPS k}{\otimes}\vv_{\EPS k})
\Big){\cdot}\vv_{\EPS k}\,\d \xx\d t
\\[-.1em]&
+\!\int_\varOmega\!\frac{\varrho_0}2|\vv_0|^2\!-\varrho_{\EPS k}(T)\vv_{\EPS k}(T){\cdot}\vv_\EPS(T)
+\frac{\!\varrho_{\EPS k}(T)}2|\vv_\EPS(T)|^2\d \xx.
\label{Euler-one-substitution}\end{align}
We further used that the $\varrho_{\EPS k}(T) $ is
also bounded in $W^{1,r}(\varOmega)$ 
and $\vv_{\EPS k}(T)$ is bounded in $ L^2(\varOmega;\R^d)$, together
with some information about the time derivative 
$\pdt{}(\varrho_{\EPS k}\vv_{\EPS k})$, cf.\ \eq{est-of-DT-rho.v}, so that
we can identify the weak limit of $\varrho_{\EPS k}(T)\vv_{\EPS k}(T)$.
Specifically, we have that
\begin{align}
&&&\varrho_{\EPS k}(T)\vv_{\EPS k}(T)\to\varrho_\EPS(T)\vv_\EPS(T)&&
\text{weakly in $\ L^2(\varOmega;\R^d)$.}&&&&
\label{Euler-weak-rho.v(T)}\end{align}

For the stress $\KK_\LAM(\FFepsk,\nabla\mm_{\EPS k})$, we will further need
strong convergence of $\nabla\mm_{\EPS k}$. This can be seen from the uniform
monotonicity of the operator $\mm\mapsto-{\rm div}(\upkappa_\LAM(\FF)\nabla\mm/\det\FF)$.
We use the Galerkin approximation \eq{Euler2-weak-Galerkin} tested by
$\wt\rr=\mm_{\EPS k}{-}\widetilde\mm_k$ with an approximation $\widetilde\mm_k:I\to Z_k$
of $\mm_\EPS$ converging to $\mm_\EPS$ strongly in $L^2(I;H^1(\varOmega;\R^d))$.
Then we can estimate
\begin{align}\nonumber
  &c_{\upkappa,\LAM}\|\nabla\mm_{\EPS k}{-}\nabla\mm_\EPS\|_{L^2(I{\times}\varOmega)}^2\le
  \int_0^T\!\!\!\int_\varOmega
\frac{\upkappa_\LAM(\FFepsk)}{\det\FFepsk}|\nabla\mm_{\EPS k}{-}\nabla\mm_\EPS|^2
\,\d\xx\d t
\\&\nonumber=\int_0^T\!\!\!\int_\varOmega\!\Big(\wh\tt_{\LAM,\EPS}(\FFepsk,\mm_{\EPS k},\theta_{\EPS k})
-\hh_{\EPS k}+\tau\rr_{\EPS k}+\HC(\FFepsk,\theta_{\EPS k})\dd_{\EPS k}
-\frac{\mm_{\EPS k}{\times}\rr_{\EPS k}}{\gamma(\FFepsk,\mm_{\EPS k},\theta_{\EPS k})}\Big){\cdot}(\mm_{\EPS k}{-}\wt\mm_k)
\\&\nonumber\qquad\qquad\qquad\
+\frac{\upkappa_\LAM(\FFepsk)\nabla\mm_{\EPS k}\!}{\det\FFepsk}{:}\nabla(\mm_\EPS{-}\wt\mm_k)
-\frac{\upkappa_\LAM(\FFepsk)\nabla\mm_\EPS\!}{\det\FFepsk}{:}\nabla(\mm_{\EPS k}{-}\mm_\EPS)\,\d\xx\d t
\stackrel{k\to\infty}{\to}0\end{align}
with some $\dd_{\EPS k}\in{\rm Dir}(\rr_{\EPS k})$, where
$c_{\upkappa,\LAM}:=\inf_{F\in\R^{d\times d}}\upkappa_\LAM(F)/\!\det F$ is positive
thanks to our definition \eq{cut-off-kappa}. This convergence to 0 for
$k\to\infty$ is due to (\ref{Euler-weak-sln}d,e). Therefore,
\begin{subequations}\begin{align}
&&&\mm_{\EPS k}\to\mm_\EPS&&\text{strongly in $L^2(I;H^1(\varOmega;\R^d))$\,,}&&
\intertext{which also improves the weak*-convergence in \eq{hyper-S-weak}
and ensures convergence
$\KK_\LAM(\FFepsk,\nabla\mm_{\EPS k})$ $\to $ $\KK_\LAM(\FFeps,\nabla\mm_\EPS)$
in $L^1(I{\times}\varOmega;\R^{d\times d}_{\rm sym})$.
Moreover, by interpolation with \eq{Euler-weak-sln-m},}
&&&\nabla\mm_{\EPS k}\to\nabla\mm_\EPS&&
\text{strongly in $L^c(I;L^2(\varOmega;\R^{d\times d}))$ for any $1\le c<\infty$.}
\label{K-strongly}\end{align}\end{subequations}
Thus $|\nabla\mm_{\EPS k}|^2\to |\nabla\mm_{\EPS}|^2$
and $\nabla\mm_{\EPS k}{\otimes}\nabla\mm_{\EPS k}\to \nabla\mm_{\EPS}{\otimes}\nabla\mm_{\EPS}$
strongly in $L^c(I;L^1(\varOmega;\R^{d\times d}))$ for any $1\le c<\infty$, which is needed
for the convergence in \eq{Euler2-weak-Galerkin} when multiplied by
$\ee(\vv_{\EPS k})\in L^p(I;L^\infty(\varOmega;\R^{d\times d}_{\rm sym}))$.

We now use the Galerkin approximation of the regularized momentum equation
\eq{Euler1-weak-Galerkin} tested by $\widetilde\vv=\vv_{\EPS k}-\widetilde\vv_k$
with 
$\widetilde\vv_k:I\to V_k$ an approximation of $\vv_\EPS$ in the sense that
$\widetilde\vv_k\to\vv_\EPS$ strongly in $L^\infty(I;L^2(\varOmega;\R^d))$
and $\Nabla^2\widetilde\vv_k\to\Nabla^2\vv_\EPS$ for $k\to\infty$
strongly in $L^p(I{\times}\varOmega;\R^{d\times d\times d})$ for $k\to\infty$.
Using also the first inequality in \eq{Euler-quasistatic-est-formal4}
and \eq{Euler-one-substitution}, we can estimate \COMMENT{CHECK the indices and the reminder $O_k$ + magnetic stresses!!!}
\begin{align}\nonumber
  &\frac1{2C_{r,\EPS}\!\!}\big\|\vv_{\EPS k}(T){-}\vv_\EPS(T)\big\|_{L^2(\varOmega;\R^d)}^2\!\!
+\DIS c_{p}\|\EE(\vv_{\EPS k}{-}\vv_\EPS)\|_{L^p(I{\times}\varOmega;\R^{d\times d})}^p\!\!
+\nu_2 c_{p}\|\nabla^2(\vv_{\EPS k}{-}\vv_\EPS)\|_{L^p(I{\times}\varOmega;\R^{d\times d\times d})}^p
\\&\nonumber
\le\int_\varOmega\frac{\varrho_{\EPS k}(T)}2\big|\vv_{\EPS k}(T){-}\vv_\EPS(T)\big|^2\,\d\xx
+\int_0^T\!\!\!\int_\varGamma\nu_\flat|\vv_{\EPS k}{-}\vv_\EPS|^p\,\d S\d t
\\[-.3em]&\nonumber\hspace*{2em}
+\int_0^T\!\!\!\int_\varOmega\!\bigg(
\DIS\big(|\ee(\vv_{\EPS k})|^{p-2}\ee(\vv_{\EPS k})-|\ee(\vv_{\EPS})|^{p-2}\ee(\vv_{\EPS})\big){:}\ee(\vv_{\EPS k}{-}\vv_{\EPS})
  \\[-.6em]&\hspace*{8em}\nonumber
 +\nu_2\big(|\nabla^2\vv_{\EPS k}|^{p-2}\nabla^2\vv_{\EPS k}
-|\nabla^2\vv_\EPS|^{p-2}\nabla^2\vv_\EPS\big)\Vdots
  \nabla^2(\vv_{\EPS k}{-}\vv_\EPS)\bigg)\,\d\xx\d t
 \\[-.5em]&=\nonumber
 \int_0^T\!\!\!\int_\varOmega\bigg(\Big(\sqrt{\frac{\rhoR\varrho_{\EPS k}\!\!}{\det\FFepsk\!\!}}\ \GRAVITY+\mu_0(\nabla\hh_{\EPS k})^\top\mm_{\EPS k}\Big){\cdot}(\vv_{\EPS k}{-}\widetilde\vv_k)
+\mu_0\hh_{\EPS k}{\cdot}\mm_{\EPS k}{\rm div}(\vv_{\EPS k}{-}\widetilde\vv_k)
\\[.1em]&\nonumber\hspace{2em}
 -\big(\TT_{\LAM,\EPS}(\FFepsk,\mm_{\EPS k},\theta_{\EPS k}){+}
\KK_\LAM(\FFepsk,\nabla\mm_{\EPS k}){+}\SS_{\LAM,\EPS}(\FFepsk,\mm_{\EPS k},\theta_{\EPS k})
\big){:}\ee(\vv_{\EPS k}{-}\widetilde\vv_k)
 \\[.2em]&\nonumber\hspace{2em}
 -\DIS|\EE(\widetilde\vv_k)|^{p-2}\EE(\widetilde\vv_k){:}\EE(\vv_{\EPS k}{-}\widetilde\vv_k)
-\nu_2\big(|\nabla^2\widetilde\vv_k|^{p-2}\nabla^2\widetilde\vv_k\big)\Vdots
 \nabla^2(\vv_{\EPS k}{-}\widetilde\vv_k)
 \\[-.1em]&\nonumber\hspace*{2em}+
\Big(\pdt{}(\varrho_{\EPS k}\vv_{\EPS k})
+{\rm div}(\varrho_{\EPS k}\vv_{\EPS k}{\otimes}\vv_{\EPS k})\Big)
{\cdot}\widetilde\vv_k
-\mathscr{S}_\LAM(\FFepsk,\mm_{\EPS k},\nabla\mm_{\EPS k})\Vdots\nabla^2(\vv_{\EPS k}{-}\widetilde\vv_k)
\bigg)\,\d\xx\d t
 \\[-.1em]&\nonumber\hspace{2em}
+\int_0^T\!\!\!\int_\varGamma\big(\TRACTION{+}
\nu_\flat|\widetilde\vv|^{p-2}\widetilde\vv_k\big)
{\cdot}(\vv_{\EPS k}{-}\widetilde\vv_k)\,\d S\d t
+\int_\varOmega\bigg(\frac{\varrho_0}2|\vv_0|^2-\varrho_{\EPS k}(T)\vv_{\EPS k}(T){\cdot}\widetilde\vv_k(T)
\\[-.3em]&\hspace*{20.5em}
+\frac{\varrho_{\EPS k}(T)}2|\widetilde\vv_k(T)|^2\bigg)\,\d\xx+\mathscr{O}_k
\ {\buildrel{k\to\infty}\over{\longrightarrow}}\ 0\,
\label{strong-hyper+}\end{align}
with $C_{r,\EPS}>0$ from \eq{Euler-est1-3} and with some $c_{p}>0$ related to the
inequality
$c_{p}|G-\widetilde G|^p\le(|G|^{p-2}G-|\widetilde G|^{p-2}\widetilde G)\Vdots(G-\widetilde G)$ holding for $p\ge2$. The remainder term $\mathscr{O}_k$ in \eq{strong-hyper+} is
\begin{align}\nonumber
\mathscr{O}_k&=\int_\varOmega\frac{\varrho_{\EPS k}(T)\!}2\,\vv_{\EPS k}(T)
{\cdot}\big(\widetilde\vv_k(T){-}\vv_\EPS(T)\big)\,\d\xx
+\int_0^T\!\!\!\int_\varGamma\nu_\flat|\vv_{\EPS k}|^{p-2}\vv_{\EPS k}{\cdot}(\widetilde\vv_k{-}\vv_\EPS)\,\d S\d t
\\[-.3em]&\nonumber\hspace*{1em}
+\int_0^T\!\!\!\int_\varOmega\!
\DIS|\EE(\vv_{\EPS k})|^{p-2}\EE(\vv_{\EPS k}){:}\EE(\widetilde\vv_k{-}\vv_\EPS)
 +\nu_2|\nabla^2\vv_{\EPS k}|^{p-2}\nabla^2\vv_{\EPS k}\Vdots
  \nabla^2(\widetilde\vv_k{-}\vv_\EPS)\,\d\xx\d t
\end{align}
and it converges to zero due to the strong approximation properties of the
approximation $\widetilde\vv_k$ of $\vv_\EPS$. Here we used
\eq{conv-of-inirtia.v}--\eq{Euler-one-substitution} and also the strong convergence
\eq{rho-conv}, \eq{z-conv}, and \eq{rho-v-conv}. Knowing already
\eq{Euler-T-strong-conv} and that $\ee(\vv_{\EPS k}{-}\widetilde\vv_k)\to0$ weakly
in $L^p(I;W^{1,p}(\varOmega;\R_{\rm sym}^{d\times d}))$, we have that
$\int_0^T\!\int_\varOmega\TT_{\LAM,\EPS}(\FFepsk,\mm_{\EPS k},\theta_{\EPS k})
{:}\ee(\vv_{\EPS k}{-}\widetilde\vv_k)\,\d\xx\d t\to0$. Thus we obtain the desired
strong convergence
\begin{subequations}\label{strong-conv}\begin{align}
&\vv_{\EPS k}\to\vv_{\EPS}&&\text{strongly in
$L^p(I;W^{2,p}(\varOmega;\R^d))$}
\intertext{and also of $\vv_{\EPS k}(T)\to\vv_{\EPS}(T)$ in $L^2(\varOmega;\R^d)$.
In fact, executing this procedure for
a current time instants $t$ instead of $T$, we obtain}
&\vv_{\EPS k}(t)\to\vv_\EPS(t)&&\text{strongly in $\,L^2(\varOmega;\R^d)\,$ for any $t\in I$.}&&
\intertext{By \eq{Euler-weak-sln-m} and the Aubin-Lions theorem, we also obtain}
&\mm_{\EPS k}\to\mm_\EPS&&\text{strongly in $\,L^c(I;L^{2^*-1/c}(\varOmega;\R^d))\,$ for any
$1\le c<\infty$.}&&
\intertext{It also implies, by continuity of the trace operator
$L^p(I;W^{2,p}(\varOmega))\to L^p(I{\times}\varGamma)$, that}
&\vv_{\EPS k}\big|_{I\times\varGamma}\to\vv_\EPS\big|_{I\times\varGamma}&&
\text{strongly in $\,L^p(I{\times}\varGamma;\R^d)$}\,.
\end{align}\end{subequations}
Having \eq{strong-conv} at disposal, the limit passage in the
Galerkin-approximation of \eq{Euler1-weak-Galerkin} to the weak solution
of \eq{Euler-thermo-reg1} is then easy. The variational inequality
\eq{Euler2-weak-Galerkin} can be converged by lower weak-semicontinuity of
its right-hand-side integral functionals. Convergence in
(\ref{Euler-weak-Galerkin}a--c) is due to Lemmas~\ref{lem1} and \ref{lem2}.

For further purposes, let us mention that the energy dissipation
balance \eq{thermodynamic-Euler-mech-disc} is inherited in the limit, i.e.
\begin{align}\nonumber
  &\hspace*{0em}\frac{\d}{\d t}
  \int_\varOmega\!\frac{\varrho_\EPS}2|\vv_\EPS|^2+
  \frac{\pi_\LAM(\FFeps)\upvarphi(\FFeps,\mm_\EPS)\!}{\det\FFeps}
  +\frac{\upkappa_\LAM(\FFeps)}{2\det\FF}|\Nabla\mm_\EPS|^2
  -\mu_0\hh_\text{\rm ext}{\cdot}\mm_\EPS\,\d\xx
\\[-.1em]&\hspace{.0em}\nonumber
\qquad+\!\int_\varOmega\!\xi_\EPS(\FFeps,\theta_\EPS;\ee(\vv_\EPS),\nabla^2\vv_\EPS,\rr_\EPS)
\,\d\xx
+\!\int_\varGamma\!\nu_\flat|\vv_\EPS|^p\,\d S
+\frac{\d}{\d t}\int_{\R^d}\frac{\mu_0}2|\Nabla u_\EPS|^2\,\d\xx
\\[-.1em]&\nonumber\hspace{.0em}
=\int_\varGamma\!\TRACTION{\cdot}\vv_\EPS\,\d S
+\int_\varOmega\bigg(\sqrt{\frac{\varrho_\EPS\rhoR}{\det_\LAM(\FFeps)\!}}\:\GRAVITY{\cdot}\vv_\EPS
-\frac{\pi_\LAM(\FFeps)\COUPLING_{\Fe}'(\FFeps,\mm_\EPS,\theta_\EPS)
\FFeps^\top}{(1{+}\EPS\theta_\EPS)\det\FFeps}{:}\ee(\vv_\EPS)
\\[-.1em]&\hspace{.0em}
\qquad\qquad-\frac{\pi_\LAM(\FFeps)\COUPLING_\mm'(\FFeps,\mm_\EPS,\theta_\EPS)\!}
{(1{+}\EPS\theta_\EPS^{1/2})\det\FFeps}
{\cdot}\big(\rr_\EPS{-}{\rm skw}(\nabla\vv_\EPS)\mm_\EPS\big)
-\mu_0\frac{\!\partial \hh_\text{\rm ext}\!}{\partial t}{\cdot}\mm_\EPS\!\bigg)\,\d\xx
\label{thermodynamic-Euler-mech+}
\end{align}
with $u_\EPS=\Delta^{-1}{\rm div}(\mm_\EPS)$. This follows by legitimacy of the
tests which led formally (in the non-regularized situation) to \eq{energy+}.
Specifically, here (\ref{Euler-thermo-reg}a--d) in the weak formulation are
in duality (and can be tested) by $|\vv_\EPS|^2/2$, $\vv_\EPS$,
$[\upvarphi(\FFeps,\mm_\EPS)/\!\det\FFeps]_{\!\FF}'$, and $\rr_\EPS$, cf.\ Step~10 below.

\medskip\noindent{\it Step 5: Strong convergence of $\rr_{\EPS k}$ and
limit passage in the thermal part for $k\to\infty$}.
For the heat equation, we will still need the strong convergence of
$\rr_{\EPS k}\to\rr_\EPS$. We cannot exploit the strong monotonicity of the operator
$\rr\mapsto\tau\rr+\HC(\FF,\theta){\rm dir}_\EPS(\rr)
$ because the exchange driving force on the right-hand side of \eq{Euler-thermo-reg4}
is not a compact lower-order term. Thus, as we already passed to the limit
in the magneto-mechanical part, we can use the ``limsup-trick'' and the
strict convexity of the potential of the mentioned operator. Specifically, 
\begin{align}\nonumber
&\int_0^T\!\!\bigg(\int_\varOmega
\xi_\EPS\big(\FFeps,\theta_\EPS;\EE(\vv_\EPS),\nabla^2\vv_\EPS,\rr_\EPS\big)\,\d\xx
+\int_\varGamma\nu_\flat|\vv_\EPS|^p\d S\bigg)\d t
\\&\nonumber
\le\liminf_{k\to\infty}\int_0^T\!\!\bigg(\int_\varOmega\xi_\EPS\big(\FFepsk,\theta_{\EPS k};\EE(\vv_{\EPS k}),\nabla^2\vv_{\EPS k},\rr_{\EPS k}\big)\,\d\xx
+\int_\varGamma\nu_\flat|\vv_{\EPS k}|^p\d S\bigg)\d t
\\&\nonumber
\le\limsup_{k\to\infty}\int_0^T\!\!\bigg(\int_\varOmega\xi_\EPS\big(\FFepsk,\theta_{\EPS k};\EE(\vv_{\EPS k}),\nabla^2\vv_{\EPS k},\rr_{\EPS k}\big)\,\d\xx
+\int_\varGamma\nu_\flat|\vv_{\EPS k}|^p\d S\bigg)\d t
\\&\nonumber
=\int_\varOmega\frac{\varrho_0}2|\vv_0|^2+\frac{\upvarphi(\FF_0,\mm_0)}{\det\FF_0}
+\frac{\upkappa(\FF_0)|\Nabla\mm_0|^2}{2\det\FF_0}
-\mu_0\hh(0){\cdot}\mm_0\,\d\xx
+\int_{\R^d}\frac{\mu_0}2|\Nabla u(0)|^2\,\d\xx
\\&\nonumber\qquad
-\liminf_{k\to\infty}
\bigg(\int_\varOmega\!\bigg(\frac{\varrho_{\EPS k}(T)}2|\vv_{\EPS k}(T)|^2\!+\frac{\pi_\LAM(\FFepsk(T))\upvarphi(\FFepsk(T),\mm_{\EPS k}(T))\!}{\det\FFepsk(T)}
\\&\nonumber\qquad\qquad\
+\frac{\upkappa_\LAM(\FFepsk(T))|\Nabla\mm_{\EPS k}(T)|^2\!}{2\det\FFepsk(T)}
-\mu_0\hh(T){\cdot}\mm_{\EPS k}(T)\,\d\xx
+\int_{\R^d}\frac{\mu_0}2|\Nabla u_{\EPS k}(T)|^2\bigg)\,\d\xx\bigg)
\\&\nonumber\qquad
+\lim_{k\to\infty}\int_0^T\!\!\bigg(
\int_\varOmega\varrho_{\EPS k}\,\GRAVITY{\cdot}\vv_{\EPS k}+\mu_0\frac{\!\partial \hh_\text{\rm ext}\!}{\partial t}{\cdot}\mm_{\EPS k}
-\pi_\LAM(\FFepsk)\frac{\COUPLING_{\Fe}'(\FFepsk,\mm_{\EPS k},\theta_{\EPS k})\FFepsk^\top\!}{(1{+}\EPS|\theta_{\EPS k}|^{\alpha})\det\FFepsk}{:}\ee(\vv_{\EPS k})
\\&\nonumber\qquad\qquad\
-\frac{\pi_\LAM(\FFepsk)\COUPLING_\mm'(\FFepsk,\mm_{\EPS k},\theta_{\EPS k})}{(1{+}\EPS|\theta_{\EPS k}|^{\alpha/2})\det\FFepsk}{\cdot}(\rr_{\EPS k}{-}{\rm skw}(\nabla\vv_{\EPS k})\mm_{\EPS k})
+\int_\varGamma\TRACTION{\cdot}\vv_{\EPS k}\d S\bigg)\d t
\\&\nonumber
\le\int_\varOmega\frac{\varrho_0}2|\vv_0|^2+\frac{\upvarphi(\FF_0,\mm_0)}{\det\FF_0}
+\frac{\upkappa(\FF_0)|\Nabla\mm_0|^2}{2\det\FF_0}
-\mu_0\hh(0){\cdot}\mm_0\,\d\xx
+\int_{\R^d}\frac{\mu_0}2|\Nabla u(0)|^2\,\d\xx
\\&\nonumber\qquad
-\int_\varOmega\!\bigg(\frac{\varrho_\EPS(T)}2|\vv_\EPS(T)|^2\!+\frac{\pi_\LAM(\FFeps(T))\upvarphi(\FFeps(T),\mm_\EPS(T))\!}{\det\FFeps(T)}
\\[-.3em]&\nonumber\qquad\qquad\
+\frac{\upkappa_\LAM(\FFeps(T))|\Nabla\mm_\EPS(T)|^2\!}{2\det\FFeps(T)}
-\mu_0\hh(T){\cdot}\mm_\EPS(T)\bigg)\,\d\xx+\int_{\R^d}\frac{\mu_0}2|\Nabla u_\EPS(T)|^2\,\d\xx
\\&\nonumber\qquad
+\int_0^T\!\!\bigg(
\int_\varOmega\varrho_\EPS\,\GRAVITY{\cdot}\vv_\EPS+\mu_0\frac{\!\partial \hh_\text{\rm ext}\!}{\partial t}{\cdot}\mm_\EPS
-\pi_\LAM(\FFeps)\frac{\COUPLING_{\Fe}'(\FFeps,\mm_\EPS,\theta_\EPS)\FFeps^\top\!}{(1{+}\EPS\theta_\EPS
)\det\FFeps}{:}\ee(\vv_\EPS)
\\&\nonumber\qquad\qquad\
-\frac{\pi_\LAM(\FFeps)\COUPLING_\mm'(\FFeps,\mm_\EPS,\theta_\EPS)}{(1{+}\EPS\theta_\EPS^{1/2}
)\det\FFeps}{\cdot}(\rr_\EPS{-}{\rm skw}(\nabla\vv_\EPS)\mm_\EPS)
+\int_\varGamma\TRACTION{\cdot}\vv_\EPS\d S\bigg)\d t
\\&\stackrel{\eq{thermodynamic-Euler-mech+}}=\!\int_0^T\!\!\bigg(\int_\varOmega
\xi_\EPS\big(\FFeps,\theta_\EPS;\EE(\vv_\EPS),\nabla^2\vv_\EPS,\rr_\EPS\big)\,\d\xx
+\int_\varGamma\nu_\flat|\vv_\EPS|^p\d S\bigg)\d t\,.
\label{limsup}\end{align}
This still gives
\begin{align}
&\rr_{\EPS k}\to\rr_\EPS&&\text{strongly in $L^2(I{\times}\varOmega;\R^d)$.}&&
\end{align}

Strong convergence (\ref{strong-conv}a,b,d) 
allows for passing to the limit also in the dissipative heat sources
and the other terms in the Galerkin-approximation 
of the heat equation \eq{Euler3-weak-Galerkin} are even easier.

\medskip\noindent{\it Step 6 -- non-negativity of temperature}:
We can now perform various nonlinear tests of the regularized but
non-discretized heat equation. The first test can be by the negative part
of temperature $\theta_\EPS^-:=\min(0,\theta_\EPS)$. Let us recall the
extension \eq{extension-negative}, which in particular gives
$\OMEGA(\FF,\mm,\theta^-)=\theta^-$ and $\OMEGA_\FF'(\FF,\mm,\theta^-)=\bm0$
and also $\OMEGA_\mm'(\FF,\mm,\theta^-)=\bm0$. Note also that
$\theta_\EPS^-\in L^2(I;H^1(\varOmega))$,
so that it is indeed a legal test for \eq{Euler-thermo-reg3}.
Here we rely on the data qualification $\DIS,\nu_2,\nu_\flat\ge0$,
$\COND=\COND(\FF,\theta)\ge0$, $\theta_0\ge0$, and $h(\theta)\ge0$ for $\theta\le0$, cf.\
(\ref{Euler-ass}f,i,n).
Realizing that $\nabla\theta^-\!=0$ wherever $\theta>0$ so that
$\nabla\theta{\cdot}\nabla\theta^-=|\nabla\theta^-|^2$ and that
$\COUPLING'_{\Fe}(\Fe,\mm,\theta)\theta^-=\COUPLING'_{\Fe}(\Fe,\mm,\theta^-)\theta^-=\bm0$
and $\COUPLING'_\mm(\Fe,\mm,\theta)\theta^-=\bm0$
and also $h(\theta)\theta^-=h(\theta^-)\theta^-=0$, this test gives
\begin{align}\nonumber
&\frac12\frac{\d}{\d t}\|\theta_\EPS^-\|_{L^2(\varOmega)}^2
\le\int_\varOmega\theta_\EPS^-\pdt{\W_\EPS}+
\COND(\FFeps,\theta_\EPS)\nabla\theta_\EPS{\cdot}\nabla\theta_\EPS^-\,\d\xx
\\[-.1em]&\nonumber
=\!\int_\varOmega\bigg(\!\W_\EPS\vv_\EPS{\cdot}\nabla\theta_\EPS^-\!
+\Big(\xi_\EPS\big(\FFeps,\theta_\EPS;\EE(\vv_\EPS),\nabla^2\vv_\EPS,\rr_\EPS\big)
+\frac{\pi_\LAM(\FFeps)\COUPLING'_{\Fe}(\FFeps,\mm_\EPS,\theta_\EPS)\FFeps^\top\!\!}
{(1{+}\EPS|\theta_\EPS|)\det\FFeps}{:}\ee(\vv_{\EPS k})
\\[-.1em]&\nonumber\ \ 
+\frac{\pi_\LAM(\FFeps)\COUPLING'_\mm(\FFeps,\mm_\EPS,\theta_\EPS)}
{(1{+}\EPS|\theta_\EPS|^{1/2})\det\FFeps}{\cdot}\big(\rr_\EPS{-}{\rm skw}(\nabla\vv_\EPS)\mm_\EPS\big)
\Big)\theta_\EPS^-\bigg)\,\d\xx
+\int_\varGamma\Big(h(\theta_\EPS){+}\frac{\nu_\flat|\vv_\EPS|^p}{2{+}\EPS|\vv_\EPS|^p}\Big)\theta_\EPS^-\,\d S
\\[-.1em]\nonumber
&\le\!\int_\varOmega\!\W_\EPS\vv_\EPS{\cdot}\nabla\theta_\EPS^-\,\d\xx
=\!\int_\varOmega\!\theta_\EPS^-\vv_\EPS{\cdot}\nabla\theta_\EPS^-\,\d\xx
=\!-\int_\varOmega\!|\nabla\theta_\EPS^-|^2{\rm div}\,\vv_\EPS\,\d\xx
\\[-.1em]
&=-\frac12\int_\varOmega|\nabla\theta_\EPS^-|^2{\rm div}\,\vv_\EPS\,\d\xx
\le\|\theta_\EPS^-\|_{L^2(\varOmega)}^2\|{\rm div}\,\vv_\EPS\|_{L^\infty(\varOmega)}^{}\,.
\label{Euler-thermo-test-nonnegative}
\end{align}
Recalling the assumption $\theta_0\ge0$ so that $\theta_{0,\EPS}^-=0$
and exploiting the information
$\vv_\EPS\in L^p(I;W^{1,p}(\varOmega;\R^d))$ with $p>d$ inherited
from \eq{Euler-quasistatic-est1-1}, by the Gronwall inequality we obtain
$\|\theta_\EPS^-\|_{L^\infty(I;L^2(\varOmega))}=0$, so that $\theta_\EPS\ge0$
a.e.\ on $I{\times}\varOmega$. 

Having proved non-negativity of temperature, we can now execute
the strategy based of the $L^1$-theory for the heat equation
which led to the estimates \eq{est-e(v)}--\eq{est+}, i.e.\ here
\begin{subequations}\label{est-eps}\begin{align}
&\|\vv_\EPS\|_{L^\infty(I;L^2(\varOmega;\R^d))\,\cap\,L^p(I;W^{2,p}(\varOmega;\R^d))}^{}\le C,\ \ \
\|\mm_\EPS\|_{L^\infty(I;H^1(\varOmega;\R^d))}^{}\le C,\label{est-e(v)-eps}
\\&\label{est+Fes-eps}\|\FFeps\|_{L^\infty(I;W^{1,r}(\varOmega;\R^{d\times d}))}\le C_r\,,
\ \ \ \:\Big\|\frac1{\det\FFeps}\Big\|_{L^\infty(I;W^{1,r}(\varOmega))}\le C_r\,,
\\&\label{est+rho-eps}\|\varrho_\EPS\|_{L^\infty(I;W^{1,r}(\varOmega))}^{}\le C_r\,,
\ \ \ \ \Big\|\frac1{\varrho_\EPS}\Big\|_{L^\infty(I;W^{1,r}(\varOmega))}\!\le C_r
\ \ \ \text{ for any $1\le r<+\infty$},
\\&\|\W_\EPS\|_{L^\infty(I;L^1(\varOmega))}^{}\le C\,,\ \ \text{ and }\ \ 
\|\theta_\EPS\|_{L^\infty(I;L^1(\varOmega))}^{}\le C\,.
\label{basic-est-of-theta-eps}
\intertext{By interpolation exploiting the Gagliardo-Nirenberg inequality
between $L^2(\varOmega)$ and $W^{2,p}(\varOmega)$, we have
$\|\cdot\|_{L^\infty(\varOmega)}^{}\le C\|\cdot\|_{L^2(\varOmega)}^{r}
\|\cdot\|_{W^{2,p}(\varOmega)}^{1-r}$ with $0<r<pd/(pd+4p-2d)$.
Using also Korn's inequality, from \eq{est-e(v)-eps}
we thus obtain the estimate}
&\label{est+v}
\|\vv_\EPS\|_{L^s(I;L^\infty(\varOmega;\R^d))}^{}\le C_s\ \ \ \text{ with }\
1\le s<\frac{p(pd{+}4p{-}2d)}{4p-2d}\,.
\intertext{By comparison from
$\pdt{}\varrho_\EPS=({\rm div}\vv_\EPS)\varrho_\EPS-\vv_\EPS{\cdot}\nabla\varrho_\EPS$, from 
$\pdt{}\FFeps=(\nabla\vv_\EPS)\FFeps-(\vv_\EPS{\cdot}\nabla)\FFeps$, and from
$\pdt{}\mm_\EPS=\rr_\EPS+{\rm skw}(\nabla\vv_\EPS)\mm_\EPS-(\vv_\EPS{\cdot}\nabla)\mm_\EPS$, we also have}
&\label{est+dF/dt}
\Big\|\pdt{\varrho_\EPS}\Big\|_{L^p(I;L^r(\varOmega))}^{}\!\le C\,,\ \ \ \ 
\Big\|\pdt{\FFeps}\Big\|_{L^p(I;L^r(\varOmega;\R^{d\times d}))}^{}\!\le C\,,
\ \text{ and }\ \Big\|\pdt{\mm_\EPS}\Big\|_{L^2(I{\times}\varOmega;\R^d)}^{}\!\le C
\,.
\end{align}\end{subequations}
The estimates \eq{basic-est-of-theta-eps} are naturally weaker than
\eq{Euler-quasistatic-est2} but, importantly,
are uniform with respect to $\EPS>0$, in contrast to
\eq{Euler-quasistatic-est2} which is not uniform in this sense.
The total energy balance \eq{thermodynamic-Euler-engr} holds for
$\EPS$-solution only as an inequality because the heat sources do not
exactly cancel; more in detail, while the regularized adiabatic heat again
cancels, the dissipative heat terms are regularized (and smaller) in
\eq{Euler-thermo-reg3} and in \eq{Euler-thermo-reg-BC-IC-2} but the
corresponding viscous stress in \eq{Euler-thermo-reg1} and force in
\eq{Euler-thermo-reg-BC-IC-1} are not regularized. This inequality still
allows to execute the above mentioned estimation.

Let us also note that the extension \eq{extension-negative} becomes now
inactive and we can work with the original data defined for
non-negative $\theta$ only.

\medskip\noindent{\it Step 7 -- further a-priori estimates}:
Furthermore, having $\rr_\EPS$ estimated in $L^2(I{\times}\varOmega;\R^d)$
uniformly with respect to $\EPS$, as in \eq{Delta-m-} we have
\begin{align}\nonumber
&\hspace{0em}\Delta\mm_\EPS=\frac{\det\FFeps}{\upkappa_\LAM(\FFeps)}\bigg(
\tau\rr_\EPS+\HC(\FFeps,\theta_\EPS){\rm dir}_\EPS(\rr_\EPS)
-\frac{\mm_\EPS{\times}\rr_\EPS}{\gamma(\FFeps,\mm_\EPS,\theta_\EPS)}-\hh_{\rm ext}-\nabla u_\EPS
\\&\hspace{3em}\nonumber
+\frac{\upvarphi_\mm'(\FFeps,\mm_\EPS)+\COUPLING_\mm'(\FFeps,\mm_\EPS,\theta_\EPS)\!}{\det\FFeps}
-\Big(\frac{\upkappa_\LAM'(\FFeps)}{\det\FFeps}-\frac{\upkappa_\LAM(\FFeps){\rm Cof}\FFeps}{\det\FFeps^2}\Big){\Vdots}(\nabla\FFeps{\otimes}\nabla\mm_\EPS)\bigg)\,,
\end{align}
so that we can estimate $\nabla^2\mm_\EPS$ by a $H^1$-regularity as
\eq{Delta-m+}, i.e.\ now
\begin{align}\label{Delta-m+eps}
&\|\nabla^2\mm_\EPS\|_{L^2(I\times\varOmega;\R^{d\times d\times d})}^{}\le C\,.
\end{align}

Furthermore,
we are to prove an estimate of $\nabla\theta_\EPS$ based on the test of the heat
equation \eq{Euler-thermo-reg3} by $\chi_\zeta(\theta_\EPS)$ with an increasing
nonlinear function $\chi_\zeta:[0,+\infty)\to[0,1]$ defined as
\begin{align}\label{test-chi}
\chi_\zeta(\theta):=1-\frac1{(1{+}\theta)^\zeta}\,,\ \ \ \ \zeta>0\,,
\end{align}
simplifying the original idea of L.\,Boccardo and T.\,Gallou\"et 
\cite{BDGO97NPDE,BocGal89NEPE} in the spirit of \cite{FeiMal06NSET},
expanding the estimation strategy in \cite[Sect.\,8.2]{KruRou19MMCM}.
Importantly, here we have $\chi_\zeta(\theta_\EPS(t,\cdot))\in H^1(\varOmega)$,
hence it is a legal test function, because 
$0\le\theta_\varepsilon(t,\cdot)\in H^1(\varOmega)$ has already been proved
and because $\chi_\zeta$ is Lipschitz continuous on $[0,+\infty)$. 

We consider $1\le \EXP<2$ and estimate the $L^\EXP$-norm  of $\nabla\theta_\varepsilon$
by H\"older's inequality as 
\begin{align}\nonumber
&\!\!\int_0^T\!\!\!\int_\varOmega|\nabla\theta_\varepsilon|^\EXP\,\d\xx\d t
\le C_1\bigg(\underbrace{\int_0^T
\big\|1{+}\theta_\varepsilon(t,\cdot)\big\|^{(1+\zeta)\EXP/(\TWO-\EXP)}
_{L^{(1+\zeta)\EXP/(\TWO-\EXP)}(\varOmega)}\,\d t}_{\displaystyle\ \ \ =:I_{\EXP,\zeta}^{(1)}(\theta_\varepsilon)}\bigg)^{1-\EXP/\TWO}
\bigg(\underbrace{\int_0^T\!\!\!\int_\varOmega\chi_\zeta'(\theta_\EPS)|\nabla\theta_\EPS|^\TWO}
_{\displaystyle\ \ \ =:I_{\zeta}^{(2)}(\theta_\varepsilon)}\bigg)^{\EXP/\TWO}.
\\[-2em]\label{8-**-+}
\end{align}
with $\chi_\zeta$ from \eq{test-chi} so that $\chi_\zeta'(\theta)=\zeta/(1{+}\theta)^{1+\zeta}$
and with a constant $C_1$ dependent on $\zeta$, $\EXP$, and $T$.
Then we interpolate the Lebesgue space $L^{(1+\zeta)\EXP/(\TWO-\EXP)}(\varOmega)$
between  $W^{1,\EXP}(\varOmega)$ and $L^1(\varOmega)$ in order to exploit the already obtained
$L^\infty(I;L^1(\varOmega))$-estimate  in \eq{basic-est-of-theta-eps}.
More specifically, by the Gagliardo-Nirenberg inequality, we obtain
\begin{align}
\big\|1{+}\theta_\varepsilon(t,\cdot)\big\|_{L^{\EXP/\sigma}(\varOmega)}^{\mu/\sigma}
\le C_2\Big(1+\big\|\nabla\theta_\varepsilon(t,\cdot)\big\|_{L^\EXP(\varOmega;\R^d)}\Big)^\EXP
\ \ \ \text{ with }\ \sigma=\frac{\TWO{-}\EXP}{1{+}\zeta}
\label{8-cond}
\end{align}
with $C_2$ depending on $\sigma$, $C_1$, and $C$ from \eq{basic-est-of-theta-eps},
so that $I_{\EXP,\zeta}^{(1)}(\theta_\varepsilon)\le C_3(1+\int_0^T\!\int_\varOmega\big|\nabla \theta_\varepsilon\big|^\EXP\,\d\xx\d t)$
with $C_3$ depending on $C_2$. Combining it with \eq{8-**-+}, we obtain
\begin{align}\|\nabla\theta_\varepsilon\|_{L^\EXP(I\times\Omega;\R^d)}^\EXP=C_1C_3\big(1+\|\nabla\theta_\varepsilon\|_{L^\EXP(I\times\Omega)}^\EXP\big)^{1-\EXP/2}I_{\EXP,\zeta}^{(2)}(\theta_\varepsilon)^{\EXP/2_{_{_{}}}}_{^{^{^{}}}}\,.
\label{8-***}
\end{align}
Furthermore, we estimate $I_{\zeta}^{(2)}(\theta_\varepsilon)$ in \eq{8-**-+}. Let us
denote by ${\cal X}_\zeta$ a primitive function to
$\theta\mapsto\chi_\zeta(\theta)\OMEGA_\theta'(\Fe,\mm,\theta)$ depending smoothly
on $\Fe$, specifically
\begin{align}
{\cal X}_\zeta(\Fe,\mm,\theta)
=\int_0^1\!\!\theta\chi_\zeta(r\theta)\OMEGA_\theta'(\Fe,\mm,r\theta)\,\d r\,.
\label{primitive+}\end{align}
Like \eq{Euler-thermodynam3-test++} but using partial
(not convective) time derivative, we have now the calculus
\begin{align}\nonumber
&\!\!\!\int_\varOmega\!\chi_\zeta(\theta)\pdt\W\,\d\xx
=\!\int_\varOmega\!\chi_\zeta(\theta)\OMEGA_\theta'(\Fe,\mm,\theta)\pdt{\theta}
+\chi_\zeta(\theta)\Big(\OMEGA_\Fe'(\Fe,\mm,\theta){:}\pdt{\Fe\!}\,
{+}\,\OMEGA_\mm'(\Fe,\mm,\theta){\cdot}\pdt{\mm}\Big)\,\d\xx
\\&\quad\nonumber
=\frac{\d}{\d t}\int_\varOmega{\cal X}_\zeta(\Fe,\mm,\theta)\,\d\xx
-\int_\varOmega\big[{\mathscr X}_\zeta\big]_\Fe'(\FF,\mm,\theta){:}\pdt{\Fe}
+\big[{\mathscr X}_\zeta\big]_\mm'(\FF,\mm,\theta){\cdot}\pdt{\mm}\,\d\xx
\\[-.0em]&\hspace{11.5em}\text{where }\
{\mathscr X}_\zeta(\FF,\mm,\theta):={\cal X}_\zeta(\Fe,\mm,\theta)
-\chi_\zeta(\theta)\OMEGA(\Fe,\mm,\theta)\,.
\label{Euler-thermodynam3-test-}\end{align}
In view of \eq{primitive+}, it holds $[{\mathscr X}_\zeta]_\Fe'(\Fe,\mm,\theta)
=\int_0^1\theta\chi_\zeta(r\theta)\OMEGA_{\Fe\theta}''(\Fe,\mm,r\theta)\,\d r
-\chi_\zeta(\theta)\OMEGA_{\Fe}'(\Fe,\mm,\theta)$
and $[{\mathscr X}_\zeta]_\mm'(\Fe,\mm,\theta)
=\int_0^1\theta\chi_\zeta(r\theta)\OMEGA_{\mm\theta}''(\Fe,\mm,r\theta)\,\d r
-\chi_\zeta(\theta)\OMEGA_{\mm}'(\Fe,\mm,\theta)$.
Altogether, testing \eq{Euler-thermo-reg3} with
\eq{Euler-thermo-reg-BC-IC-2} by $\chi_\zeta(\theta_\EPS)$ gives
\begin{align}\nonumber
&\frac{\d}{\d t}\int_\varOmega\!{\cal X}_\zeta(\FFeps,\mm_\EPS,\theta_\EPS)\,\d\xx
+\int_\varOmega
\chi_\zeta'(\theta_\EPS)\COND(\FFeps,\theta_\EPS)|\nabla\theta_\EPS|^2\,\d\xx
\\&\nonumber\  
=\!\int_\varOmega\!\bigg(\xi_\EPS\big(\FFeps,\theta_\EPS;\EE(\vv_\EPS),\nabla^2\vv_\EPS,\rr_\EPS\big)
\,\chi_\zeta(\theta_\EPS)
+\OMEGA(\FFeps,\mm_\EPS,\theta_\EPS)\chi_\zeta'(\theta_\EPS)\vv_\EPS{\cdot}\nabla\theta_\EPS
\\&\nonumber\quad 
+\big[{\mathscr X}_\zeta\big]_\Fe'(\FFeps,\mm_\EPS,\theta_\EPS){:}\pdt{\FFeps}
+\big[{\mathscr X}_\zeta\big]_\mm'(\FFeps,\mm_\EPS,\theta_\EPS){\cdot}\pdt{\mm_\EPS}
\\&\nonumber\quad
+\chi_\zeta(\theta_\EPS)\frac{\pi_\LAM(\FFeps)\COUPLING'_\mm(\FFeps,\mm_\EPS,\theta_\EPS)}
{(1{+}\EPS\theta_\EPS^{1/2})\det\FFeps}{\cdot}
\big(\rr_\EPS{-}{\rm skw}(\nabla\vv_\EPS)\mm_\EPS\big)
\\&\quad
+\chi_\zeta(\theta_\EPS)\frac{\pi_\LAM(\FFeps)\COUPLING'_{\Fe}(\FFeps,\mm_\EPS,\theta_\EPS)\FFeps^\top
{:}\ee(\vv_\EPS)}{(1{+}\EPS\theta_\EPS)\det\FFeps}
\bigg)\,\d\xx
+\!\int_\varGamma\!\Big(h_\EPS(\theta_\EPS){+}\frac{\nu_\flat|\vv_\EPS|^p}{2{+}\EPS|\vv_\EPS|^p}\Big)\chi_\zeta(\theta_\EPS)\,\d S\,.
\label{Euler-thermodynam3-test+++}\end{align}
We realize that $\chi_\zeta'(\theta)=\zeta/(1{+}\theta)^{1+\zeta}$ as used
already in \eq{8-**-+} and that ${\cal X}_\zeta(\FFeps,\mm_\EPS,\theta_\EPS)\ge
c_K\theta_\EPS$ with some $c_K$ for $\theta_\EPS\ge0$ due to
\eq{Euler-ass-adiab}; again $K$ is a compact subset of ${\rm GL}^+(d)$
related here with the already proved estimates \eq{est+Fes-eps}.
The convective term in \eq{Euler-thermodynam3-test+++} is a bit delicate.
For any $\delta>0$, it can be estimated by H\"older inequality as
\begin{align}\nonumber
\int_\varOmega\W_\EPS\chi_\zeta'(\theta_\EPS)\vv_\EPS{\cdot}\nabla\theta_\EPS\,\d\xx
&\le\frac1\delta\int_\varOmega\chi_\zeta'(\theta_\EPS)|\vv_\EPS|^2
\W_\EPS^2\,\d\xx
+\delta\int_\varOmega\chi_\zeta'(\theta_\EPS)|\nabla\theta_\EPS|^\TWO\,\d\xx
\\&=\frac1\delta\int_\varOmega\chi_\zeta'(\theta_\EPS)|\vv_\EPS|^2
\W_\EPS^2\,\d\xx+\delta I_{\zeta}^{(2)}(\theta_\EPS)\,.
\label{est-of-convectiv++e}\end{align}
Denoting by $0<\COND_0=\inf_{\FF,\theta}\COND(\FF,\theta)$, and using
\eq{Euler-thermodynam3-test+++} integrated over $I=[0,T]$,
we further estimate:
\begin{align}\nonumber
&
I_{\zeta}^{(2)}(\theta_\EPS)=
\frac1\zeta\int_0^T\!\!\!\int_\varOmega
\chi_\zeta'(\theta_\EPS)|\nabla\theta_\EPS|^\TWO\,\d\xx\d t
\le\frac1{\COND_0\zeta}\int_0^T\!\!\!\int_\varOmega\!\COND(\FFeps,\theta_\EPS)
\nabla\theta_\EPS{\cdot}\nabla\chi_\zeta(\theta_\EPS)\,\d\xx\d t 
\\&\nonumber
\le\frac1{\COND_0\zeta}\bigg(
\int_0^T\!\!\!\int_\varOmega\!\COND(\FFeps,\theta_\EPS)
\nabla\theta_\EPS{\cdot}\nabla\chi_\zeta(\theta_\EPS)
\,\d\xx\d t+\int_\varOmega\!{\cal X}_\zeta(\FFeps(T),\theta_\EPS(T))\,\d\xx\bigg)
\\&\nonumber
=\frac1{\COND_0\zeta}\bigg(\int_\varOmega\!{\cal X}_\zeta(\Fe_0,\theta_{0,\EPS})\,\d\xx
+\!\int_0^T\!\!\!\int_\varOmega\!\bigg(
\xi_\EPS\big(\FFeps,\theta_\EPS;\EE(\vv_\EPS),\nabla^2\vv_\EPS,\rr_\EPS\big)
\chi_\zeta(\theta_\EPS)
\\&\nonumber\ \ 
+\big[{\mathscr X}_\zeta\big]_\Fe'(\FFeps,\mm_\EPS,\theta_\EPS){:}\pdt{\FFeps}+
\big[{\mathscr X}_\zeta\big]_\mm'(\FFeps,\mm_\EPS,\theta_\EPS){\cdot}\pdt{\mm_\EPS}
+\OMEGA(\FFeps,\mm_\EPS,\theta_\EPS)\chi_\zeta'(\theta_\EPS)\vv_\EPS{\cdot}\nabla\theta_\EPS
\\&\nonumber\ \ 
+\chi_\zeta(\theta_\EPS)\frac{\pi_\LAM(\FFeps)\COUPLING'_\mm(\FFeps,\mm_\EPS,\theta_\EPS)}
{(1{+}\EPS\theta_\EPS^{1/2})\det\FFeps}{\cdot}
\big(\rr_\EPS{-}{\rm skw}(\nabla\vv_\EPS)\mm_\EPS\big)
\\&\nonumber\ \ 
+\frac{\pi_\LAM(\FFeps)\COUPLING'_{\Fe}(\FFeps,\mm_\EPS,\theta_\EPS)\FFeps^\top
{:}\ee(\vv_\EPS)}{(1{+}\EPS\theta_\EPS)\det\FFeps}\chi_\zeta(\theta_\EPS)
\bigg)\,\d\xx
+\!\int_0^T\!\!\!\int_\varGamma\!\Big(h_\EPS(\theta_\EPS)
+\frac{\nu_\flat|\vv_\EPS|^p}{2{+}\EPS|\vv_\EPS|^p}\Big)
\chi_\zeta(\theta_\EPS)\,\d S\d t\bigg)
\\[-.2em]&\nonumber
\!\stackrel{\eq{est-of-convectiv++e}}{\le}\!
\frac1{\COND_0\zeta}\bigg(\big\|{\cal X}_\zeta(\Fe_0,\theta_{0,\EPS})\big\|_{L^1(\varOmega)}\!
+\big\|\DIS|\EE(\vv_\EPS)|^p{+}\tau|\rr_\EPS|^2{+}\HC(\FFeps,\theta_\EPS)|\rr_\EPS|{+}\nu_2|\nabla^2\vv_\EPS|^p\big\|_{L^1(I\times\varOmega)}\!
\\&\quad\nonumber
+\int_0^T\!\!\big\|\big[{\mathscr X}_\zeta\big]_\Fe'(\FFeps,\mm_\EPS,\theta_\EPS)
\big\|_{L^{r'}(\varOmega;\R^{d\times d})}^{r'}\!
+\Big\|\pdt{\FFeps}\Big\|_{L^r(\varOmega;\R^{d\times d})}^r\d t
+\big\|h_{\max}\big\|_{L^1(I\times\varGamma)}\!
\\&\quad\nonumber
+\big\|\big[{\mathscr X}_\zeta\big]_\mm'(\FFeps,\mm_\EPS,\theta_\EPS)\big\|_{L^2(I\times\varOmega;\R^d)}^2\!+
\Big\|\pdt{\mm_\EPS}\Big\|_{L^2(I\times\varOmega;\R^d)}^2\!+\frac12\big\|\sqrt{\nu_\flat}\vv_\EPS\big\|_{L^p(I\times\varGamma;\R^d)}^p\!
\\&\quad\nonumber
+\Big\|
\pi_\LAM(\FFeps)\frac{\COUPLING'_{\Fe}(\FFeps,\mm_\EPS,\theta_\EPS)\FFeps^\top{:}\ee(\vv_\EPS)
{+}\COUPLING'_\mm(\FFeps,\mm_\EPS,\theta_\EPS){\cdot}\big(\rr_\EPS{-}{\rm skw}(\nabla\vv_\EPS)\mm_\EPS\big)}
{\det\FFeps}\Big\|_{L^1(I\times\varOmega)}\!\!
\\[-.2em]&
\quad
+\frac1\delta\|\vv_\EPS\|_{L^2(I;L^\infty(\Omega;\R^d))}^{\TWOprime}
\|\chi_\zeta'(\theta_\EPS)\OMEGA^{\TWOprime}(\FFeps,\mm_\EPS,\theta_\EPS)
\|_{L^\infty(I;L^1(\varOmega))}^{}\!\bigg)
+\frac\delta{\COND_0}I_{\zeta}^{(2)}(\theta_\EPS);
\label{+++}
\end{align}
noteworthy, we choose $\delta<\COND_0$, we can absorb the last
term in the left-hand side. Due to the assumption \eq{Euler-ass-adiab},
we can estimate the adiabatic rates
$\pi_\LAM(\FFeps)\COUPLING'_{\Fe}(\FFeps,\mm_\EPS,\theta_\EPS)\FFeps^\top{:}\ee(\vv_\EPS)/\!\det\FFeps$ and 
$\pi_\LAM(\FFeps)\COUPLING'_\mm(\FFeps,\mm_\EPS,\theta_\EPS){\cdot}\big(\rr_\EPS{-}{\rm skw}(\nabla\vv_\EPS)\mm_\EPS\big)/\!\det\FFeps$
in \eq{+++}, cf.\ \eq{Euler-est-of-rhs+}. We also use the estimates
\eq{est-eps} and \eq{est-time-derivaitves} and the assumption
\eq{Euler-ass-h} relying on the already proved non-negativity
of temperature.
By the qualification \eq{Euler-ass-primitive-c}, we have 
$|[{\mathscr X}_\zeta]_\Fe'(\FF,\mm,\theta)|\le C(1{+}\theta)$.
This allows for estimation 
\begin{align}\nonumber
&\big\|\big[{\mathscr X}_\zeta\big]_\Fe'(\FFeps,\mm_\EPS,\theta_\EPS)
\big\|_{L^{r'}(\varOmega;\R^{d\times d})}^{r'}
\le C^{r'}\big\|1+\theta_\EPS\big\|_{L^{r'}(\varOmega)}^{r'}
\\&\hspace{9em}\le C_4+C_4\|\theta_\EPS\|_{L^1(\varOmega)}^{r'(1-{\EXP^*}'/r)}\big(\|\theta_\EPS\|_{L^1(\varOmega)}^{}
\!{+}\|\nabla\theta_\EPS\|_{L^\EXP(\varOmega;\R^d)}\big)^{{\EXP^*}'/(r-1)}\!\,,
\label{est-adiabatic}\end{align}
where we use the Gagliardo-Nirenberg inequality to interpolate
$L^{r'}(\varOmega)$ between $L^1(\varOmega)$ and $W^{1,\EXP}(\varOmega)$.
Similarly, since $|[{\mathscr X}_\zeta]_\mm'(\FF,\mm,\theta)|^2\le C(1{+}\theta)$
which is again ensured by the qualification \eq{Euler-ass-primitive-c}, we can estimate
$\|[{\mathscr X}_\zeta]_\mm'(\FFeps,\mm_\EPS,\theta_\EPS)\|_{L^2(I\times\varOmega;\R^d)}^2\le
C\big\|1{+}\theta_\EPS\big\|_{L^1(I\times\varOmega)}$, which is bounded due to
\eq{basic-est-of-theta-eps}. The penultimate term in \eq{+++} is a-priori bounded
independently of $\EPS$ for $\zeta>0$ fixed because, as
$\OMEGA(\FF,\mm,\theta)=\mathscr{O}(\theta)$ due to \eq{Euler-ass-primitive-c}
and due to $\chi_\zeta'(\theta)=\mathscr{O}(1/\theta)$ uniformly for $\zeta>0$,
so that we have
$[\chi_\zeta'(\cdot)\OMEGA^\TWO(\FF,\mm,\cdot)](\theta)=\mathscr{O}(\theta)$.
Thus the estimate \eq{basic-est-of-theta-eps} guarantees
$\chi_\zeta'(\theta_\EPS)\OMEGA^{\TWOprime}(\FFeps,\mm_\EPS,\theta_\EPS)$ bounded
in $L^\infty(I;L^1(\Omega))$ while $|\vv_\EPS|^{\TWOprime}$
is surely bounded in $L^1(I;L^\infty(\Omega))$, cf.\ \eq{est+v}.

In view of \eq{est-adiabatic}, one can summarize \eq{+++} as
$I_{\zeta}^{(2)}(\theta_\EPS)\le C(1+\|\nabla\theta_\EPS\|_{L^\EXP(\varOmega;\R^d)}\big)^{{\EXP^*}'/(r-1)}$.
Combining it with \eq{8-***}, one obtain the inequality as
\begin{align}\|\nabla\theta_\varepsilon\|_{L^\EXP(I\times\Omega;\R^d)}^{}
\le C\big(1+\|\nabla\theta_\varepsilon\|_{L^\EXP(I\times\Omega;\R^d)}^{1-\EXP/2+{\EXP^*}'/(2r-2)}\big)\,.
\label{final-est-nabla-theta}
\end{align}
Reminding $\sigma:=(\TWO{-}\EXP)/(1{+}\zeta)$ from \eq{8-cond} with $\zeta>0$
arbitrarily small, one gets after some algebra, the condition $\EXP<(d{+}2)/(d{+}1)$.
Obviously, for $r$ big enough (in particular if $r>d$ as assumed),
the exponent in the left-hand side of \eq{final-est-nabla-theta} is higher than
the exponent in the right-hand side, which gives a bound for
$\nabla\theta_\varepsilon$ in $L^\EXP(I\times\Omega;\R^d)$. Altogether, we proved
\begin{subequations}\label{est-W-eps}\begin{align}
&\|\theta_\EPS\|_{L^\infty(I;L^1(\varOmega))\,\cap\,L^\EXP(I;W^{1,\EXP}(\varOmega))}^{}\le C_\EXP
\ \ \text{ with }\ 1\le\EXP<\frac{d{+}2}{d{+}1}\,.
\intertext{Next, we again exploit the calculus \eq{w=...} now omitting the
index $k$, with $\nabla\FF_{\EPS}$ bounded in $L^\infty(I;L^r(\varOmega;\R^{d\times d\times d}))$
and $\nabla\mm_{\EPS}$ bounded in $L^\infty(I;L^2(\varOmega;\R^{d\times d}))$ and
relying on the assumption \eq{Euler-ass-primitive-c}, we have also the
bound on $\nabla\W_{\EPS}$ in $L^\EXP(I;L^{\EXP^*d/(\EXP^*+d)}(\varOmega;\R^d))$, so
that}
&\|\W_\EPS\|_{L^\infty(I;L^1(\varOmega))\,\cap\,L^\EXP(I;W^{1,\EXP^*d/(\EXP^*+d)}(\varOmega))}^{}\le C_\EXP\,.
\end{align}\end{subequations}

\medskip\noindent{\it Step 8: Limit passage for $\EPS\to0$}.
We use the Banach selection principle as in Step~4, now also taking 
\eq{est-eps} and \eq{est-W-eps} into account instead of
the estimates \eq{Euler-quasistatic-est1} and \eq{Euler-quasistatic-est2}.
For some subsequence and some $(\varrho,\vv,\FF,\mm,\theta)$, we now have
\begin{subequations}\label{Euler-weak++}
\begin{align}\nonumber
&\varrho_\EPS\to\varrho&&\hspace*{-11em}
\text{weakly* in $\ L^\infty(I;W^{1,r}(\varOmega))\,\cap\,W^{1,p}(I;L^r(\varOmega))$}
\\\label{Euler-weak-rho}
&&&\hspace*{-2em}\text{and strongly in  $C(I{\times}\barOmega)$}\,,\\
&\vv_\EPS\to\vv&&\hspace*{-11em}\text{weakly* in $\
L^\infty(I;L^2(\varOmega;\R^d))\cap
L^2(I;W^{2,p}(\varOmega;\R^d))$,}\!\!&&
\label{Euler-weak-v}\\\nonumber
&\FFeps\to\FF\!\!\!&&\hspace*{-11em}\text{weakly* in $\
L^\infty(I;W^{1,r}(\varOmega;\R^{d\times d}))\,\cap\,
W^{1,p}(I;L^r(\varOmega;\R^{d\times d}))$}\!\!
\\&&&\hspace*{-2em}\text{and strongly in
$C(I{\times}\barOmega;\R^{d\times d})$, }
\label{Euler-weak-F}
\\
&\mm_\EPS\to\mm\!\!\!&&\hspace*{-11em}\text{weakly* in $L^\infty(I;H^1(\varOmega;\R^d))\,\cap\,H^1(I;L^2(\varOmega;\R^d))$},
\\&
\theta_\EPS\to\theta\!\!\!&&\hspace*{-11em}\text{weakly* in $\ L^\EXP(I;W^{1,\EXP}(\varOmega)),\
1\le\EXP<(d{+2})/(d{+}1)$.}
\intertext{Like \eq{w-conv}, by the Aubin-Lions theorem, we now have}
&\W_\EPS\to\W\!\!\!&&\hspace*{-11em}\text{strongly in $L^c(I{\times}\varOmega),\ \
1\le c<1{+}2/d$,}
\intertext{and then, using again continuity of
$(\FF,\mm,\W)\mapsto[\OMEGA(\FF,\mm,\cdot)]^{-1}(\W)$ as in \eq{z-conv}, we  have}
&\theta_\EPS\to\theta=[\OMEGA(\FF,\mm,\cdot)]^{-1}(\W)\!\!\!&&
\hspace*{-8em}\text{strongly in $\ L^c(I{\times}\varOmega),\ \ 1\le c<1{+}2/d$.}
\intertext{By the continuity of $\upvarphi_\FF'$, $\COUPLING_\FF'$, $\det$,
and $\COND$, we have also}
&\label{m-strongly+}
\COND(\FFeps,\theta_\EPS)\to\COND(\FF,\theta)
&&\hspace*{-8em}\text{strongly in $L^c(I{\times}\varOmega)$ for any $1\le c<\infty$, and}
\\\label{Euler-weak-stress}
&\TT_{\LAM,\EPS}\to\TT_\LAM=
\frac{[\pi_\LAM\upvarphi]_\FF'(\FF,\mm){+}\pi_\LAM(\FF)\COUPLING_\FF'(\FF,\mm,\theta)\!}{{\det}
\FF}\FF^\top
\hspace*{-0em}&&\hspace*{0em}\text{strongly in $L^1(I{\times}\varOmega;\R_{\rm sym}^{d\times d})$.}
\end{align}\end{subequations}

The momentum equation \eq{Euler-thermo-reg1} (still regularized by $\varepsilon$) 
is to be treated like in Step~4. Here we exploit 
the information about $\pdt{}(\varrho_\EPS\vv_\EPS)$ in 
$L^{q'}(I;W^{1,q}(\varOmega;\R^d)^*)+L^{p'}(I;W^{2,p}(\varOmega;\R^d)^*)$
obtained like in \eq{est-of-DT-rho.v}; here we used also \eq{est+v}.
By the Aubin-Lions compact-embedding theorem, we then obtain 
\begin{align}\label{rho-v-conv+}
&\varrho_\EPS\,\vv_\EPS\to\varrho\vv
&&\hspace*{-1em}\text{strongly in }L^s(I{\times}\varOmega;\R^d)\ \ \text{ with
$s$ from \eq{est+v}}\,.
\end{align}
In fact, the argumentation \eq{strong-hyper+} now with $C_r$ instead of
$C_{r,\EPS}$ is to be slightly modified by using
$(\varrho_\EPS,\vv_\EPS,\TT_{\LAM,\EPS}(\FFeps,\vv_\EPS),\theta_\EPS)$
in place of $(\varrho_{\EPS k},\vv_{\EPS k},\TT_{\LAM,\EPS}(\FFepsk,\vv_{\EPS k}),\theta_{\EPS k})$
and with $\widetilde\vv_k$ replaced by $\vv$. Also,
$\int_0^T\!\int_\varOmega\pdt{}(\varrho_{\EPS k}\vv_{\EPS k}){\cdot}\widetilde\vv_k\,\d\xx\d t$
is to be replaced by the duality
$\langle\pdt{}(\varrho_\EPS\vv_\EPS),\vv\rangle$
with $\langle\cdot,\cdot\rangle$ denoting here the duality between
$L^{q'}(I;W^{1,q}(\varOmega;\R^d)^*)+L^{p'}(I;W^{2,p}(\varOmega;\R^d)^*)$
and $L^q(I;W^{1,q}(\varOmega;\R^d))\cap L^p(I;W^{2,p}(\varOmega;\R^d))$.

Limit passage in the heat equation \eq{Euler-thermo-reg3} is then simple. 
Altogether, we proved that $(\varrho,\vv,\FF,\mm,\theta)$ solves in the weak
sense the problem \eq{Euler-thermo-reg}--\eq{Euler-thermo-reg-BC-IC}
with $\EPS=0$ and with $\TT_\LAM$ from \eq{Euler-weak-stress} in place of
$\TT_{\LAM,\EPS}$.

\medskip\noindent{\it Step 9: the original problem}.
Let us note that  the limit $\Fe$ lives in 
$L^\infty(I;W^{1,r}(\varOmega;\R^{d\times d}))\,\cap\,
W^{1,p}(I;L^\infty(\varOmega;\R^{d\times d}))$, cf.\ (\ref{est-eps}b,f),
and this space  is embedded  into 
$C(I{\times}\barOmega;\R^{d\times d})$ if $r>d$. Therefore $\Fe$ and its
determinant evolve continuously in time, being valued respectively in
$C(\barOmega;\R^{d\times d})$ and $C(\barOmega)$.
Let us recall that the initial condition $\FF_0$ complies with the bounds
\eq{Euler-quasistatic-est-formal4} and we used this $\FF_0$
also for the $\LAM$-regularized system.
Therefore $\Fe$ satisfies these bounds not only at $t=0$ but also at least
for small times. Yet, in view of the choice \eq{Euler-quasistatic-est-formal4}
of $\LAM$, this means that the $\LAM$-regularization is nonactive
and  $(\varrho,\vv,\Fe,\mm,\theta)$ solves, at least for a small time, the
original nonregularized problem
\eq{Euler-thermo-reg}--\eq{Euler-thermo-reg-BC-IC}
for which the a~priori $L^\infty$-bounds \eq{est+} hold.
By the continuation argument, we may see that the $\LAM$-regularization 
remains therefore inactive within the whole evolution of
$(\varrho,\vv,\Fe,\mm,\theta)$ on the whole time interval $I$.

\medskip\noindent{\it Step 10: energy balances}.
It is now important that the tests and then all the subsequent calculations
leading to the energy balances \eq{energy+} and
\eq{thermodynamic-Euler-engr} integrated over a current
time interval $[0,t]$ are really legitimate.

In the calculus \eq{Euler-large-thermo}, we rely on that 
$[\upvarphi(\FF,\mm)/\!\det\FF]_\FF'\in L^\infty(I{\times}\varOmega;\R^{d\times d})$
is surely in duality with
$\pdt{}\FF\in L^{p}(I;L^r(\varOmega;\R^{d\times d}))$ and
$(\vv{\cdot}\nabla)\FF\in L^s(I;L^r(\varOmega;\R^{d\times d}))$
with $s$ from \eq{est+v}. 
Moreover, $\pdt{}(\varrho\vv)\in
L^{1}(I;L^2(\varOmega;\R^d)^*)+L^{p'}(I;W^{2,p}(\varOmega;\R^d)^*)$
is in duality with
$\vv\in L^\infty(I;L^2(\varOmega;\R^d))\cap L^p(I;W^{2,p}(\varOmega;\R^d))$,
as used in \eq{calculus-convective-in-F}.
Further, the calculus \eq{calculus-convective-in-F} relies on
that $\pdt{}\varrho$ and ${\rm div}(\varrho\vv)=
\vv{\cdot}\nabla\varrho+\varrho\,{\rm div}\,\vv$ live in $L^s(I;L^{rs/(r+s)}(\varOmega))$
and thus are surely in duality with $|\vv|^2\in
L^{s/2}(I;L^\infty(\varOmega))$ with $3\le s<p(pd{+}4p{-}2d)/(4p{-}2d)$, cf.\
\eq{est+v}.
Eventually, since $\nabla^2\vv\in L^{p}(I{\times}\varOmega;\R^{d\times d\times d})$,
we have ${\rm div}^2(\nu_2|\nabla^2\vv|^{p-2}\nabla^2\vv)\in
L^{p'}(I;W^{2,p}(\varOmega;\R^d)^*)$ in
duality with $\vv$. Also ${\rm div}(\nu_1|e(\vv)|^{p-2}\ee(\vv))\in
L^{p'}(I;W^{1,p}(\varOmega;\R^d)^*)$ is in duality with $\vv$ due to the
growth condition \eq{Euler-ass-xi}.
Altogether, the calculations
\eq{Euler-large-thermo}--\eq{calculus-convective-in-F} are legitimate.

Recalling in particular
${\rm div}(\upkappa(\FF)\nabla\mm/\!\det\FF)\in L^2(I{\times}\varOmega;\R^d)$,
we can see that also the calculations \eq{test-damage}--\eq{calculus-PM4}
are legitimate. This ends the proof of Theorem~\ref{prop-Euler}.


\begin{remark}[{\sl The classical solutions}]\upshape
In fact, having $\nabla\mm$ estimated from the exchange energy,
\eq{Euler-thermo-reg5} holds even in the sense of
$L^2(I{\times}\varOmega;\R^d)$, not only in the weak sense, similarly as
(\ref{Euler-weak-Galerkin}a,b). Actually, since
${\rm div}(\upkappa(\FF)\nabla\mm/\!\det\FF)\in L^2(I{\times}\varOmega;\R^d)$,
the inequality \eq{Euler2-weak} can be formulated as a classical inclusion
\eq{Euler-thermodynam4} and \eq{Euler-thermo-reg4} a.e.\ in the sense of
$L^2(I{\times}\varOmega;\R^d)$.
\end{remark}

\begin{remark}[{\sl Importance of the exchange energy}]\upshape
At large magnets, in contrast to micromagnetism, the influence of the
exchange energy is small and, for very large magnetic continua,
eventually negligible, cf.\  \cite{DeSi93EMLF}. Yet, this energy
controls $\nabla\mm$, which ensures for ``compactness'' and strong
convergence in $\mm$ and its complete deletion would be analytically problematic
in particular because nonconvexity of $\uppsi(\FF,\cdot,\theta)$
in ferromagnetic phase .
\end{remark}

\bigskip

{\small

\baselineskip=12pt

\noindent{\it Acknowledgments.}
The author is deeply thankful to Giuseppe Tomassetti for
many inspiring discussions and comments to the manuscript.
A support from 
the Ministry of Education of the
Czech Republic project CZ.02.1.01/0.0/0.0/15-003/0000493
and the CSF/DFG project GA22-00863K,
and from the institutional support RVO:61388998 (\v CR)  is 
acknowledged.

} 


\end{document}